%% file: main.tex
\documentclass[11pt]{article}
\usepackage[top=1in,bottom=1in,left=1in,right=1in]{geometry}

\RequirePackage{amsthm,amsmath,amsfonts,amssymb}

\usepackage{thmtools}
\usepackage{thm-restate}

\usepackage{graphicx, subcaption}
\usepackage{multirow}
\usepackage{array}
\usepackage{booktabs,adjustbox}

\usepackage{algorithm}
\usepackage{algorithmic}
\usepackage[algo2e]{algorithm2e}

\usepackage{bbm}
\usepackage[normalem]{ulem}

\usepackage{hyperref}
\usepackage{cleveref}
\usepackage{dsfont}
\usepackage{enumitem}
\usepackage{stmaryrd}

\makeatletter
\newcommand{\myitem}[1]{%
\item[#1]\protected@edef\@currentlabel{#1}%
}
\makeatother

\newcommand{\argmin}{\operatornamewithlimits{argmin}}
\newcommand{\argmax}{\operatornamewithlimits{argmax}}

\input{enviroments}

\input{macroses}

\newlength\mgwidth 
\newlength\mgheight

\begin{document}

\title{Minimax Boundary Estimation and Estimation with Boundary}

\author{
	Eddie Aamari%
	\footnote{CNRS \& U. Paris Cité \& Sorbonne U., Paris, France (\url{https://perso.lpsm.paris/\~aamari/})}%
	\and
	Catherine Aaron%
	\footnote{CNRS \& U. Clermont Auvergne, Clermont-Ferrand, France (\url{https://lmbp.uca.fr/\~aaron/})}%
	\and
	Cl\'ement Levrard%
	\footnote{CNRS \& U. Paris Cité \& Sorbonne U., Paris, France (\url{http://www.normalesup.org/\~levrard/})}%
}

\date{}
\maketitle

\begin{abstract}
	\input{abstract}
\end{abstract}

	\input{introduction}

	\input{framework}

	\input{results}

	\input{conclusion}
	
	\input{proof-outline}

 	\input{simulations}

\newpage
\appendix

	\input{geom_results}

	\input{lemmas-manifolds-with-boundary}

	\input{proofs-main-detection}

    \input{tangent-space-estimation}

	\input{linear-patches}

	\input{lower-bounds-proofs}

	\input{lower-bounds-main-tools}
	
\bibliographystyle{plain}
    \bibliography{biblio}

\end{document}

%% file: enviroments.tex
\theoremstyle{plain}

\declaretheorem[name=Theorem,numberwithin=section]{theorem}

\declaretheorem[name=Lemma,numberwithin=section,sibling=theorem]{lemma}

\declaretheorem[name=Proposition,numberwithin=section,sibling=theorem]{proposition}

\declaretheorem[name=Corollary,numberwithin=section,sibling=theorem]{corollary}

\newtheorem*{theorem*}{Theorem}
\newtheorem*{proposition*}{Proposition}
\newtheorem*{lemma*}{Lemma}
\newtheorem*{corollary*}{Corollary}
\newtheorem*{definition*}{Definition}
\newtheorem*{example*}{Example}

\theoremstyle{definition}
\newtheorem{definition}[theorem]{Definition}
\newtheorem{remark}[theorem]{Remark}

%% file: macroses.tex
\newcommand{\set}[1]{\left\{{#1}\right\}}
\DeclareMathOperator{\sign}{sign}

\newcommand{\R}{\mathbb{R}}

\newcommand{\ceil}[1]{\left\lceil#1\right\rceil}

\DeclareMathOperator{\vol}{vol}
\newcommand{\Haus}{\mathcal{H}}

\newcommand{\supp}{\mathrm{Supp}}
\newcommand{\E}{\mathbb{E}}

\newcommand{\indicator}[1]{\mathbbm{1}_{#1}}

\newcommand{\TV}{\mathop{\mathrm{TV}}}

\newcommand{\inner}[2]{\left\langle {#1}, {#2} \right\rangle}
\newcommand{\norm}[1]{\left\Vert #1 \right\Vert }
\newcommand{\B}{\mathrm{B}}
\newcommand{\Sphere}{\mathcal{S}}
\newcommand{\interior}[1]{\mathring{#1}} 
\newcommand{\Bopen}{\interior{\B}}

\newcommand{\hull}{\mathrm{Hull}}

\newcommand{\normop}[1]{\left\Vert #1 \right\Vert_{\mathrm{op}} }

\newcommand{\Gr}[2]{\mathbb{G}^{#1,#2}}
\newcommand{\Span}{\mathrm{span}}
\newcommand{\Operp}{\overset{\perp}{\oplus}}

\newcommand{\dd}{\mathrm{d}}
\newcommand{\dHaus}{\mathrm{d_H}}

\newcommand{\Length}{\mathrm{Length}}
\newcommand{\Med}{\mathrm{Med}}
\newcommand{\vor}{\mathrm{Vor}}
\DeclareMathOperator{\Int}{Int}

\newcommand{\eps}[0]{\varepsilon}

\newcommand{\X}{\mathbb{X}}

\newcommand{\BO}{\hyperref[alg:boundary-labelling]{\texttt{Boundary Structure}}\xspace}

%% file: abstract.tex
We derive non-asymptotic minimax bounds for the Hausdorff estimation of $d$-dimensional submanifolds $M \subset \mathbb{R}^D$ with (possibly) non-empty boundary $\partial M$.
The model reunites and extends the most prevalent $\mathcal{C}^2$-type set estimation models: manifolds without boundary, and full-dimensional domains.
We consider both the estimation of the manifold $M$ itself and that of its boundary $\partial M$ if non-empty.
Given $n$ samples, the minimax rates are of order $O\bigl((\log n/n)^{2/d}\bigr)$ if $\partial M = \emptyset$ and $O\bigl((\log n/n)^{2/(d+1)}\bigr)$ if
$\partial M \neq \emptyset$, up to logarithmic factors.
In the process, we develop a Voronoi-based procedure that allows to identify enough points $O\bigl((\log n/n)^{2/(d+1)}\bigr)$-close to $\partial M$ for reconstructing it.
Explicit constant derivations are given, showing that these rates do not depend on the ambient dimension $D \gg d$.

%% file: introduction.tex
\section{Introduction}
	\label{sec:introduction}

Topological data analysis and geometric inference techniques have significantly grown in importance in the high-dimensional statistics area, both in its theoretical and practical aspects~\cite{Wasserman18,Chazal17}.
Unlike Lasso-type methods~\cite{Hastie09} which strongly rely on a specific coordinate system, geometric inference techniques naturally yield features that are invariant through rigid transformations of the ambient space.

A central problem in this field is manifold estimation~\cite{Belkin06,Genovese12a,Genovese12b,Aamari19b,Puchkin19,Fefferman19}.
Assuming that data $\X_n = \set{X_1,\ldots,X_n}$ originate from some unknown distribution $P$ on $\R^D$, these works study the estimation of its support $M = \supp(P) \subset \R^D$, assumed to be a submanifold of dimension $d \ll D$.
This provides a non-linear dimension reduction, that can allow to mitigate the curse of dimensionality, and helps for data visualization~\cite{Lee07}.
Manifold estimation is also of crucial importance for inferring other geometric features of $M$, as it appears as a critical intermediate step in a growing series of plugin strategies. See for instance~\cite{Chazal13} for persistent homology,~\cite{Berenfeld20} for the reach, or~\cite{Divol21} for density estimation.

\subsection{Support estimation}

\subsubsection*{Overview}
\label{par:overview}

So far, the statistical study of support estimation in Hausdorff distance has been carried out within two somehow orthogonal settings:
Full dimensional domains $\dim(M)=D$ on one hand --- which necessarily have non-empty boundary $\partial M\neq \emptyset$ ---, and low-dimensional submanifolds $\dim(M)=d <D$ without boundary $\partial M = \emptyset$ on the other hand. 
More precisely:

\begin{enumerate}[label=(\roman*)
,leftmargin=2em]
\item \label{item:intro-boundary-full}
Assuming that $M = \supp(P) \subset \R^D$ is \emph{full-dimensional} $\dim(M) = D$ (i.e. roughly everywhere of non-empty interior) and that $P$ has enough mass in every neighborhood
of its support,~\cite{Dumbgen96,Cuevas04} derive error bounds of order $(\log n /n)^{1/D}$. 
Here, a rate-optimal estimator simply consists of the sample set $\hat{M} = \X_n$ itself.
Even under the additional geometric restriction of $M$ being convex, this rate is still the best possible, due to the possible outward corners a convex set may contain.
Beyond convexity, faster rates can actually be attained with additional smoothness constraints: if the (topological) boundary $\bar{\partial} M$ of the convex $M$ is $\mathcal{C}^2$-smooth,~\cite{Dumbgen96} derives a convergence rate of order $(\log n/n)^{2/(D+1)}$ by considering the convex hull $\hat{M} = \hull(\X_n)$, which also allows to estimate $\bar{\partial} M$ with $\bar{\partial} \hat{M}$ at the same rate.
This phenomenon was also exhibited by~\cite{Tsybakov95,Casal07,Aaron16b} in similar convexity-type settings.
Let us also mention the recent work of \cite{Calder22}, which proposed a computationally efficient (yet unfortunately rate-suboptimal) boundary labelling method.

Note that despite a nearly quadratic gain in the rate for smooth cases, this framework still remains a hopeless scenario for high dimensional datasets, as it heavily suffers from the curse of dimensionality, both statistically and computationally.
\\
\noindent
This paper extends these results for $M$ possibly of lower dimension $d\ll D$ and curved.

\item \label{item:intro-boundaryless-low}
To overcome the curse of dimensionality, assuming that $M = \supp(P) \subset \R^D$ is a $\mathcal{C}^2$ submanifold of dimension $\dim(M) = d < D$ \emph{with empty (differential) boundary} $\partial M = \emptyset$,~\cite{Genovese12a,Kim2015} show that the minimax rate of estimation of $M$ is of order $(\log n /n)^{2/d}$.
The estimator of~\cite{Genovese12a} being intractable in practice,~\cite{Aamari18} later proposed an optimal algorithm that outputs a triangulation of the data points which is computable in polynomial time. 
Using local polynomials, faster estimation rates of order $(\log n / n)^{k/d}$ were also shown to be achievable over $\mathcal{C}^k$-smooth submanifolds~\cite{Aamari19b}.
Although insensitive to the ambient dimension $D \gg d$, these results highly rely on the fact that $\partial M = \emptyset$.
\\
\noindent
This paper extends these results for $M$ possibly with non-empty boundary.
\end{enumerate}

\subsubsection*{Background}

\begin{figure}[h]
\centering
\includegraphics[width=0.4\textwidth]{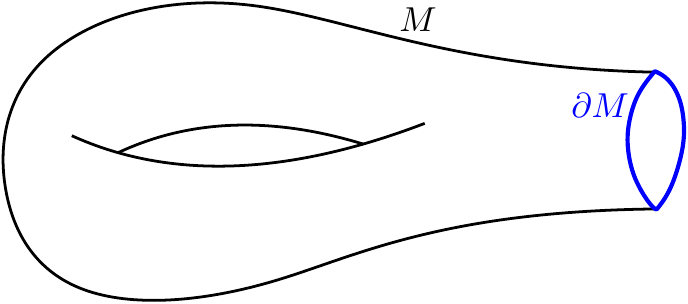}
\caption{A surface $M$ with non-empty boundary $\partial M$. Note that $\dim (\partial M) = \dim (M) -1$, so that sample points from some roughly uniform distribution $P$ on $M$ almost surely never belong to $\partial M$.
However, points close to $\partial M$ should be processed differently in the analysis of such a sample, since they have an unbalanced neighborhood: they may cause boundary effects.}
\label{fig:boundary_example}
\end{figure}

By definition, a submanifold $M \subset \R^D$ of dimension $d$ is a smooth subspace that can be parametrized locally by $\R^d$. Hence, neighborhoods of points in $M$ 
all look like $d$-dimensional balls. In contrast, a manifold \emph{with boundary} is a smooth space that can be parametrized locally by $\R^d$ or by $\R^{d-1} \times \R_+$. 
If not empty, the boundary of $M$, denoted by $\partial M$, is the set of points nearby which $M$ can only be parametrized by $\R^{d-1} \times \R_+$.
Informally, the class of manifolds with boundary allows to take into account the possible ``rims'' a surface may contain (see \Cref{fig:boundary_example}).

As mentioned above, most of the existing manifold estimation techniques require that $\partial M = \emptyset$, which is very restrictive in view of real data \cite{Wasserman18}.
When the empty boundary condition is dropped, the location of $\partial M$ is often assumed to be known via an \emph{oracle}, able to correctly label points that lie close to the boundary \cite{Rineau08}.
Prior to the present paper, a theoretically grounded construction of such an oracle given unlabeled data was not known, since the optimal detection and estimation rates of $\partial M$ in arbitrary dimension had not been studied.
This is mainly due to the technicalities that the presence of a boundary usually gives rise to. 
For instance, the restricted Delaunay triangulation to a surface with boundary may not even be homeomorphic to the surface~\cite{Dey09b}. 
Hence, Delaunay-based reconstructions are not good candidates to handle boundary, which contrasts sharply with the boundaryless case~\cite{Boissonnat14,Aamari18}.
Despite these barriers, a few interesting works on boundary inference can be found in the literature.

For surfaces in space ($d=2$, $D = 3$), the so-called peeling algorithm consists in pruning an ambient triangulation (the $\alpha$-shape of the point cloud) to handle boundary~\cite{Dey09b}.
This method leverages boundary triangles being flatter than inner triangles. 
Unfortunately, such a method is limited to low dimensions, for the same instability problems described in~\cite{Boissonnat09}.

On the other hand, in full dimension ($d=D$), \cite{Cuevas04} proposed a plugin estimator based on an estimator of $M$ itself: under technical
constraints, if $\hat{M}$ approximates $M$, then $\partial \hat{M}$ approximates $\partial M$.
Such an plugin strategy provides a wide range of very general consistent boundary estimators: see for instance~\cite{Casal07,Aaron16b} for convergence rates under additional assumptions.
Note that naturally, such an approach is very costly --- as acknowledged by the authors themselves ---, and does not generalize easily to non-linear low-dimensional cases.

More recently, \cite{Aaron16a} designed an asymptotic boundary detection scheme based on local barycenter displacements: if a point $x \in M$ is close to $\partial M$, then the ball $M \cap \B(x,r)$ around $x$ will not be balanced, and its barycenter would shift away from $\partial M$. 
This naturally yields a criterion to decide whether $x$ belongs to $\partial M$ or not.
Unfortunately, this method requires the sampling density $f$ over $M$ to be Lipschitz and fails otherwise, as discontinuities of $f$ far from $\partial M$ may create artificial local barycenter  shifts, and hence false positives.
Let us also mention that this local barycenter shift has also been used in the context of density estimation on a manifold with boundary: \cite{Berry17} proposed a method for estimating the distance and direction of the boundary in order to correct the extra bias of a kernel density estimator near $\partial M$.

\subsection{Contribution}

This paper studies the minimax rates of estimation of $d$-dimensional $\mathcal{C}^2$-submanifolds $M \subset \R^D$ with possible $\mathcal{C}^2$ boundary $\partial M$ (\Cref{def:manifold_with_boundary}), and the estimation of the boundary itself if not empty.
As now standard in the literature, the loss is given by the Hausdorff distance $\dHaus$ (a sup-norm between sets, see \Cref{def:hausdorff_distance}), and $\mathcal{C}^2$ regularity of sets is measured through their reach $\tau_M,\tau_{\partial M}>0$ (a generalized convexity parameter, see \Cref{def:reach}).

Informally, we extend the known full-dimensional $\mathcal{C}^2$ support estimation rates to the case of low-dimensional curved $M$ with $\mathcal{C}^2$ boundary.
Indeed, if $M$ is contained in a $d$-dimensional affine subspace of $\R^D$ and has a $\mathcal{C}^2$ boundary, its estimation boils down to the full dimensional case (\Cref{par:overview}~\ref{item:intro-boundary-full}), and can be done with rate $(\log n / n)^{2/(d+1)}$ 
\cite{Casal07,Aaron16b}.
The present article proves that even if $M$ is curved, the same rate drives the estimation hardness of $M$ and $\partial M$. In addition, the estimator adapts automatically to the possible emptiness of $\partial M$, in which case $M$ can be estimated at rate $(\log n / n)^{2/d}$ (see \Cref{par:overview}~\ref{item:intro-boundaryless-low}).
More precisely, we show that for $n$ large enough independent of $D$,
\begin{align*}
\tag{\Cref{thm:boundary-main-upper-bound,thm:boundary-main-lower-bound}
}
\inf_{\hat{B}_n}
\sup_{
	\substack{
		\partial M \neq \emptyset
		\\
		\tau_{M} \geq \tau_{\min} 
		\\
		\tau_{\partial M} \geq \tau_{\partial, \min}
	}
}
\vartheta_n(\partial M)^{-1}
\E\,\big[\dHaus\bigl(\partial M,\hat{B}_n \bigr)\bigr]
=
\tilde{\Theta}(1)
,
\end{align*}
\vspace{-2em}
\begin{align*}\tag{\Cref{thm:manifold-main-upper-bound,thm:manifold-main-lower-bound}
}
\inf_{\hat{M}_n}
\sup_{
	\substack{
		\tau_{M} \geq \tau_{\min} 
		\\
		\tau_{\partial M} \geq \tau_{\partial, \min}
	}
}
\vartheta_n(\partial M)^{-1}
\E\,\big[\dHaus\bigl(M,\hat{M}_n \bigr)\bigr]
=
\tilde{\Theta}(1)
,
\end{align*}
where
$\hat{B}_n$ and $\hat{M}_n$ range among all the possible estimators based on $n$ samples, and
\begin{align*}
\vartheta_n(\partial M) \asymp 
\begin{cases}
(1 / n)^{2/(d+1)}
&
\text{if } \partial M \neq \emptyset,
\\
(1 / n)^{2/d}
&
\text{if } \partial M = \emptyset.
\end{cases}
\end{align*}%
These rates, given up to $\log n$ factors through $\tilde{\Theta}(1)$, do not depend on $D$.

\subsection{Outline}
We first describe the geometric framework and statistical setting we consider (\Cref{sec:framework}). Then, we state the main boundary detection and estimation results (\Cref{sec:main-results}) and discuss them (\Cref{sec:conclusion}). 
We present the principal steps of the proofs (\Cref{sec:proof-outline}).
Finally, we discuss complexity, heuristics for parameter selection, and provide illustrations of the method on synthetic data (\Cref{sec:experiments}).
For space constraints, the minor intermediate lemmas and most technical parts of the proofs are deferred to the Appendix.

%% file: framework.tex
\section{Framework}
	\label{sec:framework}
	
Throughout, $D \geq 1$ is referred to as the ambient dimension and $\mathbb{R}^D$ is endowed with the Euclidean inner product $\inner{\cdot}{\cdot}$ and the associated norm $\norm{\cdot}$. The closed Euclidean ball of center $x$ and radius $r$ is denoted by $\B(x,r)$, and its open counterpart by $\Bopen(x,r)$.
Given a linear subspace $T \subset \R^D$, we also write $\B_T(0,r) := T \cap \B(0,r)$ for the $r$-ball of $T$ centered at $0 \in T$.
\subsection{Geometric setting}

\subsubsection{Submanifolds with boundary}
By definition, the $d$-dimensional submanifolds $M \subset \R^D$ with boundary are the subsets of $\R^D$ that can locally be parametrized either
by the Euclidean space $\R^d$, or the half-space $\R^{d-1} \times \R_{+}$ \cite[Chapter~2]{Lee11}.
\begin{definition}[Submanifold with Boundary, Boundary, Interior]
\label{def:manifold_with_boundary}

A closed subset $M \subset \R^D$ is a $d$-dimensional \emph{$\mathcal{C}^2$-submanifold with boundary} of $\R^D$, if for all $p \in M$ and all small enough open neighborhood $V_p$ of $p$ in $\R^D$,
there exists an open neighborhood $U_0$ of $0$ in $\R^D$ and a $\mathcal{C}^2$-diffeomorphism $\Psi_p : U_0 \to V_p$ with $\Psi_p(0) = p$, such that either:
\begin{enumerate}[label=(\roman*)]
\item
\label{item:interior_point}
$\Psi_p\left(U_0 \cap \left(\R^{d} \times \{0\}^{D-d}\right)\right) = M \cap V_p$.\\
Such a $p \in M$ is called an \emph{interior} point of $M$, the set of which is denoted by $\Int M$.
\item
\label{item:boundary_point}
$\Psi_p\left(U_0 \cap \left(\R^{d-1}\times \R_+ \times \{0\}^{D-d}\right)\right) = M \cap V_p$.\\
Such a $p \in M$ is called a \emph{boundary} point of $M$, the set of which is denoted by $\partial M$.
\end{enumerate}
\end{definition}

\begin{remark}[Boundaries]
\label{rem:def_manifold_with_boundary}
The \emph{geometric} (or \emph{differential}) boundary $\partial M$ is not to be confused with the ambient \emph{topological} boundary defined as $\bar{\partial} S := \bar{S} \setminus \mathring{S}$ for $S \subset \R^D$, where the closure and interior are taken with respect to the ambient topology of $\R^D$.
Indeed, one easily checks that if $d<D$, then $\bar{\partial} M = M$. 
On the other hand, the two sets $\bar{\partial} M$ and $\partial M$ coincide when $d = D$.
\end{remark}

Then, submanifolds \textit{without} boundary are those $M$ that fulfill $\partial M = \emptyset$, i.e. that are everywhere locally parametrized by $\R^d$, and nowhere by $\R^{d-1}\times \R_+$.
From this perspective --- as confusing as this standard terminology can be ---, submanifolds without boundary are special cases of submanifolds with boundary.
Note that key instances of manifolds without boundary are given by boundaries of manifolds, as expressed by the following result.

\begin{proposition}[{\cite[p.30]{Hirsch76}}]
\label{prop:boundary_is_manifold}
If $M \subset \R^D$ is a $d$-dimensional $\mathcal{C}^2$-submanifold with nonempty boundary $\partial M$, then $\partial M$ is a $(d-1)$-dimensional $\mathcal{C}^2$-submanifold without boundary.
\end{proposition}
\begin{remark}
If non-empty, this fact will allow us to estimate $\partial M$ using the estimator designed for manifolds without boundary from \cite{Aamari18},
that we will build on top of some preliminarily filtered  \emph{boundary observations} (see \Cref{sec:boundary_points}).
\end{remark}

\subsubsection{Tangent and normal structures}
\label{sec:tangent_and_normal_structures}
\begin{figure}[h]
\centering
\includegraphics[width=0.6\textwidth,page=2]{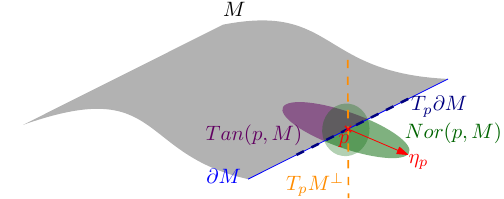}
\caption{
Tangent and normal structure of a surface ($d=2$) in space ($D=3$) at a boundary point.
}
\label{fig:normal-vector}
\end{figure}

In the present $\mathcal{C}^2$-smoothness framework, the difference between boundary and interior points sharply translates in terms of local first order approximation properties of $M$ either by its so-called tangent cones or tangent spaces, which we now define (see \Cref{fig:normal-vector}).

\begin{definition}[Tangent and Normal Cones and Spaces]
\label{def:tangent_normal_cones}
Let $p \in M$,  and $\Psi_p$ its local parametrization from Definition \ref{def:manifold_with_boundary}. 
\begin{itemize}[leftmargin=*]
\item
The \emph{tangent cone} $Tan(p,M)$ of $M$ at $p$ is defined as
\[
Tan(p,M) := \left \{ \begin{array}{@{}cc}
d_0 \Psi_p ( \R^d \times \{0\}^{D-d}) & \text{~~\:if } p \in \Int M, \\
d_0 \Psi_p \left(\R^{d-1}\times \R_{+} \times \{0\}^{D-d}\right) &  \text{if } p \in \partial M,
\end{array} \right .
\]
where $d_0 \Psi_p$ denotes the differential of $\Psi_p$ at $0$.\\
The \emph{tangent space} $T_p M$ is then defined as the linear span $T_p M := \Span(Tan(p,M))$.
\item
The \emph{normal cone} $Nor(p,M)$ of $M$ at $p$ is the dual cone of $Tan(p,M)$:
\[
Nor(p,M) := \{ v \in \R^D \mid \forall u \in Tan(p,M), \left\langle u,v \right\rangle \leq 0 \}.
\]
The \emph{normal space} of $M$ at $p$ is defined accordingly by $N_p(M) := \Span(Nor(p,M))$.
\end{itemize}
\end{definition} 
Whenever $p \in \Int M$, it falls under the intuition that $Tan(p,M) = T_p M$ and $Nor(p,M) = N_pM$, while when $p \in \partial M$, $N_p M$ and $T_p M$ share one direction which is orthogonal to $T_p \partial M$. 
These properties are summarized in the following proposition (see \Cref{fig:normal-vector}).
\begin{proposition}[Outward-Pointing Vector]
\label{prop:cones_of_manifold}
Let $M$ be a $\mathcal{C}^2$-submanifold with boundary.
\begin{itemize}[leftmargin=*]
\item 
If $p \in \Int M$, then $Tan(p,M) = T_p M$ and $Nor(p,M) = N_p M$ are orthogonal linear spaces spanning $\R^D$.
\item
If $p \in \partial M$, then $Tan(p,M)$ and $Nor(p,M)$ are complementary half-spaces, in the sense that $T_p M + N_p M = \R^D$ and $T_p M \cap N_p M$ is one-dimensional. The unique unit vector $\eta_p$ in $Nor(p,M)\cap T_p M$
is called the \emph{outward-pointing} vector.
It satisfies 
\begin{align*}
\hfill
Tan(p,M) =
T_p M \cap \set{ \inner{\eta_p}{.} \leq 0}
,
\hspace{1em}
Nor(p,M) 
= 
N_p M \cap \set{\inner{\eta_p}{.} \geq 0}
,
\hfill
\end{align*}
and 
\begin{align*}
T_p \partial M \Operp \Span(\eta_p) = T_p M
,
\end{align*}
where $\Operp$ denotes the orthogonal direct sum relation.
\end{itemize}
\end{proposition}
The proof of Proposition \ref{prop:cones_of_manifold} derives from elementary differential calculus and is omitted. 
The above purely differential definition of the tangent and normal cones coincides with that of the general framework of sets with positive reach \cite {Federer59} (to follow in \Cref{sec:geometric-assumptions-and-statistical-model}).
This general framework will enable us to quantify how well $M$ is locally approximated by its tangent cones.

\subsection{Geometric assumptions and statistical model}
\label{sec:geometric-assumptions-and-statistical-model}

Any $\mathcal{C}^2$-submanifold $M$ of $\R^D$ admits a tubular neighborhood in which any point has a unique nearest neighbor on $M$ \cite[p.93]{Bredon93}.
However, the width of this tubular neighborhood might be arbitrarily small. This scenario occurs when $M$ exhibits high curvature or nearly self-intersecting areas~\cite{Aamari19}. In this case, the estimation of $M$ gets more difficult, since such locations require denser sample to be reconstructed accurately.
The width of such a tubular neighborhood is given by the so-called \emph{reach} (\cite[Defintion~4.1]{Federer59}), whose formal definition goes as follows.

Given a closed set $S \subset \R^D$, the \emph{medial axis} $\Med(S)$ of $S$ is the set of ambient points that do not have a unique nearest neighbor on $S$.
More precisely, if 
$$\dd(z,S) := \min_{x \in S} \norm{z-x}$$ 
stands for the \emph{distance function} to $S$, then
\begin{equation}
\label{eq:medialaxis}
\Med(S)
:=
\set{
z \in \R^D
|
\exists x\neq y \in S, \norm{z-x} = \norm{z-y} = \dd(z,S)
}
.
\end{equation}
The reach of $S$ is then defined as the minimal distance from $S$ to $\Med(S)$.
\begin{definition}[Reach]
\label{def:reach}
The \emph{reach} of a closed set $S \subset \R^D$ is
\begin{align*}
\tau_S
:=
\min_{x \in S} \dd\left(x,\Med(S) \right)
=
\inf_{z \in \Med(S)} \dd\left(z,S \right).
\end{align*}
\end{definition}

By construction of the medial axis \Cref{eq:medialaxis}, the projection on~$S$
$$
\pi_S(z) := \argmin_{x \in S} \norm{x-z}
$$ 
is well defined (exactly) on $\R^D \setminus \Med(S)$.
In particular, $\pi_S$ is well defined on any $r$-neighborhood of $S$ of radius $r < \tau_S$.
\begin{remark}
One easily checks that $S$ is convex if and only if $\tau_S = \infty$ \cite[Remark 4.2]{Federer59}.
In particular, for the empty set $S = \emptyset$, we have $\tau_\emptyset = \infty$.
\end{remark}

Requiring a lower bound on the reach of a manifold amounts to bound its curvature~\cite[Proposition~6.1]{Niyogi08}, and prevents quasi self-intersection at scales smaller than the reach~\cite[Theorem~3.4]{Aamari19}. 
Moreover, it allows to assess the quality of the linear approximation of the manifold by its tangent cones. 
In fact,~\cite[Theorem~4.18]{Federer59} shows that for all closed set $S \subset \R^D$ with reach $\tau_S>0$, its tangent cone $Tan(x,S)$ is well defined at all $x \in S$, and $\dd(y-x,Tan(x,S)) \leq \norm{y-x}^2/(2\tau_S)$ for all $y \in S$.
This motivates the introduction of our geometric model below.   

\begin{definition}[Geometric Model]
\label{def:geometric_model}
Given integers $1 \leq d \leq D$ and positive numbers $\tau_{\min},\tau_{\partial,\min}$, we let $\mathcal{M}^{d,D}_{\tau_{\min},\tau_{\partial,\min}}$ denote the set of compact connected $d$-dimensional $\mathcal{C}^2$-submanifolds $M \subset \R^D$ with boundary, such that
$$\tau_M \geq \tau_{\min}
\text{~~and~~}
\tau_{\partial M} \geq \tau_{\partial,\min}
.
$$
\end{definition}

\begin{remark}
Let us emphasize the following properties of the model:
\begin{itemize}[leftmargin=*]
\item
The model $\mathcal{M}^{d,D}_{\tau_{\min},\tau_{\partial,\min}}$ includes both submanifolds with empty and non-empty boundary $\partial M$, the main requirement being that $\tau_{\partial M} \geq \tau_{\partial,\min}$. If $\partial M = \emptyset$, this requirement is always fulfilled since $\tau_\emptyset = \infty$.
Note also that \Cref{def:geometric_model} does not exclude the case $d = D$, in which case $M$ consists of a domain of $\R^D$ with non-empty interior.
Furthermore, since the boundary $\partial M$ of a submanifold $M$ is either empty or itself a submanifold without boundary, a non-empty $\partial M$ cannot be convex~\cite[Theorem 3.26]{Hatcher02}.
As a result, $\mathcal{M}^{d,D}_{\tau_{\min},\infty}$ is exactly the set of submanifolds $M \in \mathcal{M}^{d,D}_{\tau_{\min},\tau_{\partial,\min}}$ that have empty boundary. 
In particular, \Cref{def:geometric_model} encompasses the model of \cite{Genovese12a,Kim2015,Aamari18}.

\item
Similarly, since $\tau_M = \infty$ if and only if $M$ is convex, $\mathcal{M}^{d,D}_{\infty,\tau_{\partial,\min}}$ is exactly the set of submanifolds $M \in \mathcal{M}^{d,D}_{\tau_{\min},\tau_{\partial,\min}}$ that are convex (and hence have non-empty boundary).
In particular, \Cref{def:geometric_model} encompasses the model of \cite{Dumbgen96}.

\item
In full generality, the two lower bounds on the respective reaches of $M$ and $\partial M$ are \emph{not} redundant with one another.
As shown in \Cref{fig:reach_boundary_comparison}, $\tau_M$ and $\tau_{\partial M}$ are not related when $d<D$.
However, for $d = D$, $\partial M$ is the topological boundary of $M$ (\Cref{rem:def_manifold_with_boundary}).
In this case,~\cite[Remark~4.2]{Federer59} and an elementary connectedness argument show that $\tau_M \geq \tau_{\partial M}$.
Said otherwise, this means that the reach regularity of a full-dimensional domain is no worse than that of its boundary.
Hence,  
$
\mathcal{M}^{D,D}_{\tau_{\min},\tau_{\partial,\min}} 
= 
\mathcal{M}^{D,D}_{\tau_{\partial,\min},\tau_{\partial,\min}}
$ for all $\tau_{\min} \leq \tau_{\partial,\min}$, so that for $d = D$, one may set $\tau_{\min} = \tau_{\partial,\min}$ without loss of generality.

\begin{figure}[!htp]
	\centering
	\begin{subfigure}{0.31\textwidth}
		\centering
		\includegraphics[width=1\textwidth, page = 1]{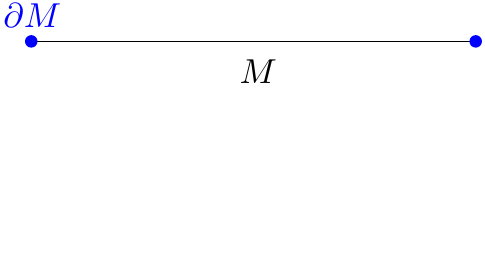}		
		\caption{$\tau_{\partial M} < \tau_M = \infty$.}
		\label{subfig:reach_boundary_comparison_1}
	\end{subfigure}
	\hspace{0.01\textwidth}
	\begin{subfigure}{0.31\textwidth}
		\centering
		\includegraphics[width=1\textwidth, page = 2]{reach_boundary_comparison}		
		\caption{$\tau_{\partial M} = \tau_M$.}
		\label{subfig:reach_boundary_comparison_2}
	\end{subfigure}
	\hspace{0.01\textwidth}
	\begin{subfigure}{0.31\textwidth}
		\centering
		\includegraphics[width=1\textwidth, page = 3]{reach_boundary_comparison}	
		\caption{$\tau_{\partial M} > \tau_M$.}
		\label{subfig:reach_boundary_comparison_3}
	\end{subfigure}
	\caption{
		For $d < D$, the reach of a submanifold $M$ and that of its boundary $\partial M$ are not related.
	}
	\label{fig:reach_boundary_comparison}
\end{figure}

\end{itemize}
\end{remark}
The geometric model $\mathcal{M}^{d,D}_{\tau_{\min}, \tau_{\partial,\min}}$ being settled, we are now in position to define a generative model on such manifolds.
In what follows, we let $\Haus^d$ denote the $d$-dimensional Hausdorff measure on $\R^D$ (see e.g. \cite[Section~2.10.2]{Federer69}).

\begin{definition}[Statistical Model]
\label{def:statistical_model}
Given $0<f_{\min} \leq f_{\max} < \infty$, we let $\mathcal{P}^{d,D}_{\tau_{\min},\tau_{\partial,\min}}(f_{\min},f_{\max})$ denote the set of Borel probability distributions $P$ on $\R^D$ such that:
\begin{itemize}[leftmargin=*]
\item
$M = \supp(P) \in \mathcal{M}^{d,D}_{\tau_{\min},\tau_{\partial,\min}}$,
\item
$P$ has a density $f$ with respect to the volume measure $\vol_M = \indicator{M} \Haus^d$ on $M$, such that $f_{\min} \leq f(x) \leq f_{\max}$ for all $x \in M$.
\end{itemize}
\end{definition}
From now on, we assume that we observe an i.i.d. $n$-sample $X_1, \hdots, X_n$ with unknown common distribution $P \in \mathcal{P}^{d,D}_{\tau_{\min},\tau_{\partial,\min}}(f_{\min},f_{\max})$, and denote the sample point cloud by 
$$
\mathbb{X}_n := \{X_1, \hdots, X_n\}
.
$$
Based on $\X_n$, the performance of the estimators of $M$ and $\partial M$ will be assessed in Hausdorff distance, which plays the role of a $L^\infty$-distance between compact subsets of $\R^D$.
\begin{definition}[Hausdorff Distance]
    \label{def:hausdorff_distance}
        Given two compact subsets $S,S' \subset \R^D$, the \emph{Hausdorff distance} between them is
        \begin{align*}
            \dHaus(S,S')
            &:=
            \max\bigl\{
            \max_{x \in S} \dd(x,S') 
            ,
            \max_{x'\in S'} \dd(x',S) 
            \bigr\}
           .
        \end{align*}
    \end{definition}

%% file: results.tex
\section{Main results}
	\label{sec:main-results}
	
This section gathers the main results of this article: construction of estimators of $\partial M$ and $M$, bounds on their Hausdorff performance, and nearly matching minimax lower bounds.
To cope with the possible presence of a boundary, our first step is to determine which data points lie close to the boundary, if any.

\subsection{Detecting boundary observations}\label{sec:boundary_points}

\subsubsection{Intuition}
\label{sec:intuition_detection}

\begin{figure}[!ht]
\centering
\includegraphics[width=0.7\textwidth]{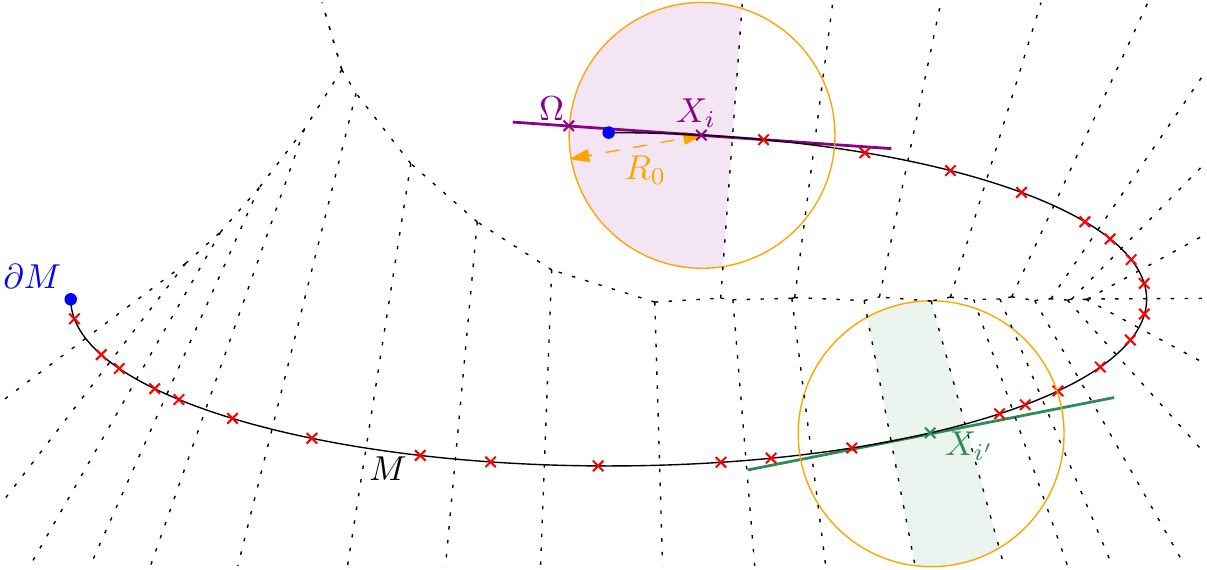}
\caption{
An ambient Voronoi diagram built on top of observations $\mathbb{X}_n$ lying on an open plane curve ($d=1$, $D=2$). The denser $\mathbb{X}_n$ in $M$, the narrower the Voronoi cell of the $X_i$'s in the tangent directions $T_{X_i} M$. 
Observations close to $\partial M$ yield cells that extend in the outward pointing direction.
Localization radius $R_0>0$ prevents global foldings of $M$ that would mix different ambient neighborhoods of $M$ when projecting onto $T_{X_i} M$.
}
\label{fig:voronoi_insight}
\end{figure}

In the full-dimensional case ($d=D$), data points close to the boundary may be identified by how (macroscopically) large their Voronoi cells tend to be~\cite{Casal07}.
That is, if $\rho>0$ is a detection radius, the \emph{boundary observations} may be defined as 
\[
\mathcal{Y}_\rho = \{ X_i \in \X_n \mid \exists O \in \R^D, \norm{O - X_i} \geq \rho \text{ and } \Bopen(O,\norm{O-X_i}) \cap \mathbb{X}_n = \emptyset \}.
\]
If $X_i$ belongs to $\mathcal{Y}_\rho$ with associated $O\in \R^D$, then $\hat{\eta}_i:=\frac{O-X_i}{\norm{O-X_i}}$ appears to provide a consistent estimator of the unit outer normal vector of $\partial M$ at $\pi_{\partial M}(X_i)$ \cite{Aaron20}. 
The present work leverages the above intuition and extends it to the case where $M$ is a $d$-dimensional manifold with $d<D$. 
In fact, the manifold $M$ not being full-dimensional raises the following additional subtleties:
\begin{itemize}[leftmargin=*]
\item
Even if $X_i$ is far from $\partial M$, its Voronoi cell is large in the directions of $T_{X_i}M^\perp$, as it actually contains at least $X_i+\B_{T_{X_i}M^\perp}(0,\tau_{\min})$.
To detect points close to the boundary \emph{only}, we shall hence avoid these normal non-informative directions and solely focus on the tangential components of the Voronoi cells. For instance, by first projecting points onto (an estimate of) $T_{X_i} M$.
\item
If $X_i$ is close to $\partial M$ but $M$ is folded over $X_i$, then the Voronoi cell of $X_i$ in the Voronoi diagram of the projected sample might be small (see \Cref{fig:voronoi_insight}).
To detect enough points close to the boundary, not all the sample should thus be projected, but rather just a neighborhood $\mathbb{X}_n \cap \B(X_i,R_0)$ of $X_i$, for some localization radius $R_0>0$ to be tuned.
\end{itemize}

These two remarks lead to the following first detection procedure: for a collection of estimated tangent spaces $\hat{T}_i$'s, one \emph{may} label $X_i$ as being a \emph{boundary observation} if it has a large Voronoi cell within its $R_0$-neighborhood, when projected onto $X_i + \hat{T}_i$. 
That is, if there exists $O\in \hat{T}_i$ such that $\norm{O}\geq \rho$ and $\Bopen(O,\norm{O})\cap \pi_{\hat{T}_i}(\X_n\cap \B(X_i,R_0)-X_i)=\emptyset$.
Unfortunately, when $1 < d < D$, this intuitive detection method is not sufficient to provably detect enough observations close to the boundary.
This issue can be overcome by investigating \emph{all} the Voronoi cells of $\pi_{\hat{T}_j}(X_i)$ for $X_j \in \B(X_i,r) \cap \X_n$, where $r$ is a small scale parameter.
The details of this detection procedure are given in \Cref{sec:the_detection_method}.

As it is now clear how critical the knowledge of tangent spaces is to build a Voronoi-based boundary detection scheme, let us first briefly detail how we estimate them.
\subsubsection{Tangent space estimation}
Following the ideas of \cite{Aamari18}, we will estimate tangent spaces using local principal component analysis. 

\begin{definition}[Tangent Space Estimator]
\label{def:tangent_spaces_estimator}
For $i \in  \set{1,\ldots,n}$ and $h>0$, we introduce the local covariance matrix 
\[
\hat{\Sigma}_i(h) := \frac{1}{n-1}\sum_{j \neq i} (X_j-X_i)(X_j-X_i)^t\indicator{\B(X_i,h)}(X_j),
\]
and define $\hat{T}_i$ as the linear span of the first $d$ eigenvectors of $\hat{\Sigma}_i(h)$.
\end{definition}

Note that $\hat{T}_i$ is a local estimator, in the sense that it is $\bigl( (X_j-X_i) \indicator{X_j \in \B(X_i,h)}\bigr)_{1 \leq j \leq n}$-measurable (i.e. it only depends on the observations that are $h$-close to $X_i$). For a suitable choice of $h$, the following proposition provides guarantees on the principal angle between $T_{X_i} M$ and $\hat{T}_i$.
In what follows, given two linear subspaces $T,T'\subset \R^D$, the \emph{principal angle} between them is
$$
\angle(T,T')
:=
\normop{\pi_{T} - \pi_{T'}},
$$ where $\normop{A} := \sup_{\norm{x} \leq 1} \norm{Ax}$ stands for the operator norm of $A \in \R^{n \times n}$.

\begin{restatable}[Tangent Space Estimation]{proposition}{proptangentspaceestimation}
\label{prop:tangent_space_estimation}
Let $h = \bigl( C_d \frac{f_{\max}^4}{f_{\min}^5} \frac{\log n }{n-1} \bigr)^{\frac{1}{d}}$, for a large enough constant $C_d$. For $n$ large enough so that $h \leq \frac{\tau_{\min}}{32}\wedge \frac{\tau_{\partial, \min}}{3} \wedge \frac{\tau_{\min}}{\sqrt{d}}$, with probability larger than $1-2 \bigl( \frac{1}{n} \bigr)^{\frac{2}{d}}$, we have
\begin{align*}
\max_{1 \leq i \leq n}
\angle (T_{X_i} M,\hat{T}_i) \leq C_d \frac{f_{\max}}{f_{\min}} \frac{h}{\tau_{\min}}.
\end{align*}
\end{restatable}
A proof of \Cref{prop:tangent_space_estimation} can be found in \Cref{sec:proof_tangent_space_estimation}. In what follows, we shall always choose $h$ and $n$ large enough as in \Cref{prop:tangent_space_estimation}.

\subsubsection{Detection method and normal vector estimation}
\label{sec:the_detection_method}
Now, for a local (though macroscopic) scale $R_0>0$, a detection radius $\rho>0$ and a local bandwidth $r>0$, we compute the $d$-dimensional Voronoi diagrams of $(\pi_{\hat{T}_i}(\B(X_i,R_0) \cap \mathbb{X}_n - X_i ) )_{1 \leq i \leq n}$ and define our boundary observations detection procedure as follows.
See \Cref{sec:intuition_detection} for a heuristic, and \Cref{fig:boundary_points_illustrated_2d} for an illustration associated with this definition.

\begin{definition}[Boundary Observations]\label{def:boundary_points}
For $i \in \set{1,\ldots,n}$, we let $J_{R_0,r,\rho}(X_i)$ be the set of $r$-neighbors $X_j$ of $X_i$ for which $X_i$ has a $\rho$-large Voronoi cell in the projected Voronoi diagram at $X_j$. 
That is, writing
\begin{align*}
\vor^{(j)}_{R_0}(X_i)
:=
\set{
O \in \hat{T}_j
\left|
\Bopen\bigl(O,\Vert O-\pi_{\hat{T}_j}(X_i-X_j)\Vert\bigr)
\cap 
{\pi}_{\hat{T}_j}(\B(X_j,R_0)\cap \mathbb{X}_n-X_j)
=\emptyset
\right.
}
,
\end{align*}
we define
\begin{align*}
J_{R_0,r,\rho}(X_i)
:=
\set{
X_j \in \B(X_i,r) \cap \X_n 
\left| 
	\vor^{(j)}_{R_0}(X_i)
	\cap
	\Bopen_{\hat{T}_j}(\pi_{\hat{T}_j}(X_i-X_j),\rho)^c
	\neq
	\emptyset
\right.
}.
\end{align*}

\begin{figure}[!ht]
\centering
\includegraphics[width=0.8\textwidth]{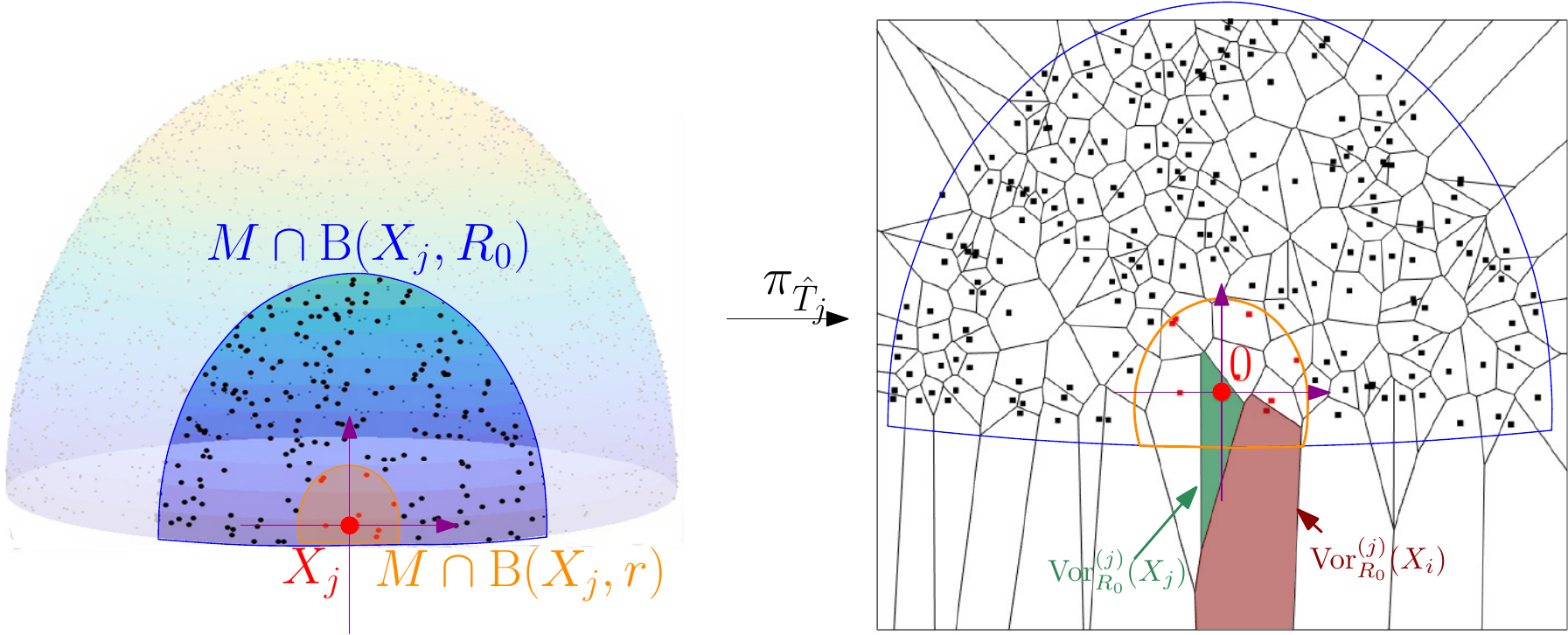}
\caption{
Illustration of \Cref{def:boundary_points} over a half-sphere ($d=2$,$D=3$).
Although the central point $X_j$ (red) does not have a large Voronoi cell in $\hat{T}_j$, its neighbor $X_i$ does. Therefore, $X_j$ belongs to $J_{R_0,r,\rho}(X_i)$. In particular, $X_i$ is labelled as a boundary point.
Note that throughout the process, the Voronoi diagram is only computed in the $d$-planes $\{\hat{T}_j\}_{1 \leq j \leq n}$, not in the ambient space $\R^D$.
}
\label{fig:boundary_points_illustrated_2d}
\end{figure}

The set of \emph{boundary observations} $\mathcal{Y}_{R_0,r,\rho} \subset \X_n$ is then defined as the set of data points that have at least one such large Voronoi cell:
\begin{equation}\label{defboundObs}
\mathcal{Y}_{R_0,r,\rho}
:=
\{X_i \in \X_n \mid J_{R_0,r,\rho}(X_i)\neq \emptyset \}.
\end{equation}
\end{definition} 

\begin{remark}
Detecting boundary observations requires to compute $n$ Voronoi diagrams in dimension $d$. Note that this step does not depend on the ambient dimension $D$, and can run in parallel. 
\end{remark}

This strategy also provides a natural way to estimate unit normal outward-pointing vectors. For this, given a boundary observation $X_i \in \mathcal{Y}_{R_0,r,\rho}$, we simply consider directions in which $\vor^{(j)}_{R_0}(X_i)$ is $\rho$-wide (see \Cref{fig:voronoi_insight}).
A formal definition goes as follows.

\begin{definition}[Normal Vector Estimator]\label{defi:normal_vector_estimate}
For $X_i \in \mathcal{Y}_{R_0,r,\rho}$ and $X_j \in J_{R_0,r,\rho}(X_i)$,~let 
\[
\Omega^{(j)}_{R_0,r,\rho}
\in \argmin
	\set{
	\norm{\Omega-{\pi}_{\hat{T}_j}(X_i-X_j)}
	\left|
	\Omega \in \vor^{(j)}_{R_0}(X_i)
	\cap
	\Bopen_{\hat{T}_j}(\pi_{\hat{T}_j}(X_i-X_j),\rho)^c
	\right.
	}
	.
\]
The estimator of the unit normal outward-pointing vector in $\hat{T}_j$ is defined as
\[
 \tilde{\eta}_i^{(j)}
 :=
 \frac{\Omega^{(j)}_{R_0,r,\rho}-{\pi}_{\hat{T}_j}(X_i-X_j)}{\norm{\Omega^{(j)}_{R_0,r,\rho}-{\pi}_{\hat{T}_j}(X_i-X_j)}}.
 \]
The final estimator of the unit outward-pointing normal vector at $X_i$ is then defined as
\begin{equation}\label{eq:defmeaneta}
\tilde{\eta}_i
:=
\frac{1}{\#J_{R_0,r,\rho}(X_i)}\sum_{j \in J_{R_0,r,\rho}(X_i)} \tilde{\eta}_i^{(j)}.
\end{equation}
\end{definition}

\begin{remark}
Let us mention that the choice of $\Omega^{(j)}_{R_0,r,\rho}$ in \Cref{defi:normal_vector_estimate} has been made to ensure measurability.
As will be clear in the proofs (see \Cref{lem:deterministic_protection_2}), \emph{any} choice of $\Omega \in \vor^{(j)}_{R_0}(X_i)
	\cap
	\Bopen_{\hat{T}_j}(\pi_{\hat{T}_j}(X_i-X_j),\rho)^c$ witnessing to the $\rho$-width of the Voronoi cell would lead to the same normal estimation rates as $\tilde{\eta}_i$.
\end{remark}

As expected, when localization radii are chosen properly, \Cref{thm:detection} below provides quantitative bounds for boundary detection and normal estimation.

\begin{theorem}[Guarantees for Boundary Detection and Normals]\label{thm:detection}
Take $R_0 \leq \frac{\tau_{\min}\wedge \tau_{\partial, \min}}{40}$.
Define 
\[
r_-
:=
\sqrt{(\tau_{\min}\wedge \tau_{\partial, \min}) R_0} \left ( c_d \frac{f_{\max}^5 \log n }{f_{\min}^6 n R_0^d} \right )^\frac{1}{d+1}
, 
r_+ :=\frac{R_0}{12}
, \text{ and } 
\rho_-
:=
\frac{R_0}{4}
=:
\frac{\rho_+}{2}
.
\]
Then, for $n$ large enough, with probability at least $1-4n^{-\frac{2}{d}}$, we have that for all $\rho \in [\rho_-, \rho_+]$ and $r \in [r_-, r_+]$:
\begin{enumerate}[label=(\roman*),leftmargin=*]
\item\label{item:thm_detection_emptyboundary}
If $\partial M = \emptyset$, then $\mathcal{Y}_{R_0,r,\rho} = \emptyset$;
\item
If $\partial M \neq \emptyset$ then:
\begin{enumerate}
\item 
\label{item:thm:detection1}
For all $X_i \in \mathcal{Y}_{R_0,r,\rho}$, 
\[
\dd(X_i, \partial M) \leq \frac{2r^2}{\tau_{\min} \wedge \tau_{\partial, \min}}
;
\]
\item 
\label{item:thm:detection3}
For all $x \in \partial M$,
\[
\dd(x,\mathcal{Y}_{R_0,r,\rho}) \leq 3r;
\]
\item 
\label{item:thm:detection2}
For all $X_i \in \mathcal{Y}_{R_0,r,\rho}$, 
\[
\Vert \eta_{\pi_{\partial M}(X_i)} - \tilde{\eta}_i \Vert  \leq \frac{20r}{\sqrt{R_0 (\tau_{\min}\wedge \tau_{\partial, \min})}} 
.
\]
\end{enumerate}
\end{enumerate}
\end{theorem}

\begin{remark}
Key quantities in \Cref{thm:detection} are the scale $R_0$ and the local bandwith $r$, that need to be carefully tuned in practice. Whenever prior information on the reaches $\tau_{\min}$ and $\tau_{\partial, \min}$ is at hand, we may choose $R_0$ as large as $\frac{\tau_{\min}\wedge \tau_{\partial, \min}}{40}$. Then, an optimal choice $r=r_-$ leads to the bounds: 
\begin{itemize}
\item[\ref{item:thm:detection1}] For all $X_i \in \mathcal{Y}_{R_0,r,\rho}$, 
\[
\dd(X_i, \partial M) \leq (\tau_{\min} \wedge \tau_{\partial, \min}) \left (C_d \frac{ f_{\max}^5}{f_{\min}^5} \frac{\log n }{nf_{\min} (\tau_{\min} \wedge \tau_{\partial, \min})^d} \right )^{\frac{2}{d+1}},
\]
\item[\ref{item:thm:detection3}] For all $x \in \partial M$, 
\[
\dd \left ( x, \mathcal{Y}_{R_0,r,\rho} \right ) \leq 
(\tau_{\min} \wedge \tau_{\partial, \min})
\left (C_d \frac{ f_{\max}^5}{f_{\min}^5} \frac{\log n }{nf_{\min} (\tau_{\min} \wedge \tau_{\partial, \min})^d} \right )^{\frac{1}{d+1}},
\]
\item[\ref{item:thm:detection2}] For all $X_i \in \mathcal{Y}_{R_0,r,\rho}$,
\[
\Vert \eta_{\pi_{\partial M}(X_i)} - \tilde{\eta}_i \Vert \leq   \left (C_d \frac{ f_{\max}^5}{f_{\min}^5} \frac{\log n }{nf_{\min} (\tau_{\min} \wedge \tau_{\partial, \min})^d} \right )^{\frac{1}{d+1}}.
\]
\end{itemize}
\end{remark}

The proof of \Cref{thm:detection} is given in \Cref{sec:proof_thm_detection}. In a nutshell, \Cref{item:thm_detection_emptyboundary} guarantees that no false positive occur if $\partial M = \emptyset$. On the other hand, if $\partial M \neq \emptyset$, for $\varepsilon \asymp (\log n/n)^{1/(d+1)}$ and optimal choices of $r_-$ and $R_0$, \Cref{item:thm:detection1,item:thm:detection3} ensure that $\mathcal{Y}_{R_0,r,\rho}$ is an $O(\varepsilon)$-covering of $\partial M$
that consists of points 
$O(\varepsilon^2)$-close to $\partial M$. 

In the convex case $\tau_{\min} = \infty$, taking the convex hull of $\mathcal{Y}_{R_0,r,\rho}$ --- similarly to~\cite{Dumbgen96} --- would result in an $O(\varepsilon^2)$-approximation of $M$, and the boundary of this convex hull in an $O(\varepsilon^2)$-approximation of $\partial M$.
Finally, \Cref{item:thm:detection2} asserts that the estimated normals at boundary observations are $O(\varepsilon)$-precise.

The intuition behind the respective rates $O(\varepsilon)$ and $O(\varepsilon^2)$ is the same as in the convex case of \cite{Dumbgen96}: for a fixed boundary point $x \in \partial M$, the ``curved rectangle'' $\{u \in \B(x,\varepsilon) \cap M \mid \dd(u,\partial M) \leq \varepsilon^2\}$ has volume of order $\varepsilon^{d-1} \times \varepsilon^2 = \varepsilon^{d+1}$. 
Hence, the choice $\varepsilon \asymp (\log n/n)^{1/(d+1)}$ ensures that these curved rectangular regions are occupied by sample points with high probability, uniformly over the choice of $x$ on a grid.
Our procedure then guarantees that these close-to-boundary sample points will be identified as such. 

\begin{remark}
The above argument may be pushed further to gain insights on the number $|\mathcal{Y}_{R_0,r,\rho}|$ of detected points. 
\begin{itemize}
    \item
To derive an upper bound, use \Cref{item:thm:detection1} to get that $|\mathcal{Y}_{R_0,r,\rho}| \leq \sum_{i=1}^n \indicator{\dd(X_i, \partial M) \lesssim \varepsilon^2}$ with high probability, where $\varepsilon \asymp (\log n/n)^{1/(d+1)}$. 
Since $P(\{u \in M \mid \dd(u,\partial M) \leq \varepsilon^2\}) \lesssim \mathrm{Vol}_{d-1}(\partial M) \varepsilon^2$ for $\varepsilon$ small enough, this ensures that 
$$|\mathcal{Y}_{R_0,r,\rho}| \lesssim n \varepsilon^2 \asymp \log n^{2/(d+1)}n^{(d-1)/(d+1)}
.
$$
In particular, for $d=1$, optimal choices of the parameters guarantees that the number of detected points should be no more than roughly $\log n$. In this $1$-dimensional case, it falls under the intuition that the 'optimal' number of detected points should be $2$, corresponding to extremal points drawn on a curve. 
\item
On the other hand,  Items~\ref{item:thm:detection1} and \ref{item:thm:detection3} combined provide a lower bound on $|\mathcal{Y}_{R_0,r,\rho}|$. 
Letting $N(\varepsilon)$ denote the $\varepsilon$-covering number of $\partial M$, standard volume arguments show that $N(\varepsilon) \simeq \mathrm{Vol}_{d-1}(\partial M)/\varepsilon^{d-1}$.
Furthermore, the fact that $\sup_{x \in \partial M} \dd(x,\mathcal{Y}_{R_0,r,\rho}) \lesssim \varepsilon$ means that $\mathcal{Y}_{R_0,r,\rho}$ is a $O(\varepsilon)$-covering of $\partial M$, so that $|\mathcal{Y}_{R_0,r,\rho}| \geq N(\varepsilon)$. Hence, we obtain 
$$
|\mathcal{Y}_{R_0,r,\rho}| \geq N(\varepsilon) \gtrsim (\log n)^{-(d-1)/(d+1)}n^{(d-1)/(d+1)}.
$$ 
\end{itemize}
These two matching bounds (up to $\log n$ factors) back the intuition that as $d$ grows large, most of the mass (and hence sample) is concentrated nearby the boundary.
\end{remark}

If no prior information on $\tau_M$ and $\tau_{\partial M}$ are available, choosing $R_0 = (\log n)^{-1}$ would meet the requirements of \Cref{thm:detection} for $n$ large enough. As well, choosing $r= \sqrt{R_0 \log n } \left( \log n / (nR_0^d) \right )^{1/(d+1)}$ would asymptotically meet the requirements of \Cref{thm:detection}. Both of these choices  incur an extra $\log n $ factor in the bounds.

Still based on $\mathcal{Y}_{R_0,r,\rho}$, we extend this ``hull'' construction to the non-convex case by leveraging the additional tangential (\Cref{prop:tangent_space_estimation}) and normal (\Cref{thm:detection}~\ref{item:thm:detection2}) estimates, to provide estimators of $M$ and $\partial M$.

\subsection{Boundary estimation}

Assume that $\partial M \neq \emptyset$. Then $\partial M$ is a $(d-1)$-dimensional $\mathcal{C}^2$-submanifold without boundary.
Therefore, using manifold estimators of \cite{Aamari18, Aamari19b, Maggioni16,Fefferman19} designed for the empty boundary case with input points $\mathcal{Y}_{R_0,r,\rho}$ seems relevant. We choose to focus on the manifold estimator proposed in~\cite{Aamari18}, based on the Tangential Delaunay Complex~\cite{Boissonnat14}, as it also provides a topologically consistent estimation. This procedure, as well as the aforementioned two others, takes as input boundary points but also estimates of the tangent spaces (of the boundary). Thus, a preliminary step is to provide estimators for the boundary tangent spaces at points of $\mathcal{Y}_{R_0,r,\rho}$.
\begin{definition}[Boundary's Tangent Space Estimator]\label{defi:boundary_tangent_space}
For all $X_i \in \mathcal{Y}_{R_0,r,\rho}$, $\hat{T}_{\partial,i}$ is defined as the orthogonal complement of $\pi_{\hat{T}_i}(\tilde{\eta}_i)$ in $\hat{T}_i$. That is,
\[
\hat{T}_{\partial,i} 
:= 
(\pi_{\hat{T}_i}(\tilde{\eta}_i))^\perp \cap \hat{T}_i.
\] 
\end{definition}
A straightforward consequence of \Cref{prop:tangent_space_estimation} and \Cref{thm:detection} is that the estimator $\hat{T}_{\partial,i}$ is a $O\bigl((\log n/n)^{1/(d+1)}\bigr)$-approximation of $T_{\pi_{\partial M}(X_i)} \partial M$, for any $X_i \in \mathcal{Y}_{R_0,r,\rho}$.

\begin{restatable}[Boundary's Tangent Space Estimation]{corollary}{corboundarytangentspaceestimation}\label{cor:boundary_tangent_space_estimation}
Under the assumptions of \Cref{prop:tangent_space_estimation} and \Cref{thm:detection} we have, for $n$ large enough, with probability larger than $1-4n^{-\frac{2}{d}}$,
\[
\max_{X_i \in \mathcal{Y}_{R_0,r,\rho}}
\angle(
T_{\pi_{\partial M}(X_i)} \partial M
,
\hat{T}_{\partial,i}
)
\leq \frac{20r}{\sqrt{(\tau_{\min}\wedge \tau_{\partial, \min})R_0}}.
\]
Thus, choosing $R_0=\frac{\tau_{\min}\wedge \tau_{\partial, \min}}{40}$ and $r=r_-$ yields
\[
\max_{X_i \in \mathcal{Y}_{R_0,{r_-},\rho}}
\angle(
T_{\pi_{\partial M}(X_i)} \partial M
,
\hat{T}_{\partial,i}
)
\leq \left (C_d \frac{ f_{\max}^5}{f_{\min}^5} \frac{\log n }{nf_{\min} (\tau_{\min} \wedge \tau_{\partial, \min})^d} \right )^{\frac{1}{d+1}}.
\]
\end{restatable}
A short proof can be found in \Cref{sec:proof_cor_boundary_tangent_space_estimation}, that connects $\angle(
T_{\pi_{\partial M}(X_i)} \partial M
,
\hat{T}_{\partial,i}
)$ to $\angle(T_{X_i}M,\hat{T}_i)$ and $\angle{(\eta_{\pi_{\partial M}(X_i)},}{\tilde{\eta}_i)}$.  
The estimation rate for  $T_{\pi_{\partial M}(X_i)} \partial M$ is then driven by the larger of these quantities, i.e. $\angle(\eta_{\pi_{\partial M}(X_i)},\tilde{\eta}_i)$ according to \Cref{prop:tangent_space_estimation} and \Cref{thm:detection}.

Equipped with \Cref{cor:boundary_tangent_space_estimation}, we are now in position to provide an estimator for $\partial M$. 
Following~\cite{Aamari18}, we let $\varepsilon=C\frac{\tau_{\partial, \min}}{R_0}r$, where $r$ and $R_0$  are chosen as 
in \Cref{thm:detection}, and let $\mathbb{Y}_\partial$ denote an $\varepsilon$-sparsification of $\mathcal{Y}_{R_0,r,\rho}$, i.e.
a subset of $\mathcal{Y}_{R_0,r,\rho}$ that forms an $\varepsilon$-covering of $\mathcal{Y}_{R_0,r,\rho}$ with $\varepsilon$-separated points. 
Such a sparsification can be obtained by running the farthest point sampling algorithm over $\mathcal{Y}_{R_0,r,\rho}$, and it results in a $2\varepsilon$-covering of $\partial M$, according to \Cref{thm:detection}.
We also denote by $\mathbb{T}_\partial$ the collection of $\hat{T}_{\partial, i}$'s, for $X_i \in \mathbb{Y}_\partial$, and define our estimator of $\partial M$ as the (weighted) 
Tangential Delaunay Complex~\cite{Boissonnat14} based on $(\mathbb{Y}_\partial, \mathbb{T}_\partial)$:
\[
\widehat{\partial M} := \mathrm{Del}^{\omega_*}(\mathbb{Y}_\partial, \mathbb{T}_\partial).
\]
Since $\partial M$ has no boundary, \cite[Theorem~4.4]{Aamari18} applies and yields the following reconstruction result.
\begin{theorem}[Boundary Estimation: Upper Bound]
	\label{thm:boundary-main-upper-bound}
Provided that $\partial M \neq \emptyset$ and under the assumptions of \Cref{prop:tangent_space_estimation} and \Cref{thm:detection}, we have for $n$ large enough, with probability larger than $1-4n^{-\frac{2}{d}}$,
\begin{enumerate}[label=(\roman*)]
\item
$\dHaus(\partial M, \widehat{\partial M}) 
\leq C_d \frac{\tau_{\partial, \min}}{R_0^2}r^2$,
\item $\partial M$ and $\widehat{\partial M}$ are ambient isotopic.
\end{enumerate}
As a consequence, for $n$ large enough, choosing $R_0 = \frac{\tau_{\min}\wedge \tau_{\partial, \min}}{40} $ and $r=r_-$, we have
\[
\E_{P^n}\left[ \dHaus(\partial M,\widehat{\partial M}) \right]
\leq
C_d  \tau_{\partial, \min} \left ( \frac{ f_{\max}^5}{f_{\min}^5} \frac{\log n }{nf_{\min} (\tau_{\min} \wedge \tau_{\partial, \min})^d} \right )^{\frac{2}{d+1}}.
\]
\end{theorem}

The proof derives from a direct application of the reconstruction result of~\cite[Theorem~4.4]{Aamari18}, the assumptions of which hold with high probability, according to the distance bounds of \Cref{thm:detection}~\ref{item:thm:detection1} and~\ref{item:thm:detection3} and the angle bounds of \Cref{cor:boundary_tangent_space_estimation}.

Note that the ambient dimension~$D$ plays no role in \Cref{thm:boundary-main-upper-bound}, neither in the assumptions, the rate nor the constants.
Interestingly, it assesses the topological correctness of our estimator $\widehat{\partial M}$,
showing the particular interest of estimators based on simplicial complexes. 
Choosing the largest possible $R_0$, i.e. 
$R_0 = \frac{\tau_{\min}\wedge \tau_{\partial, \min}}{40} $, and $r=r_-$, 
 \Cref{thm:boundary-main-upper-bound} provides an upper bound on $\dHaus(\partial M,\widehat{\partial M})$ with high probability, uniformly over the class
 $\mathcal{P}^{d,D}_{\tau_{\min}, \tau_{\partial, \min}}(f_{\min}, f_{\max})$ introduced in \Cref{def:statistical_model}. 
This uniform convergence rate is in line with the estimation rate $O\bigl((\log n /n)^{2/(d+1)}\bigr)$ for boundary estimation given by~\cite{Casal07,Dumbgen96}, under convexity-type assumptions in the full dimensional case.
Letting $\tau_{\min} = \infty$, the convex case can even be seen of as a sub-case of our class of distributions, since $\mathcal{P}^{d,D}_{\tau_{\min}, \tau_{\partial, \min}}(f_{\min}, f_{\max}) \supset \mathcal{P}^{d,D}_{\infty, \tau_{\partial, \min}}(f_{\min}, f_{\max})$.
In fact, even in this simpler case, we can show that the rate $O\bigl((\log n /n)^{2/(d+1)}\bigr)$ is minimax over the class of convex submanifolds.

\begin{theorem}[Boundary Estimation: Lower Bound]
	\label{thm:boundary-main-lower-bound}
	Assume that $f_{\min} \leq c_d/\tau_{\partial, \min}^d$, and that $c'_d/\tau_{\partial, \min}^d \leq f_{\max}$ for some small enough $c_d,(c'_d)^{-1}>0$. Then for all $n \geq 1$, 
		\begin{align*}
	\inf_{\hat{B}}
	\sup_{P \in \mathcal{P}^{d,D}_{\infty,\tau_{\partial, \min}}(f_{\min},f_{\max})}
	\E_{P^n}
	\left[
	\dHaus\bigl(\partial M,\hat{B}\bigr)
	\right]
	&\geq
	C_d \tau_{\partial, \min}
	\left\{
		1
		\wedge
		\left(
			\frac{1}{f_{\min} \tau_{\partial, \min}^d n}
		\right)^{\frac{2}{d+1}}
	\right\}
	.
	\end{align*}	
\end{theorem}
A proof of \Cref{thm:boundary-main-lower-bound} is given in \Cref{sec:proof-lower-bounds-main} and relies on standard Bayesian arguments.

Since for all $\tau_{\min}>0$, $\mathcal{P}^{d,D}_{\infty, \tau_{\partial, \min}}(f_{\min}, f_{\max}) \subset \mathcal{P}^{d,D}_{\tau_{\min}, \tau_{\partial, \min}}(f_{\min}, f_{\max})$, \Cref{thm:boundary-main-lower-bound} and \Cref{thm:manifold-main-lower-bound} together ensure that our boundary estimation procedure is minimax over the model $\mathcal{P}^{d,D}_{\tau_{\min}, \tau_{\partial, \min}}(f_{\min}, f_{\max})$, up to $\log n $ factors. 
From a statistical viewpoint, these two results show that estimating the boundary under reach conditions on $M$ is not more difficult than estimating the boundary in the convex case. 

\subsection{Boundary-adaptive manifold estimation}

If $\partial M = \emptyset$, it is known that $M$ can be estimated optimally by local linear patches~\cite{Aamari19b}.
That is, choosing $\varepsilon_{\interior{M}} = \left ( C_d \frac{f_{\max}^4 \log n }{f_{\min}^5 n} \right )^{1/d} $, and estimating $M$ via the union of tangential balls $\hat{M}=\bigcup_{i=1}^n X_i + \B_{\hat{T}_i}(0,\varepsilon_{\interior{M}}  )$ leads to $\dHaus(M,\hat{M}) \leq C_d f_{\max} \varepsilon_{\interior{M}}^2/(f_{\min}\tau_{\min})$ \cite[Theorem~6]{Aamari19b}, recovering the minimax rate $O\bigl((\log n /n)^{2/d}\bigr)$ over the class of $\mathcal{C}^2$ manifolds without boundary \cite{Kim2015}. 

If $\partial M \neq \emptyset$ and $X_i$ is close to $\partial M$, a tangential ball $X_i + \B_{\hat{T}_i}(0,\varepsilon_{\interior{M}})$ may go past $\partial M$ along the normal direction $\eta_{\pi_{\partial M}(X_i)}$, leading to a poor approximation of $M$ in terms of Hausdorff distance. In this case, replacing $X_i+\B_{\hat{T}_i}(0,\varepsilon_{\interior{M}} )$ by a tangential half-ball oriented at the opposite of the outward-pointing normal vector $\eta_{\pi_{\partial M}(X_i)}$ seems more appropriate. 
We formalize this intuition as follows. 

\begin{figure}[!ht]
\centering
\includegraphics[width=0.6\textwidth]{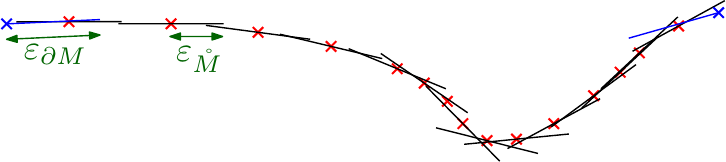}
\caption{
The local linear estimator $\hat{M}$ from \Cref{defi:manifold_estimator} for $d=1$ and $D=2$. The boundary estimator $\hat{M}_{\partial}$ corresponds to the union of the two blue segments, and $\hat{M}_{\Int}$ to that of the black segments.
}
\label{fig:boundary_linear_patches}
\end{figure}

Let $\mathcal{Y}_{R_0,r,\rho}$ denote the detected boundary observations of \Cref{def:boundary_points}. These points will generate half-balls, with radius $\varepsilon_{\partial M}$, that will roughly approximate the inward slab $M \cap \B(\partial M, \varepsilon_{\partial M})$ of radius $\varepsilon_{\partial M}$. To approximate the remaining part of $M$, we further define the $\varepsilon_{\partial M}$-\emph{inner points} as 
\begin{align}\label{eq:inner_points}
\interior{\mathcal{Y}}_{\varepsilon_{\partial M}} 
:= 
\set{
X_i \in \X_n \mid \dd(X_i, \mathcal{Y}_{R_0,r,\rho}) \geq  \varepsilon_{\partial M}/2 
}.
\end{align}
Then, the manifold $M$ may be reconstructed as follows (see \Cref{fig:boundary_linear_patches}).
\begin{definition}[Boundary-Adaptive Manifold Estimator]\label{defi:manifold_estimator} 
Given some scale parameters $\varepsilon_{\interior{M}}$ and $\varepsilon_{\partial M}$, the manifold estimator
$
\hat{M}
:=
\hat{M}_{\Int}
\cup
\hat{M}_{\partial},
$
is defined as
\begin{align*}
\hat{M}_{\Int}
&
:=
\bigcup_{X_i \in \interior{\mathcal{Y}}_{\varepsilon_{\partial M}}} 
X_i + \B_{\hat{T}_i}(0,\varepsilon_{\interior{M}}),
\\
\hat{M}_{\partial}
&
:=
\bigcup_{ X_i \in \mathcal{Y}_{R_0,r,\rho}} 
\left( 
X_i + \B_{\hat{T}_i}(0, \varepsilon_{\partial M}) 
\right) 
\cap 
\{z,\langle z-X_i,\tilde{\eta}_i \rangle \leq 0 \}
,
\end{align*}
with
\begin{itemize}
\item
the $\hat{T}_i$'s being the estimated tangent spaces from \Cref{prop:tangent_space_estimation},
\item
the $\tilde{\eta}_i$'s being the estimated of the outward-pointing normals from \Cref{thm:detection}.
\end{itemize}
\end{definition}

Note that $\hat{M}$ is adaptive in the sense that it does not require information about emptiness of $\partial M$. 
If $\partial M = \emptyset$, then $\mathcal{Y}_{R_0,r,\rho} = \emptyset$ with high probability (\Cref{thm:detection}~\ref{item:thm_detection_emptyboundary}). 
In this case $\hat{M}$ coincides (with high probability) with the estimator from \cite{Aamari19b}, which is minimax over the class of boundariless $\mathcal{C}^2$-manifolds. 
\Cref{thm:manifold-main-upper-bound} below extends the error bound for $\hat{M}$ whenever $\partial M \neq \emptyset$. 
\begin{theorem}[Estimation with Boundary: Upper Bound]\label{thm:manifold-main-upper-bound}
	Choose $(R_0,r,\rho)$ as in \Cref{thm:detection}, set 
	$$
	\varepsilon_{\interior{M}} = \left ( C_d \frac{\log n }{f_{\min} n} \right )^{\frac{1}{d}}
	\text{ and }
	\varepsilon_{\partial M} =  18r
	.
	$$
	Then for $n$ large enough, with probability larger than $1-4n^{-\frac{2}{d}}$, we have 
	\begin{align*}
	\dHaus(M,\hat{M}) 
	\leq  C_d 
	\begin{cases}	
	(f_{\max}/f_{\min})^{\frac{4}{d}+1} \varepsilon_{\interior{M}}^2/\tau_{\min}
	&
	\text{if } \partial M = \emptyset,
	\\
	\varepsilon_{\partial M}^2/R_0
	&
	\text{if } \partial M \neq \emptyset.
	\end{cases}
\end{align*}
As a consequence, for $n$ large enough, with $R_0 = \frac{\tau_{\min}\wedge \tau_{\partial, \min}}{40} $ and $r=r_-$, it holds
\[
\E_{P^n} \left[ \dHaus(M,\hat{M}) \right]
\leq 
C_d 
\begin{cases}
\tau_{\min} 
\left( 
\dfrac{f_{\max}^{2+d/2}}{f_{\min}^{2+d/2}}
\dfrac{\log n}{f_{\min} \tau_{\min}^d n} \right )^{\frac{2}{d}}
&
\text{if } \partial M = \emptyset,
\\
\\
(\tau_{\min}\wedge \tau_{\partial, \min}) \left ( \dfrac{f_{\max}^5 }{f_{\min}^5} \dfrac{\log n}{f_{\min} (\tau_{\min} \wedge \tau_{\partial, \min})^d n} \right )^{\frac{2}{d+1}}
&
\text{if } \partial M \neq \emptyset
.
\end{cases}
\]	
\end{theorem}

A proof of \Cref{thm:manifold-main-upper-bound} is given in \Cref{sec:proof_thm_manifold_upper_bound}.
Again, note that \Cref{thm:manifold-main-upper-bound} is completely oblivious to the ambient dimension $D$.
In the empty boundary case, $\hat{M}$ achieves the rate $O\bigl((\log n / n)^{2/d} \bigr)$, which is minimax \cite{Kim2015}.
Whenever $\partial M$ is not empty, the given convergence rate of $\hat{M}$ coincides with that of $\widehat{\partial M}$ for boundary estimation (\Cref{thm:boundary-main-upper-bound}), as well as that of \cite[Corollary~1]{Dumbgen96} for convex domains, and that of \cite[Theorem~3]{Casal07} for $r$-convex domains.
Note that these last two convexity-type assumptions are stronger than the bounded reach assumption for $M$ and $\partial M$, so that \Cref{thm:manifold-main-upper-bound} generalizes \cite{Dumbgen96,Casal07}.
As for the boundary estimation problem, we show that this rate $O\bigl((\log n / n)^{2/d} \bigr)$ is in fact minimax optimal over the class of $d$-dimensional convex domains (i.e. $\tau_{\min} = \infty$), up to $\log n$ factors.

\begin{theorem}[Manifold Estimation: Lower Bounds]
	\label{thm:manifold-main-lower-bound}
	~
	\begin{enumerate}[label=(\roman*)]
	\myitem{(Boundaryless)}
		\label{item:minimax_lower_bound_boundaryless}
	Assume that $f_{\min} \leq c_d/\tau_{\min}^d$ and that $c'_d/\tau_{\min}^d \leq f_{\max}$, for some small enough $c_d,(c'_d)^{-1}>0$. If $d \leq D-1$, then for all $n \geq 1$,
	\begin{align*}
	\inf_{\hat{M}}
	\sup_{P \in \mathcal{P}^{d,D}_{\tau_{\min},\infty}(f_{\min},f_{\max})}
	\E_{P^n}
	\left[
	\dHaus\bigl(M,\hat{M}\bigr)
	\right]
	&\geq
	C_d \tau_{\min}
	\left\{
		1
		\wedge
		\left(
			\frac{1}{f_{\min} \tau_{\min}^d n}
		\right)^{\frac{2}{d}}
	\right\}
	.
	\end{align*}
	\myitem{(Convex)}
		\label{item:minimax_lower_bound_with_boundary}
	Assume that $f_{\min} \leq c_d/\tau_{\partial, \min}^d$ and $c'_d/\tau_{\partial, \min}^d \leq f_{\max}$, for some small enough $c_d,(c'_d)^{-1}$ $>0$. Then for all $n \geq 1$,
	\begin{align*}
	\inf_{\hat{M}}
	\sup_{P \in \mathcal{P}^{d,D}_{\infty,\tau_{\partial, \min}}(f_{\min},f_{\max})}
	\E_{P^n}
	\left[
	\dHaus\bigl(M,\hat{M}\bigr)
	\right]
	&\geq
	C_d \tau_{\partial, \min}
	\left\{
		1
		\wedge
		\left(
			\frac{1}{f_{\min} \tau_{\partial, \min}^d n}
		\right)^{\frac{2}{d+1}}
	\right\}
	.
	\end{align*}
	\end{enumerate}
\end{theorem}
The proof of \Cref{thm:manifold-main-lower-bound} relies on the same bayesian arguments as \Cref{thm:boundary-main-lower-bound} (see \Cref{sec:proof-lower-bounds-main}).
The first point is a slight refinement of the $\mathcal{C}^2$ case of \cite[Theorem~7]{Aamari19b}, as it exhibits the dependency on $\tau_{\min}$ and $f_{\min}$ of
the minimax rates over the class of $\mathcal{C}^2$ manifolds without boundary.
Note also that in this case, the assumption $d \leq D-1$ clearly is necessary for the model not to be empty.

Interestingly, this shows that the upper bound given in \Cref{thm:manifold-main-upper-bound} for the empty boundary case is sharp with respect to $\tau_{\min}$.  
The second point of \Cref{thm:manifold-main-lower-bound} provides the minimax rate for manifold estimation over the class of convex domains whose boundary has bounded reach. 
In terms of sample size, this shows that our estimator has the best possible convergence rate $O\bigl((\log n / n)^{2/(d+1)} \bigr)$ (up to $\log n$ factors) in the convex case, as well as the two procedures of \cite{Dumbgen96, Casal07}. 
As for the boundary estimation problem, this result intuitively carries the message that estimating a manifold with boundary under reach conditions is not more difficult than estimating a $d$-dimensional convex $\mathcal{C}^2$-domain.
In other words, for $\partial M \neq \emptyset$ and a fixed boundary's convexity radius $\tau_{\partial, \min}$, no additional gain can be expected from requiring a large convexity radius for the manifold (driven by $\tau_{\min}$).
At last, \Cref{thm:manifold-main-upper-bound} shows that the given dependency on the reach boundary $\tau_{\partial, \min}$ is sharp, at least in the case where $\tau_{\partial, \min} \leq \tau_{\min}$. Whether the tradeoff between $\tau_{\min}$ and $\tau_{\partial, \min}$ exhibited in \Cref{thm:manifold-main-upper-bound} is sharp in general remains an open question.

%% file: conclusion.tex
\section{Conclusion and further perspectives}
	\label{sec:conclusion}
	
Both generalizing over full dimensional $\mathcal{C}^2$ domains and boundaryless $\mathcal{C}^2$-submanifolds, this work derives nearly tight minimax upper and lower bounds for $\mathcal{C}^2$-submanifold estimation with possibly non-empty $\mathcal{C}^2$ boundary.
Both the boundary estimator and the manifold estimator exhibit rates that are independent of the ambient dimension, which is of critical interest in the regime $d \ll D$ to achieve efficient dimensionality reduction.
To our knowledge, this is the first instance of a statistical study dealing with general submanifold with boundary.

This work is the first minimax estimation study on manifolds with boundary. Hence, the focus has not been put on computational aspects.
Yet, the proposed method is fully constructive and can easily be implemented using PCA and computational geometric algorithms.
Given the space constraints, we refer the interested reader to \Cref{sec:experiments}, which discusses computational complexity, parameter tuning, and provides a few numerical examples.

On the geometric side, a significant further direction of research pertains to manifold estimation with boundary in smoother models than $\mathcal{C}^2$, such as those introduced in~\cite{Aamari19b}.
Beyond Hausdorff minimax optimality, an interesting feature of the boundary estimator of \Cref{thm:boundary-main-upper-bound} is its topological exactness.
This property is  made possible by the fact that $\partial(\partial M) = \emptyset$ and the existence of constructive triangulations that reconstruct boundaryless submanifolds (see~\cite[Theorem~4.4]{Aamari18}).
In contrast, topologically exact reconstruction methods of manifolds with boundary are only known in the specific case of \emph{isomanifolds} (see~\cite[Theorem~43]{Boissonnat20}), which led us to stick to an unstructured estimator with linear patches in this case (see \Cref{thm:manifold-main-upper-bound}).

On the statistical side, a major limitation of this work is the absence of noise. The proposed method would exhibit the same rates if noise of amplitude $\sigma \ll (\log n / n)^{2/d} \indicator{\partial M  = \emptyset} + (\log n / n)^{2/(d+1)} \indicator{\partial M  \neq \emptyset}$ is added, but it is likely to fail otherwise as it is based on the data points themselves.
Such instabilities are common in the geometric inference literature~\cite{Cuevas04,Aamari19,Berenfeld20,Divol20}, and noise is often assumed to vanish as $n$ goes to $\infty$.
However, a recent line of works in the boundariless case exhibited various iterative denoising procedures that tend to relax this assumption. See for instance \cite{Fefferman19,Puchkin19,Aizenbud21}.
Whether such algorithms could be adapted for $\partial M \neq \emptyset$ is of particular interest.

\section*{Acknowledgments}
We are grateful to the members of the \textit{Laboratoire de Probabilit\'es, Statistique et Mod\'elisation} and the \textit{Laboratoire de Math\'ematiques Blaise Pascal} for their insightful comments.

%% file: proof-outline.tex
\section{Proofs outline}
	\label{sec:proof-outline}

Due to space constraints, the geometric results necessary to the proofs given below are deferred to the appendix (\Cref{sec:geom_results}).

\subsection{Proof of Theorem~\ref{thm:detection}}\label{sec:proof_thm_detection}

The main boundary detection result is based on the following geometric and purely deterministic result.

\begin{theorem}[Deterministic Layout for Boundary Detection and Normals]\label{thm:deterministic}
Let
\begin{align*}
R_0 \leq \frac{\tau_{\min}}{32}
,
r_0 \leq \frac{R_0 \wedge \tau_{\partial,\min}}{4}
,
r \leq \frac{R_0}{12}
,
\\
\theta \leq \frac{1}{24}
,
\eps_1 \leq \frac{r}{4}
,
\eps_2 \leq \frac{r_0}{120} \wedge \frac{r^2}{{\tau_{\min} \wedge \tau_{\partial, \min}}}
,
\\
\text{and } 
3r \leq \rho_- < \rho_+ \leq \frac{{\tau_{\min} \wedge \tau_{\partial, \min}}}{80}.
\end{align*}
Assume that we have:
\begin{enumerate}
\item\label{item:thm:deterministic:sampling}
 A point cloud $\mathcal{X}_n \subset M$ such that $\dHaus(M,\mathcal{X}_n)\leq \eps_1$,
\item\label{item:thm:deterministic:tangents}
Estimated tangent spaces $T_j$ such that $\max_{1 \leq j \leq n} \angle(T_{X_j} M,T_j) \leq \theta$.
\end{enumerate}
For $x\in \partial M$ and $j\in \set{1,\ldots,n}$ such that $\angle(T_xM,T_j)<1$, write $\eta^*_j(x)$ for the unit vector of $Nor(x,M)\cap T_j$ (see \Cref{prop:existencesestimnorM}).
Defining $\mathcal{Y}_{j}:=\pi_{T_j}(\B(X_j,R_0)\cap \mathcal{X}_n - X_j)$, assume furthermore that:
\begin{enumerate}\setcounter{enumi}{2}
\item\label{item:thm:deterministic:normals} 
For all $x\in \partial M$ and $X_j \in \mathcal{X}_n \cap \B(x,2r)$,
for all  $\rho\geq \rho_-$ and $\Omega\in T_j$ such that 
 $\Vert\Omega-(\pi_{T_j}\left( x-X_j \right )-r_0\eta_j^*(x))\Vert\leq r_0+\rho-\eps_2$ we have
 $\B(\Omega,\rho)\cap \mathcal{Y}_{j}\neq \emptyset$.
\end{enumerate}
Then for all $\rho \in [\rho_-,\rho_+]$, using notation of \Cref{def:boundary_points,defi:normal_vector_estimate}, the following holds:
\begin{enumerate}[label=(\roman*),leftmargin=*]
\item\label{deter:empty:boundary}
If $\partial M = \emptyset$, then $\mathcal{Y}_{R_0,r,\rho}=\emptyset$.
\item
If $\partial M \neq \emptyset$, then,
\begin{enumerate}
\item\label{deter:protec:distance}
For all $X_i \in \mathcal{Y}_{R_0,r,\rho}$,
\begin{equation*}
\dd(X_i,\partial M)
\leq 
\frac{2r^2}{{\tau_{\min} \wedge \tau_{\partial, \min}}}
,
\end{equation*}
\item\label{deter:detec}
For all $x\in \partial M$,
\begin{equation*}
\dd(x,\mathcal{Y}_{R_0,r,\rho})
\leq 
3r.
\end{equation*}
\item\label{deter:protec:angle}
For all $X_i \in \mathcal{Y}_{R_0,r,\rho}$ with associated $X_j \in J_{R_0,r,\rho}(X_i)$ and witness $\Omega \in \vor^{(j)}_{R_0}(X_i)
	\cap
	\Bopen_{{T}_j}(\pi_{{T}_j}(X_i-X_j),\rho)^c$,
\begin{align*}
\Vert \eta_{\pi_{\partial M}(X_i)} - \tilde{\eta}_i^{(j)} \Vert 
& \leq
4
\theta
+
8
\sqrt{
\frac{\tau_{\min} \wedge \tau_{\partial, \min}}{\rho \wedge r_0}}
\frac{r}{\tau_{\min} \wedge \tau_{\partial, \min}}
,
\end{align*}
where $\tilde{\eta}_i^{(j)}=
\left ( \Omega-{\pi}_{T_j}(X_i-X_j) \right )/{\norm{\Omega-{\pi}_{T_j}(X_i-X_j)}}$.
\end{enumerate}

\end{enumerate}

\end{theorem}

A proof of \Cref{thm:deterministic} is given in the following \Cref{sec:proof_thm_deterministic}. 
\Cref{fig:parameter-explanation} below illustrates the role of the different parameters involved in this result. 
In the assumptions, \Cref{item:thm:deterministic:sampling,item:thm:deterministic:tangents} require that $M$ is sampled densely enough, and that the tangent spaces at sample points have been estimated with precision $\theta$. 
In light of \Cref{prop:tangent_space_estimation}, these assumptions will be satisfied at scales $\eps_1, \theta = O((\log n /n)^{1/d})$ for a random $n$-sample and a standard tangent space estimator, with high probability.
Then, \Cref{item:thm:deterministic:normals} basically requires that the ``curved rectangles'' $\{u \in \B(x,\sqrt{\varepsilon_2}) \cap M \mid \dd(u,\partial M) \leq \varepsilon_2\}$ nearby all $x \in \partial M$ are occupied by sample points (see \Cref{fig:parameter-explanation}). This assumption is key for identifying boundary observations and estimating normals accurately.
The volume heuristic given below \Cref{thm:detection} suggests that for a random $n$-sample, this assumption is satisfied with high probability at scale $\eps_2 = O((\log n / n)^{2/(d+1)})$.

\begin{figure}[!ht]
\centering
\includegraphics[width=0.4\textwidth]{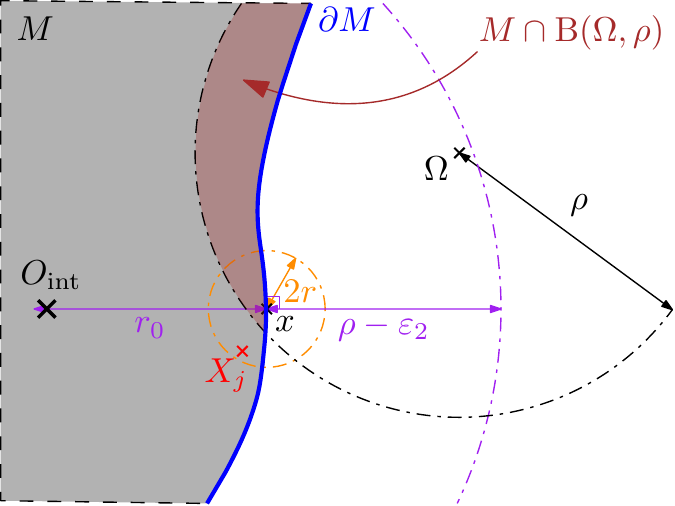}
\caption{
Illustration of Item 3 of \Cref{thm:deterministic} in full dimension $d=D=2$, yielding the simplification that $\pi_{T_j} = \mathrm{Id}_{d}$ for all $j \in \{1,\ldots,n\}$.
Here, we denoted $O_{\mathrm{int}} := \pi_{T_j}\left( x-X_j \right )-r_0\eta_j^*(x)$.
The assumption $\B(\Omega,\rho)\cap \mathcal{Y}_{j}\neq \emptyset$ means that the (brown) zone $\B(\Omega,\rho)\cap M$ contains sample points.
In full generality, when $d < D$, a similar layout can be drawn, with $M,\partial M$ and $M \cap \B(\Omega,\rho)$ replaced by $\pi_{T_j}(M),\pi_{T_j}(\partial M)$ and $\pi_{T_j}(M) \cap \B(\Omega,\rho)$.
}
\label{fig:parameter-explanation}
\end{figure}
The precise statement ensuring that the conditions of \Cref{thm:deterministic} are fulfilled with high probability for a random $n$-sample goes as follows.
\begin{restatable}{proposition}{propprobaconditioncat}\label{prop:proba_condition_cat}
Fix $R_0 \leq \frac{\tau_{\min}\wedge \tau_{\partial, \min}}{40} $, define $\rho_- = r_0 = \frac{R_0}{4}$, $\rho_+ = \frac{R_0}{2}$, and set
\[
\varepsilon_2 = r_0 \left (C_d \frac{f_{\max}^5}{f_{\min}^5} \frac{\log n}{f_{\min}(n-1)r_0^d} \right )^{\frac{2}{d+1}}.
\]
Then for $n$ large enough, the following statements hold with probability larger than $1- 3 n^{-\frac{2}{d}}$: for all $i \in \set{1,\ldots,n}$,
\begin{enumerate}[label=(\roman*),leftmargin=*]
\item $\displaystyle \angle(T_{X_i}M,\hat{T}_i) 
\leq 
\frac{1}{\tau_{\min}}
\left ( C_d \frac{f_{\max}^{4+d}}{f_{\min}^{5+d}} \frac{\log n}{n} \right )^\frac{1}{d} \leq 1/24$;
\item for all $(x,\Omega) \in (\B(X_i, r_0) \cap \partial M) \times \hat{T}_i$, 
 \[
 \| \Omega - (\pi_{\hat{T}_i}\left (x - X_i \right ) - r_0 \eta_i^*(x))\| \leq r_0 + \rho - \varepsilon_2 \quad \Rightarrow \quad \B(\Omega,\rho) \cap \mathcal{Y}_i \neq \emptyset,
 \]
 where $\eta_i^*(x)$ denotes the unique unit vector of $Nor(x,M) \cap \hat{T}_i$ (see \Cref{prop:existencesestimnorM}) .
 \end{enumerate}
\end{restatable}
A proof of \Cref{prop:proba_condition_cat} is given in \Cref{sec:proof_prop_proba_condition_cat}. 

\begin{proof}[Proof of \Cref{thm:detection}]
Combining \Cref{prop:proba_condition_cat} and \Cref{lem:covering} ensures that the requirements of \Cref{thm:deterministic} are fulfilled, with probability larger than $1-4n^{-2/d}$ for $n$ large enough, by choosing $T_i = \hat{T_i}$ and the following set-up:
$$
R_0 \leq \frac{\tau_{\min}\wedge \tau_{\partial, \min}}{40} \quad ,  \quad \frac{R_0}{2} = \rho_+ \geq \rho \geq \rho_- = r_0 = \frac{R_0}{4},$$
$$ \varepsilon_1 = \left ( C_d \frac{\log n}{f_{\min}n} \right ) ^\frac{1}{d} \quad ,  \quad \varepsilon_2 = r_0 \left (C_d \frac{f_{\max}^5}{f_{\min}^5} \frac{\log n}{f_{\min}n r_0^d} \right )^{\frac{2}{d+1}}, 
$$
$$
 \sqrt{(\tau_{\min}\wedge \tau_{\partial, \min})\varepsilon_2} = r_- \leq r \leq r_+ = \frac{R_0}{12} \quad ,  \quad \theta = \frac{1}{\tau_{\min}}\left ( C_d \frac{f_{\max}^{4+d}}{f_{\min}^{5+d}} \frac{\log n}{n} \right )^\frac{1}{d} \leq \frac{r_-} {\tau_{\min} \wedge \tau_{\partial, \min}}.
\qedhere
$$
\end{proof}

\subsection{Proof of Theorem~\ref{thm:deterministic}}
\label{sec:proof_thm_deterministic}
We decompose the proof into three intermediate results. 
As a first step, we prove that the sample points witnessing for boundary observations --- i.e. points $X_j$ making $J_{R_0,r, \rho}(X_i)$ nonempty, see \eqref{defboundObs} ---, must be close to $\partial M$. 
In fact, we show that they must be among the points $X_j$'s on which the Assumption~\ref{item:thm:deterministic:normals} of \Cref{thm:deterministic} holds.
See \Cref{sec:proof-of-lemdeterministicprotection} for the proof.

\begin{restatable}{lemma}{lemdeterministicprotection}
\label{lem:deterministic_protection}
Under the assumptions of \Cref{thm:deterministic}, if $X_j \in J_{R_0,r,\rho}(X_i)$, then $\partial M \neq \emptyset$ and 
\begin{equation*}
\dd(X_j,\partial M)
\leq 
2r
.
\end{equation*}
\end{restatable}

The next step builds upon \Cref{lem:deterministic_protection}, to guarantee that the detected boundary observations --- i.e. points $X_i$ such that $J_{R_0,r, \rho}(X_i) \neq \emptyset$ --- are close to the boundary $\partial M$, and that the associated estimated normals are close to the true normals at boundary points. 
In other words, we prove \Cref{thm:deterministic}~\ref{deter:protec:distance} and~\ref{deter:protec:angle}.
\begin{lemma}[{\Cref{thm:deterministic}~\ref{deter:protec:distance} and~\ref{deter:protec:angle}}]
\label{lem:deterministic_protection_2}
Under the assumptions of \Cref{thm:deterministic},  for all $X_i \in \mathcal{Y}_{R_0,r,\rho}$,
\begin{align*}
\dd(X_i,\partial M)\leq 2 \eps_2,
\end{align*}
and for all witness $X_j \in J_{R_0,r,\rho}(X_i)$,
\begin{align*}
\Vert \eta_{\pi_{\partial M}(X_i)} - \tilde{\eta}_i^{(j)} \Vert 
& \leq
4
\left(
\theta
+
\sqrt{  \left ( \frac{1}{\rho} + \frac{1}{r_0} \right ) \eps_2}
+
\frac{4r}{\tau_{\min}}
\right)
.
\end{align*}
\end{lemma}

\begin{proof}[Proof of \Cref{lem:deterministic_protection_2}]
To begin with, note that as $X_i \in \mathcal{Y}_{R_0,r,\rho}$ has some witness $X_j \in J_{R_0,r,\rho}(X_i)$, \Cref{lem:deterministic_protection} entails that $\partial M \neq \emptyset$. 
Also, since $\norm{X_i-X_j}\leq r \leq \tau_{\min}/48$, \Cref{prop:tangent_variation_geodesic} and \Cref{lem:dist_geod_dist_eucl} yield that 
\begin{equation}\label{deter:protec:eq1}
\angle(T_{X_i}M,T_j)\leq \theta + \frac{2r}{\tau_{\min}}\leq \frac{1}{24}+\frac{1}{24}\leq \frac{1}{12}.
\end{equation}
Furthermore, \Cref{lem:deterministic_protection} and triangle inequality gives
\begin{equation*}
\dd(X_i,\partial M) 
\leq 
\norm{X_i - X_j}
+
\dd(X_j,\partial M)
\leq
3r
,
\end{equation*}
so that $x':= \pi_{\partial M}(X_i) \in \partial M$ satisfies $\norm{x'-X_i} \leq 3r \leq R_0$.
\begin{figure}[ht!]
\centering
\includegraphics[width=0.8\textwidth]{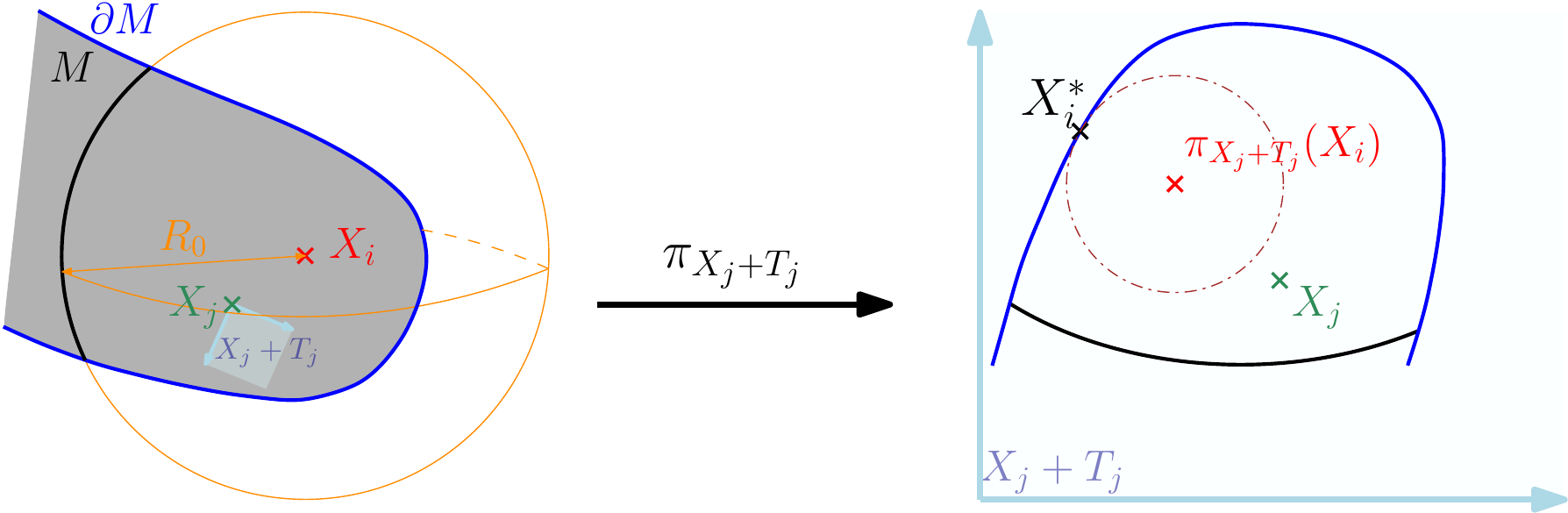}
\caption{
Layout for the proof of \Cref{lem:deterministic_protection_2}.
}
\label{fig:Xast-lift}
\end{figure}

Consider $X_i^*\in\argmin_{z \in \pi_{X_j+T_j}(\partial M \cap \B(X_i,R_0))}\Vert z-\pi_{X_j+T_j}(X_i)\Vert$ (see \Cref{fig:Xast-lift}).
As $x'\in \partial M\cap \B(X_i,R_0)$, $\pi_{X_j+T_j}(x')$ lies in the set where the $\argmin$ defining $X_i^*$ ranges, and hence
 \begin{equation}\label{eq:deter:protec:eq2}
 \norm{X_i^*-\pi_{X_j+T_j}(X_i)}
 \leq
 \norm{\pi_{X_j+T_j}(x')-\pi_{X_j+T_j}(X_i)}
 \leq
 \norm{x'-X_i}
 \leq
 3r
 .
\end{equation}
Introduce now $x\in \partial M \cap \B(X_i,R_0)$ such that $\pi_{T_j}(x-X_j)=X_i^*$.
From $\norm{x-X_i} \leq R_0$ only, \Cref{prop:projectionsestimtangentM},  
\eqref{deter:protec:eq1} and \eqref{eq:deter:protec:eq2} actually guarantee that
\begin{align*}
\norm{x-X_i}
\leq 
\frac{
\norm{X_i^*-\pi_{X_j+T_j}(X_i)}
}
{
1-\frac{1}{12}-\frac{\norm{x-X_i}}{2\tau_{\min}}
}
\leq
4r
.
\end{align*} 
Applying \Cref{prop:tangent_variation_geodesic} and \Cref{prop:geodesic_vs_euclidean} yields that 
\begin{align}\label{eq:deter:protec:eq1angle2}
\angle(T_{x}M,T_{j})
&\leq
 \angle(T_{x}M,T_{X_i}M)
 +
 \angle(T_{X_i}M,T_{X_j}M) 
 +
 \angle(T_{X_i}M,T_j) 
 \notag
 \\
 &\leq 
 \frac{2 \norm{X_i-x}}{\tau_{\min}}
 +
 \frac{2 \norm{X_j-X_i}}{\tau_{\min}}
 +
 \theta
 \notag
 \\
 &\leq
 \frac{10r}{\tau_{\min}}
 +
 \theta
 .
\end{align}
In particular, $\angle(T_{x}M,T_{j})\leq 1/8$, so that \Cref{cor:projorthobord} asserts that
\begin{equation*}
\pi_{X_j+T_j}(X_i)-X_i^*=-\norm{X_i^*-\pi_{X_j+T_j}(X_i)}\eta_j^*(x),
\end{equation*}
where $\eta_j^*(x)$ is \emph{the} unit vector of $T_j \cap Nor(x,M)$ (see \Cref{prop:projectionsestimnorM}).
Now, we write $O := X_i^*-r_0\eta_j^*(x)$. Recall that by definition, since $X_j \in J_{R_0,r,\rho}(X_i)$, there exists $\Omega=\pi_{T_j}(X_i-X_j)+\rho \tilde{\eta}_i^{(j)} \in T_j$ such that
$\B(\Omega,\rho) \cap \mathcal{Y}_{j}=\emptyset$.

On one hand, since $\B(\Omega,\rho) \cap \mathcal{Y}_j = \emptyset$, Assumption~\ref{item:thm:deterministic:normals} of \Cref{thm:deterministic} implies that
$\norm{\Omega-O}\geq r_0+\rho-\eps_2.$
On the other hand, we can develop
\begin{align*}
&\norm{\Omega-O}
=\norm{\bigl(r_0-\norm{X_i^*-\pi_{X_j+T_j}(X_i)}\bigr)\eta_j^*(x)+\rho \tilde{\eta}_i^{(j)}}
\\
&=
\sqrt{\bigl(\rho+r_0-\norm{X_i^*-\pi_{X_j+T_j}(X_i)}\bigr)^2-2\rho\bigl(r_0-\norm{X_i^*-\pi_{X_j+T_j}(X_i)})(1-\langle \eta_j^*(x),\tilde{\eta}_i^{(j)}\rangle}\bigr)
\\
&\leq
\rho+r_0-\norm{X_i^*-\pi_{X_j+T_j}(X_i)}-\frac{\rho\bigl(r_0-\norm{X_i^*-\pi_{X_j+T_j}(X_i)}\bigr)(1-\langle \eta_j^*(x),\tilde{\eta}_i^{(j)}\rangle)}{\rho+r_0-\norm{X_i^*-\pi_{X_j+T_j}(X_i)}}.
\end{align*}
Hence, combining the two above bounds on $\norm{\Omega-O}$ solves to
\begin{equation}\label{eq:protec_deter_final}
\norm{X_i^*-\pi_{X_j+T_j}(X_i)}
+
\frac{\rho\bigl(r_0-\norm{X_i^*-\pi_{X_j+T_j}(X_i)}\bigr)(1-\langle \eta_j^*(x),\tilde{\eta}_i^{(j)}\rangle)}{\rho+r_0-\norm{X_i^*-\pi_{X_j+T_j}(X_i)}}
\leq 
\eps_2.
\end{equation}
From \Cref{eq:protec_deter_final}, we can now conclude readily.
\begin{itemize}[leftmargin=*]
\item
To bound $\dd(X_i,\partial M)$, note that \eqref{eq:protec_deter_final} gives $\norm{X_i^*-\pi_{X_j+T_j}(X_i)}\leq \eps_2$. 
Therefore, \Cref{prop:projectionsestimtangentM} yields \Cref{thm:deterministic}~\ref{deter:protec:distance} by writing
\begin{align*}
\dd(X_i,\partial M)
\leq 
\norm{X_i-x}
&\leq
\frac{\norm{X_i^*-\pi_{X_j+T_j}(X_i)}}
{
1 - \angle(T_x M,T_j) - \norm{X_i-x}/(2\tau_{\min})
}
\leq  2\eps_2
.
\end{align*}
\item
To bound $\Vert\eta_{\pi_{\partial M}(X_i)} - \tilde{\eta}_i^{(j)} \Vert$, note that \eqref{eq:protec_deter_final} and the fact that $r_0 \geq 2 \eps_2$ also yield
$$
1-\langle\eta_j^*(x),\tilde{\eta}_i^{(j)}\rangle
\leq 
\frac{\rho+r_0}{\rho(r_0-\eps_2)}\eps_2 
\leq 2 \left (\frac{1}{\rho} + \frac{1}{r_0} \right ) \varepsilon_2.
$$
As $\eta_j^*(x)$ and $\tilde{\eta}_i^{(j)}$ are both unit vectors, this leads to
\begin{equation}\label{angle1}
\Vert\eta_j^*(x)-\tilde{\eta}_i^{(j)}\Vert
=
\sqrt{2(1-\langle\eta_j^*(x),\tilde{\eta}_i^{(j)}\rangle
)}
\leq
2
\sqrt{
\left (\frac{1}{\rho} + \frac{1}{r_0} \right ) \varepsilon_2
}
.
\end{equation}
In addition, \Cref{prop:estimeta} and bound \eqref{eq:deter:protec:eq1angle2} combine to
\begin{align}\label{angle2}
\Vert \eta_j^*(x) - \eta_{x} \Vert 
&\leq 
\sqrt{2}\angle (T_{x}M,T_j)
\leq
 \sqrt{2}
 \left(
 \frac{10r}{\tau_{\min}}
 +
 \theta
 \right)
.
\end{align}
Finally, triangle inequality yields
$$
\norm{x-\pi_{\partial M}(X_i)}\leq \norm{X_i-x}+ \dd (X_i,\partial M) \leq 4\eps_2 \leq (\tau_{\min} \wedge \tau_{\partial, \min})/32
,
$$
so that \Cref{prop:normal_vector_stability} asserts that
\begin{equation}\label{angle3}
 \Vert \eta_{\pi_{\partial M}(X_i)} - \eta_{x} \Vert \leq \frac{36}{\tau_{\min} \wedge \tau_{\partial, \min}} \varepsilon_2
\leq
2
\sqrt{
\left (\frac{1}{\rho} + \frac{1}{r_0} \right ) \varepsilon_2
}
.
\end{equation}
Combining  \Cref{angle1,angle2,angle3} with triangle inequality concludes the proof of \Cref{thm:deterministic}~\ref{deter:protec:angle} and that of \Cref{lem:deterministic_protection_2}.
\qedhere
\end{itemize}
\end{proof}
The last point \ref{deter:detec} of \Cref{thm:deterministic} derives from the following lemma.

\begin{lemma}[{\Cref{thm:deterministic}~\ref{deter:detec}}]\label{lem:deter_detection}
Under the assumptions of \Cref{thm:deterministic}, if $\partial M \neq \emptyset$, then for all $x\in \partial M$, there exists $X_i \in \mathcal{Y}_{R_0,r,\rho}$ such that
\begin{equation*}
\dd(x,\mathcal{Y}_{R_0,r,\rho})
\leq 
3r.
\end{equation*}
\end{lemma}

\begin{proof}[Proof of \Cref{lem:deter_detection}]
Let $x\in \partial M$, and assume without loss of generality that we have $\|x-X_1\| = \min_{1 \leq i \leq n} \| x - X_i\|$. We thus have $\norm{x-X_{1}}\leq \eps_1\leq R_0$. 
Similarly to the proof of \Cref{lem:deterministic_protection_2}, define
\[
X_1^*
\in
\argmin_{z \in \pi_{X_1+T_1}(\partial M \cap \B(X_1,R_0))}
\Vert z-X_1 \Vert
,
\] 
and take $y\in \partial M\cap \B(X_1,R_0)$ such that $\pi_{X_1+T_1}(y)=X_1^*$. As $x \in \partial M \cap \B(X_1,R_0)$, we have $\norm{X_1^*-X_1} \leq \norm{\pi_{X_1+T_1}(x-X_1)} \leq \norm{x-X_1} \leq \eps_1$, so that \Cref{prop:projectionsestimtangentM} entails
\[
\norm{y-X_1}\leq \frac{\eps_1}{1-\theta-\frac{R_0}{2\tau_{\min}}}\leq 2\eps_1.
\]
Since $\theta \leq 1/24$ and $\varepsilon_1 \leq \tau_{\min}/120$,  \Cref{prop:tangent_variation_geodesic,prop:geodesic_vs_euclidean} yield that 
\[
\angle(T_{y}M,T_{1}) 
\leq 
\angle(T_{y}M,T_{X_1} M) 
+
\angle(T_{X_1}M,T_{1}) 
\leq
\frac{2 \norm{X_{1}-y}}{\tau_{\min}}
+
\theta 
\leq 1/8.
\]
Hence, let $\eta_1^*(y)$ be \emph{the} unit vector of $Nor(y,M) \cap T_1$ (see \Cref{prop:projectionsestimnorM}). 
In turn, \Cref{lem:GeoInter_bleupart2} applied at $y$ asserts that
 \[
 \Bopen_{y+T_1}
 \bigl(
 y+2\rho_+ \eta_{1}^*(y)
 ,
 2\rho_+
 \bigr)
 \cap 
 \pi_{y+T_{1}}(\B(y,\tau_{\min}/16)\cap M)
 =
 \emptyset
 .
 \]
Since $R_0 \leq \tau_{\min}/32$, $\B(X_{1},R_0)\subset \B(y,\tau_{\min}/16)$. Moreover, $\pi_{X_1+T_{1}}(\B(X_{1},R_0)\cap M) = (X_1-y)^\perp + \pi_{y+T_{1}}(\B(X_{1},R_0)\cap M)$ and $(X_1-y)^\perp = (X_1^*-y)^\perp$, and hence
 \begin{align}				
 \label{eq:vacuity_2rho}
 \Bopen_{X_1+T_1}(X_1^*+2\rho_+ \eta_{1}^*(y),2\rho_+)\cap \pi_{X_1+T_{1}}(\B(X_{1},R_0)\cap M)=\emptyset.
 \end{align}
 Since $\rho\leq 2\rho_+$, we deduce that
  $
  \Bopen_{X_1+T_1}(X_1^*+\rho \eta_{1}^*(y),\rho)\cap (X_1+ \mathcal{Y}_{1})=\emptyset
  .
  $
Now, consider
$$
\delta
:=
\min
\set{
	t>0
	,
	\B_{X_1+T_1}\bigl(X_1^*+(\rho -t) \eta_{1}^*(y),\rho \bigr)
	\cap 
	(X_1 + \mathcal{Y}_{1})
	\neq 
	\emptyset
}
,
$$
  Since for all $t \geq \eps_2$,
  the point
  $
  \Omega_t 
  :=
  X_1^*-X_1+(\rho -t) \eta_{1}^*(y)
  \in T_1
  $
  satisfies
  $$
  \norm{\Omega_t- (\pi_{T_1}(y-X_1) -r_0\eta_{1}^*(y))}
  =
  r_0+\rho -t\leq r_0+\rho -\eps_2
  ,
  $$
  Assumption~\ref{item:thm:deterministic:normals} of \Cref{thm:deterministic} forces to have $\B(\Omega_t,\rho) \cap \mathcal{Y}_1 \neq \emptyset$, and hence $\delta \leq \eps_2$.
  
  By construction of $\delta$, there exists
  $z = \pi_{X_1+T_1}(X_{i_0}) \in \partial \B_{X_1+T_1}(X_1^*+(\rho - \delta ) \eta_{1}^*(y),\rho )\cap (X_1+\mathcal{Y}_{1})$.
  We may decompose $z$ as $z=X_1^*+\alpha  \eta_{1}^*(y) + \beta v$, where $v$ is a unit vector of $T_{1}\cap \Span(\eta_{1}^*(y))^{\perp}$.
  Since $z \in \partial \B_{X_1+T_1}(X_1^*+(\rho - \delta ) \eta_{1}^*(y),\rho )$ and $z \in (X_1+\mathcal{Y}_1) \subset \Bopen_{X_1+T_1}(X_1^*+2\rho_+ \eta_{1}^*(y),2\rho_+)^c$ from \eqref{eq:vacuity_2rho}, we have
  \begin{itemize}
  \item $\norm{z-(X_1^*+(\rho -\delta )\eta_{1}^*(y))}=\rho $, and thus $(\alpha-\rho + \delta)^2+\beta^2=\rho^2$;
\item $\norm{z-(X_1^*+2\rho_+\eta_{1}^*(y))}\geq 2\rho_+$, and thus $(\alpha-2\rho_+)^2+\beta^2 \geq 4\rho_+^2$.
  \end{itemize}
Therefore, after developing the above, we get that $\norm{X_1^*-z}^2 = \alpha^2 + \beta^2$ satisfies
 \[
 \begin{cases}
 \norm{X_1^*-z}^2
 =
 2\rho \delta-\delta^2
 +
 2\alpha(\rho -\delta)
 \leq
 2\rho \delta + 2 \rho_+ \alpha
 ,
 \\
  \norm{X_1^*-z}^2
  \geq
  4 \rho_+ \alpha
  =
  2 (2 \rho_+ \alpha)
  ,\\
 \end{cases}
 \]
which yields 
$
 \norm{X_1^*-z}^2
 \leq
 4\rho
 \delta
 \leq 
 4 \rho \eps_2.
$

Also by construction, we have
$
\Bopen_{T_1}(\Omega_\delta,\rho) 
\cap 
\mathcal{Y}_1
= 
\emptyset
$ 
and $\norm{
 \Omega_\delta - \pi_{X_1+T_1}(X_{i_0}-X_1)
 } 
 = 
 \norm{\Omega_\delta - z}
 =
 \rho
$
As a result, it is clear that if $\norm{X_{i_0}-X_{1}}\leq r$, then $\Omega_\delta \in \vor^{(1)}_{R_0}(X_{i_0})$, which yields $X_1 \in J_{R_0,\rho,r}(X_{i_0})$ and hence $X_i \in \mathcal{Y}_{R_0,r,\rho}$. 
Therefore, it remains to prove that $\norm{X_{i_0}-X_{1}}\leq r$ to conclude the proof.
For this, simply write 
$$
\norm{\pi_{X_1+T_{1}}(X_i)-X_{1}}\leq \| z - X_1^* \| + \| X_1^* - X_1 \| \leq  \left(\eps_1+2\sqrt{\rho \eps_2}\right)
,
$$
and since $X_{i_0}\in \B(X_1,R_0)$, \Cref{prop:projectionsestimtangentM} applied at $X_1$ yields 
$$\norm{X_i-X_{1}}\leq
2\left(\eps_1+2\sqrt{\rho \eps_2}\right) \leq r.$$
As a result, we can conclude the proof of \Cref{lem:deter_detection} (and \Cref{thm:deterministic}~\ref{deter:detec}) by noting that $$
\dd(x,\mathcal{Y}_{R_0,\rho,r})
\leq
\norm{x-X_{i_0}}
\leq
\norm{x-X_1}
+
\norm{X_1-X_{i_0}}
\leq 3r
.
\qedhere
$$
\end{proof}

\input{proof-probabilistic-layout}

\subsection{Proof of Theorem~\ref{thm:manifold-main-upper-bound}}\label{sec:proof_thm_manifold_upper_bound} 
The proof of \Cref{thm:manifold-main-upper-bound} is based on the following deterministic result, whose proof is deferred to \Cref{sec:appendix:linear-patches}.

\begin{restatable}[Estimation with Local Linear Patches]{theorem}{thmlocallinearpatchdeterministic}
\label{thm:local_linear_patch_deterministic} 
Write $r_0 := (\tau_{\min} \wedge \tau_{\partial, \min})/40$,
let $\eps_0,a,\delta \geq 0$, and $0 \leq \theta,\theta' \leq 1/16$. 
Assume that we have:
 \begin{enumerate}
  \item A point cloud $\mathcal{X}_n\subset M$ such that $\dHaus(M,\mathcal{X}_n)\leq \varepsilon_0$,
  \item Estimated tangent spaces $(T_i)_{1\leq i \leq n}$ such that $ \max_{1 \leq i \leq n} \angle (T_{X_i}M,T_i) \leq \theta$,
  \item A subset of boundary observations $\mathcal{X}_{\partial} \subset \mathcal{X}_n$ such that
 \[
 \max_{x\in \partial M} \dd(x,\mathcal{X}_{\partial})\leq \delta
 \text{ and }
 \max_{x\in \mathcal{X}_{\partial}} \dd(x,\partial M) \leq a \delta^2,
 \]
 from which we build interior observations
 $$
 \interior{\mathcal{X}}_{\varepsilon_{\partial M}} 
 := 
 \{ X_i \in \mathcal{X}_n \mid \dd(X_i, \mathcal{X}_{\partial}) \geq \varepsilon_{\partial M}/2 \}
 .
 $$
 \item Estimated unit normal vectors $(\eta_i)_{1\leq i \leq n}$ on $\mathcal{X}_{\partial}$ such that 
 $
 \max_{X_i \in \mathcal{X}_\partial} \Vert \eta_i  - \eta_{\pi_{\partial M}(X_i)}\Vert \leq \theta'.
 $
 \end{enumerate}
Let $\mathbb{M} = \mathbb{M}(\mathcal{X}_n,\mathcal{X}_\partial,T,\eta)$ be defined as
$
\mathbb{M}
:=
\mathbb{M}_{\Int}
\cup
\mathbb{M}_{\partial},
$
with
\begin{align*}
\mathbb{M}_{\Int}
&
:=
\bigcup_{X_i \in \interior{\mathcal{X}}_{\varepsilon_{\partial M}}} 
X_i + \B_{{T}_i}(0,\varepsilon_{\interior{M}}),
\\
\mathbb{M}_{\partial}
&
:=
\bigcup_{ X_i \in \mathcal{X}_{\partial}}
\left( 
X_i + \B_{{T}_i}(0, \varepsilon_{\partial M}) 
\right) 
\cap 
\{z,\langle z-X_i,{\eta}_i \rangle \leq 0 \}
,
\end{align*} 
Then if $\eps_{\partial M} \leq r_0/2$,
$\eps_0 \leq \eps_{\mathring{M}} \leq \eps_{\partial M}/6$
,
and $\max\set{\delta,a\delta^2} \leq \eps_{\partial M}/6$, we have
\begin{align*}
	\dHaus\bigl(M,\mathbb{M}\bigr)
	\leq 
	\begin{cases}	
	\varepsilon_{\mathring{M}}\left(\theta +\varepsilon_{\mathring{M}}/\tau_{\min}\right)
	&
	\text{if } \partial M = \emptyset,
	\\
  2a\delta^2
  +
   8\eps_{\partial M}
\left(
\theta + \theta' + \eps_{\partial M}/r_0
\right)
	&
	\text{if } \partial M \neq \emptyset.
	\end{cases}
\end{align*}
\end{restatable}

Equipped with \Cref{thm:local_linear_patch_deterministic}, choose, for $i \in \set{1,\ldots,n}$, $T_i = \hat{T}_i$ as in \Cref{prop:tangent_space_estimation},  ${\eta}_i = \tilde{\eta}_i$ as in \Cref{thm:detection}, and $\mathcal{X}_\partial = \mathcal{Y}_{R_0,r,\rho}$. Then  we define 
\[
\hat{M} := \mathbb{M} (\X_n,\mathcal{Y}_{R_0,r,\rho}, \hat{T}, \tilde{\eta}).
\]
Combining \Cref{prop:tangent_space_estimation}, \Cref{cor:boundary_tangent_space_estimation}, \Cref{thm:detection} and \Cref{lem:covering} ensure that the requirements of \Cref{thm:local_linear_patch_deterministic} are satisfied with probability at least $1-4n^{-\frac{2}{d}}$ for $n$ large enough, with the following choices of parameters: $\varepsilon_{\partial M}= 6\delta$,
\[
\begin{array}{lll}
\delta 
=
3r,
&
\varepsilon_0 
=
\left ( C_d \frac{\log n}{f_{\min}n} \right )^{\frac{1}{d}},
&
\varepsilon_{\interior{M}}
=
\left ( C_d \frac{\log n}{f_{\min}n} \right )^{\frac{1}{d}}, 
\\
\theta
=
\left ( C_d \frac{f_{\max}^{4+d}}{f_{\min}^{5+d}} \frac{\log n }{(n-1)\tau_{\min}^d} \right )^{\frac{1}{d}}, 
&
\theta'
=
\frac{20r}{\sqrt{(\tau_{\min}\wedge \tau_{\partial, \min})R_0}},
&
a
=
(4( \tau_{\min}\wedge \tau_{\partial, \min}))^{-1},
\end{array}
\]
which concludes the proof of the first bound in \Cref{thm:manifold-main-upper-bound}. 

To get the bound in expectation, let $K$ denote the diameter of $M$, and note that there exists $X_{i_0} \in \X_n$ such that $\{X_{i_0}\} \subset \hat{M}$, so that 
$
\sup_{x \in M} \dd(x,\hat{M}) \leq K, 
$
almost surely. 
Conversely, since $\hat{M} \subset M + \B(0, \varepsilon_{\partial M} \vee \varepsilon_{\interior{M}})$, we deduce that
$
\sup_{x \in \hat{M}} \dd(x,M) \leq K$
for $n$ large enough. 
Finally, noticing that for $n$ large enough, the result follows by writing
\[
\bigl( 4n^{-\frac{2}{d}} \bigr) K 
\leq 
C_d
(\tau_{\min}\wedge \tau_{\partial, \min})
\left [\left ( \dfrac{f_{\max}^{2+d/2}}{f_{\min}^{2+d/2}}
\dfrac{\log n}{f_{\min} \tau_{\min}^d n} \right )^{\frac{2}{d}}
\wedge 
\left ( \dfrac{f_{\max}^5 }{f_{\min}^5} \dfrac{\log n}{f_{\min} (\tau_{\min} \wedge \tau_{\partial, \min})^d n} \right )^{\frac{2}{d+1}} \right ]
.
\]

%% file: proof-probabilistic-layout.tex
\subsection{Proof of Proposition~\ref{prop:proba_condition_cat}}\label{sec:proof_prop_proba_condition_cat}

\begin{proof}[Proof of \Cref{prop:proba_condition_cat}]
Without loss of generality we fix $i=1$, and work conditionally on $X_1$. Let $A_1$ denote the event 
\[
A_1: = \left \{ \angle(T_{X_1}M,\hat{T}_1) \leq C_d \left ( \frac{f_{\max}^{4+d}\log n}{f_{\min}^{5+d} \tau_{\min}^d(n-1)}\right)^{1/d} \right \}, 
\]
which has probability larger than $1-2 \left ( 1/n \right )^{1 + \frac{2}{d}}$ from \Cref{prop:tangent_space_estimation}. 
Note that $A_1$ is $\sigma(Y_2, \hdots, Y_n)$-measurable, where $Y_i = X_i \indicator{X_i \in \B(X_1,h)}$.
We further assume that $n$ is large enough so that we have $\angle(T_{X_1}M,\hat{T}_1) \leq 1/24$ on $A_1$. In particular, note that Item (i) is satisfied on $A_1$. 

Let us now bound the probability that Item (ii) does not occur. As in \Cref{lem:proba_occupation_losange}, we assume $\varepsilon_2 := \left ( A \frac{f_{\max}^4 \log n}{f_{\min}^5 (n-1)} \right )^{\frac{2}{d+1}}$, where $A$  is to be fixed later. For  $x \in \B(X_1, r_0) \cap \partial M$, denote by $O_x^{int} = \pi_{\hat{T}_1}(x-X_1) - r_0 \eta_1^*(x)$. 

Recall here that $\mathcal{Y}_1$ is defined by $\mathcal{Y}_{1}:=\pi_{T_1}(\B(X_1,R_0)\cap \mathcal{X}_n - X_1)$.
If $\Omega \in \hat{T}_1$ is such that $\B(\Omega,\rho) \cap \mathcal{Y}_1 = \emptyset$ and $\|\Omega - O_x ^{int}\| \leq \rho + r_0 - \varepsilon_2$ for some $\rho \geq \rho_-$ and $\rho_- + r_0 > \varepsilon_2 >0$, then choosing $\Omega_0 = \Omega + (\rho - \rho_-) \frac{O_{x}^{int} - \Omega}{\rho + r_0 - \varepsilon_2 }$ yields that
\begin{align*}
\left\{
\begin{array}{l}
\B(\Omega_0, \rho_-) \cap \mathcal{Y}_1 
\subset 
\B(\Omega, \rho) \cap \mathcal{Y}_1 = \emptyset, 
\\
\| \Omega _0 - O_x^{int}\| \leq r_0 + \rho_- - \varepsilon_2.
\end{array}
\right.
\end{align*} 
But as $\|x-X_1\| \leq r_0$, \Cref{lem:GeoInter_bleupart2} ensures that on the event $A_1$ we have
\begin{align*}
\B(O_x^{in},r_0) \cap \hat{T}_1 \subset \pi_{\hat{T}_1}(\B(X_1,5r_0/2 + r_0) \cap M - X_1) \subset \pi_{\hat{T}_1}(\B(X_1,R_0) \cap M-X_1)
.
\end{align*}
Thus, if we let
\begin{multline*}
\mathcal{Q}_{r,\rho,\varepsilon}
:=
\left\{
(O,\Omega) \in \B_{\hat{T}_1}(0,2r_0) \times \B_{\hat{T}_1}(0,4r_0) 
\bigl|
\|\Omega-O\| \leq r + \rho - \varepsilon 
\right.
\\
\left.
\text{ and }
\B_{\hat{T}_1}(O,r) \subset \pi_{\hat{T}_1}(\B(X_1,R_0) \cap M - X_1)
\right\}
,
\end{multline*}
then for all $\rho \geq \rho_-$, we have the inclusion of events
\begin{align*}
& \bigl\{ \exists (x,\Omega) \in \B(X_1,r_0)\times \hat{T}_1 \mid \| \Omega - O_x^{int}\| \leq r_0 + \rho - \varepsilon_2 \text{ and } \B(\Omega,\rho) \cap \mathcal{Y}_1 = \emptyset 
\bigr\} \cap A_1 \\
&\quad \subset \bigcup_{(O,\Omega) \in \mathcal{Q}_{r_0,\rho_-,\varepsilon_2}} \{ \B(\Omega,\rho_-) \cap \mathcal{Y}_1 = \emptyset \} \cap A_1.
\end{align*}
This union of events being infinite, we now discretize space by considering an $(\varepsilon_2/8)$-covering $\mathcal{C}(\varepsilon_2)$ of $\B_{\hat{T}_1}(0,4r_0)$. For all $(\Omega,O) \in \mathcal{Q}_{r_0,\rho_-,\varepsilon_2}$, we also let $\Omega'$ and $O'$ denote the closest elements in $\mathcal{C}(\varepsilon_2)$ to $\Omega$ and $O$ respectively.
Letting $r'_0 := r_0 - \varepsilon_2/8$ and $\rho'_0 := \rho_-- \varepsilon_2/8$, triangle inequality yields that on $A_1$, 
\begin{align*}
\left\{
\begin{array}{l}
\B_{\hat{T}_1}(O',r'_0) \subset  \pi_{\hat{T}_1}(\B(X_1,R_0) \cap M)-X_1,
\\
\B(\Omega', \rho'_0) \cap \mathcal{Y}_1 = \emptyset, 
\\
\|\Omega' - O' \| \leq r'_0 + \rho'_0 - \varepsilon_2/2.
\end{array}
\right.
\end{align*}
As a result, provided that $n$ is large enough so that $\varepsilon_2 \leq 4 r_0$, the previous event union satisfies
\begin{align*}
\bigcup_{
\mathcal{Q}_{r_0,\rho_-,\varepsilon_2}} \{ \B(\Omega,\rho_-) \cap \mathcal{Y}_1 = \emptyset \} \cap A_1
&\subset 
\bigcup_{
\mathcal{Q}_{\frac{r_0}{2},\frac{\rho_-}{2},\frac{\varepsilon_2}{2}} \cap \mathcal{C}(\varepsilon_2)^2} \left \{ \B \left (\Omega,\frac{\rho_-}{2} \right ) \cap \mathcal{Y}_1 = \emptyset \right \} \cap A_1
.
\end{align*}
Let $(O,\Omega) \in \mathcal{P}\left ( \frac{r_0}{2},\frac{\rho_-}{2},\frac{\varepsilon_2}{2} \right)$ be now fixed. 
Recalling that $Y_i = X_i \indicator{X_i \in \B(X_1,h)}$, and that the event $A_1$ is $\sigma(Y_2, \hdots, Y_n)$-measurable, we may write
\begin{align*}
\mathbb{P} &\left ( A_1 \cap  \left \{ \B \left (\Omega,\frac{\rho_-}{2} \right ) \cap \mathcal{Y}_1 = \emptyset \right \} \right) 
\\
&= 
\mathbb{E} \left [ \mathbb{P} \left ( A_1 \cap  \left \{ \B \left (\Omega,\frac{\rho_-}{2} \right ) \cap \mathcal{Y}_1 = \emptyset \right \} \mid (Y_2, \hdots, Y_n) \right)  \right ]
\\
&= 
\mathbb{E} \left [ \indicator{A_1} \mathbb{P} \left (  \left \{ \B \left (\Omega,\frac{\rho_-}{2} \right ) \cap \mathcal{Y}_1 = \emptyset \right \} \mid (Y_2, \hdots, Y_n) \right) \right ] 
\\
&
\leq 
\mathbb{E} \left [ \indicator{A_1} \mathbb{P} \left ( 
\min_{ X_i \in \X_n \cap \B(X_1, R_0)} \bigl\|\pi_{\hat{T}_1}(X_i-X_1)-\Omega\bigr\| > \frac{\rho_-}{2}
\mid (Y_2, \hdots, Y_n)  \right) \right ]\\
&\leq 
\mathbb{E} \left [ \indicator{A_1} \mathbb{P} \left ( 
\min_{ X_i \in \X_n \cap (\B(X_1, R_0)\setminus \B(X_1,h))} \bigl\|\pi_{\hat{T}_1}(X_i-X_1)-\Omega\bigr\| > \frac{\rho_-}{2}
\mid (Y_2, \hdots, Y_n)  \right) \right ]
.
\end{align*}
Furthermore, as the family $(\pi_{\hat{T}_1}(X_i))_{X_i \notin \B(X_1,h)}$ is i.i.d conditionally on $(Y_2, \hdots, Y_n)$, \Cref{lem:proba_occupation_losange} yields
\begin{align*}
\mathbb{E} 
&
\left [ \indicator{A_1} \mathbb{P} \left ( 
\min_{ X_i \in \X_n \cap (\B(X_1, R_0)\setminus \B(X_1,h))} \bigl\|\pi_{\hat{T}_1}(X_i-X_1)-\Omega\bigr\| > \frac{\rho_-}{2}
\mid (Y_2, \hdots, Y_n)  \right) \right ] \\
&\leq 
\mathbb{E} \left [ \indicator{A_1} \left ( 1 - A r_0^{\frac{d-1}{2}} C_d  \frac{f_{\max}^4 \log n}{f_{\min}^4(n-1)} \right )^{n - |\mathbb{X}_n \cap \B(X_1,h)|}\right ] 
\\
&\leq 
\sum_{k=0}^{n-1} \binom{n-1}{k} (C_d f_{\max} h^d)^k \left ( 1 - A r_0^{\frac{d-1}{2}} C_d  \frac{f_{\max}^4 \log n}{f_{\min}^4(n-1)} \right )^{n - k} 
\\
&\leq 
\left ( 1 - A r_0^{\frac{d-1}{2}} C_d  \frac{f_{\max}^4 \log n}{f_{\min}^4(n-1)} + \frac{C_d f_{\max}^5 \log n}{f_{\min}^{5}(n-1)} \right )^{n-1}
.
\end{align*}
Choosing $A:= C_d \frac{f_{\max}}{f_{\min}}r_0^{\frac{1-d}{2}} \geq C_d r_0^{\frac{1-d}{2}}$, yields that
\begin{align*}
\left ( 1 - A r_0^{\frac{d-1}{2}} C_d  \frac{f_{\max}^4 \log n}{f_{\min}^4(n-1)} + \frac{C_d f_{\max}^5 \log n}{f_{\min}^{5}(n-1)} \right )^{n-1} \leq n^{-C_d},    
\end{align*}
so that, for $C_d$ large enough,
\begin{align*}
\left | \mathcal{C}(\varepsilon_2) \right |^2 \mathbb{P} \left ( A_1 \cap  \left \{ \B \left (\Omega,\frac{\rho_-}{2} \right ) \cap \mathcal{Y}_1 = \emptyset \right \} \right) \leq \left ( \frac{1}{n} \right )^{1 + \frac{2}{d}},
\end{align*}
for $n$ large enough. 
Thus, a union bound gives the result of \Cref{prop:proba_condition_cat}, since we have set $\varepsilon_2 = r_0 \left (C_d \frac{f_{\max}^5}{f_{\min}^5} \frac{\log n}{f_{\min}(n-1)r_0^d} \right )^{\frac{2}{d+1}}$ for $C_d$ large enough.

\end{proof}

%% file: simulations.tex
\section{Computational considerations and experimental illustrations}
\label{sec:experiments}

\subsection{Pseudo-code and computational complexity}
\label{sec:computational-complexity}

        \begin{algorithm}
            \begin{algorithmic}[1]
                \REQUIRE ~\\
                $\X_n = \{X_1,\ldots,X_n\} \subset \R^D$: Sample points \\
                $h$: Bandwidth for principal component analysis 
                \hfill \textit{\{write $N_{\mathrm{PCA}} = nh^d$\}}
                \\
                $R_0$: Macroscopic localization scale for projections
                \hfill \textit{\{write $N_{\mathrm{loc}} = nR_0^d$\}}
                \\
                $r$: Neighborhood radius of boundaryness witnesses
                \hfill \textit{\{write $N_{\mathrm{wit}} = nr^d$\}}
                \\
                $\rho$: Minimal width of witnessing Voronoi cells
                \\
                \ENSURE $h,r, R_0,\rho \geq 0$
                \STATE Initialization with an empty boundary structure
                \\
                $\mathcal{Y}_{R_0,r,\rho} \leftarrow \emptyset$
                \\
                $\mathbf{\eta} \leftarrow \emptyset$
                \STATE 
                Computation of pairwise point distances
                \\
                $\mathrm{dist} \leftarrow \bigl(\norm{X_j-X_i}\bigr)_{1 \leq i,j \leq n}$ 
                \hfill \textit{\{time $O(Dn^2)$\}}
                \FORALL{$j \in \{1,\ldots,n\}$}
                    \STATE 
                    Estimation of $T_{X_j}M$ via local PCA on a  $h$ neighborhood of $X_j$
                    \\
                    $\pi_{\hat{T}_j} \leftarrow \mathrm{PCA}_d(\X_n \cap \B(X_j,R_0))$
                    \hfill 
                    \textit{\{time $O(N_{\mathrm{PCA}}Dd)$ with fast-PCA~\cite{Sharma07}\}}
                    \STATE 
                    Projection of the $R_0$-neighborhood of $X_j$ onto $\hat{T}_j$
                    \\
                    $\hat{\X}_{n}^{(j)} \leftarrow \pi_{\hat{T}_j}\bigl(\X_n \cap \B(X_j,R_0)\bigr)$
                    \hfill
                    \textit{\{time $O(N_{\mathrm{loc}}Dd)$\}}
                    \STATE Computation of the Voronoi cells of the projected sample
                    \\
                    $\vor^{(j)}_{R_0} \leftarrow \bigl(\vor^{(j)}_{R_0}(X_i)\bigr)_{\hat{X}_i \in \hat{\X}_n^{(j)}}$
                    \hfill 
                    \textit{\{time $O(N_{\mathrm{loc}}2^d K)$~\cite{Sheehy15}\}}
                \FORALL{$\hat{X}_i^{(j)} \in \hat{\X}_n^{(j)} \cap \B(X_j,r)$}
                    \STATE
                    List the vertices $v_{i,1}, \ldots, v_{i,K}$ of the Voronoi cell $\vor^{(j)}_{R_0}(X_i)$ of $\hat{X}_i = \pi_{\hat{T}_j}(X_i)$ 
                    \\
                    \hfill \textit{\{Typical $K = O\bigl(d^{\frac{d}{2}-1}\bigr)$~\cite[Theorem 7.2, case $s = 0$]{Moller89}%
                    \}%
                    }
                    \\
                    \STATE
                    Computation of the radius and prominent direction of $\vor^{(j)}_{R_0}(X_i)$ when centered at $\pi_{\hat{T}_i}(X_j)$
                    \\
                    $\hat{v}_{i}^{(j)} 
                    \leftarrow 
                    \argmax_{1 \leq k \leq K}\|v_{i,k} - \pi_{\hat{T}_j}(X_i)\|
                    $
                    \\
                    $
                    \rho_{i}^{(j)} 
                    \leftarrow 
                    \|\hat{v}_{i}^{(j)} - \pi_{\hat{T}_j}(X_i)\|
                    $
                    \hfill \textit{\{time $O(dK)$\}}
                    \IF{$\rho_{i}^{(j)}>\rho$} 
                        \STATE Addition of $X_j$ to the boundary observations with associated unit normal
                        \\
                        $\mathcal{Y}_{R_0,r,\rho} \leftarrow \mathcal{Y}_{R_0,r,\rho} \cup \{X_j\}$
                        \\
                        $\mathbf{\eta}_j
                        \leftarrow
                        \frac{\hat{v}_{i}^{(j)} - \pi_{\hat{T}_j}(X_i)}{\|\hat{v}_{i}^{(j)} - \pi_{\hat{T}_j}(X_i)\|}
                        $
                    \ENDIF
                \ENDFOR \hfill \textit{\{$O(1+N_{\mathrm{wit}})$ iterations\}}
                \ENDFOR \hfill \textit{\{$n$ iterations\}}
                \RETURN
                $(\mathcal{Y}_{R_0,r,\rho},\mathbf{\eta})$
            \end{algorithmic} 
            \caption{\BO \hfill \textit{\{Average time complexity\}}}
            \label{alg:boundary-labelling}
        \end{algorithm}

The algorithm \BO displays the pseudo-code for building the filtered set $\mathcal{Y}_{R_0,r,\rho}$ from \Cref{def:boundary_points}, alongside with the estimated unit normals from \Cref{defi:normal_vector_estimate}.
Note that it only uses standard algorithmic sub-blocks, such as PCA and Voronoi diagrams.

Its time complexity critically depends on the number of vertices $K$ of the individual $d$-dimensional Voronoi cells.
In worst case, each vertex corresponds to a $d$-simplex in the Delaunay triangulation, so that $K = O(n^{\lceil \frac{d}{2} \rceil})$, which corresponds to the maximum number of simplices in a Delaunay triangulation based on $n$ points).
Yet, for random point clouds, the typical number of such vertices reduces to $K = O\bigl(d^{\frac{d}{2}-1}\bigr)$~\cite[Theorem 7.2, case $s = 0$]{Moller89}.

At the end of the day, the average time complexity of \BO does not exceed (up to logarithmic factors)
$$
O\left(
Dn^2
+
n
\biggl(
    N_{\mathrm{PCA}}Dd
    +
    N_{\mathrm{loc}}Dd
    +
    N_{\mathrm{loc}} 2^d K
    +
    (1+N_{\mathrm{wit}})
    d K
\biggr)    
\right)
,
$$
where $N_{\mathrm{PCA}} = nh^d$, $N_{\mathrm{loc}} = nR_0^d$, and $N_{\mathrm{wit}} = n r^d$.

The parameters leading to optimal rates in our theoretical results (see \Cref{prop:tangent_space_estimation} and \Cref{thm:detection}) correspond to $h \asymp (1/n)^{1/d}$, $R_0 \asymp \rho \asymp 1$, and $r \asymp (1/n)^{1/(d+1)}$ up to $\log n$ factors.
Hence, the above time complexity bound boils down to 
\begin{align*}
O\left(
n^2
\bigl(
    Dd
    +
    d^{\frac{d}{2}}
\bigr)    
\right)
.
\end{align*}

Overall, note that the dependency on the ambient dimension $D$ is limited to a linear factor.
On the other hand, the leading factor in terms of sample size arises from the computation of the whole distance matrix of $\X_n$.
This dependency could be mitigated to $O(Dn N_{\mathrm{loc}})$ by only computing the local distances $\left( \norm{X_i-X_j} \right)_{\norm{X_i-X_j} \leq R_0}$. This can be done approximately, for instance via standard greedy exploratory geometric algorithms~\cite{Har11}.
The factor $O(n^2)$ factor is also contributed from the fact that $N_{\mathrm{loc}} \asymp n$, since $R_0 \asymp 1$. 
Strategies consisting in taking $R_0 = o(1)$ could lead to a computation-precision tradeoff (see \Cref{thm:boundary-main-upper-bound}).
Similarly, the choice $r = 0$ that we used in practice yield $N_{\mathrm{wit}} = 1$, further reducing algorithmic complexity.

\subsection{Heuristics for data-driven parameters calibration}
\label{sec:heuristics-parameters}

\subsubsection{Bandwidth $h$ for PCA}
\label{sec:heuristic_h}

Intuitively, the bandwidth $h$ from \Cref{prop:tangent_space_estimation} should be taken so that all the balls $\B(X_i,h)$ contain at least $d \vee \log n$ sample points.
In practice, we set 
$$
h(k) 
=
\inf \{r > 0 \mid \forall i \in \{1,\ldots,n\}, \left | \B(X_i,r) \cap \X_n \right |  \geq k \},
$$
with $k = d \log n$. 
This particular scale is chosen so that, on average, the PCA's are computed using $\log n$ neighbors per estimated direction.

Let us note that in the noise-free case considered here, the choice of $h$ (or $k$) does not impact significantly the tangent space estimation step provided $k \geq d+1$ and $h(k) = O(\dHaus(M,\X_n))$.
This latter Hausdorff distance can be approached by the smallest $h$ such that $\B(\X_n,h)$ is connected.
In the same spirit, pointwise choices of $h = h_i$ based on $k$-nearest neighbors could be a way to adapt to possible non-uniformity of sampling ($f_{\min} \ll f_{\max}$).

\subsubsection{Macroscopic localization scale $R_0$ for projection}
\label{sec:heuristic_R0}

Throughout the theoretical analysis of the method, the scale $R_0$ of \Cref{def:boundary_points} is chosen so that the (approximate) tangent projections $\pi_{\hat{T}_i} : M \cap \B(X_i,R_0) \to \hat{T}_i$ do not distort the metric significantly.
Assessing this information empirically may be performed graphically, via the scatter-plot of all pairs $$\bigl(\norm{X_i-X_j}, \norm{\hat{T}_{i}(X_i-X_j))}\bigr)_{1 \leq i,j \leq n}
.
$$
With this graphical representation in mind, a natural choice of $R_0$ is the largest radius $R$ such that the plot remains close to the diagonal $y=x$ over $[0,R]$, where ``close'' needs to be properly defined.

This heuristic is used in the experiments of \Cref{sec:simulations}.
As metric distortion is a multiplicative quantity, we consider
$$
\gamma(R)
:=
\min_{\|X_i-X_j\|\leq R}
\frac{\norm{\hat{T}_{i}(X_i-X_j))}}{
\norm{X_i-X_j}}
,
$$
and we choose the largest $R$ such that $\gamma(R)$ remains close to $1$.
That is, we pick
$$
R_0
=
\max
\left\{
R>0
\left|
\forall X_i,X_j \in \X_n \text{~s.t.~} \norm{X_i - X_j} \leq R,
\left|
\frac{\norm{\hat{T}_{i}(X_i-X_j))}}{\norm{X_i - X_j}}
-1
\right|
\leq
\delta
\right.
\right\}
,
$$
where $\delta<1$ is a metric distortion tolerance parameter.
We believe that any choice of $\delta \leq 1/2$ would have the method work.

\subsubsection{Neighborhood radius $r$ of boundaryness witnesses}
\label{sec:heuristic_r}

We strongly believe that parameter $r$ introduced in \Cref{def:boundary_points} is a purely theoretical artifact.
That is, the overall method is likely to have the same theoretical guarantees when applied with $r=0$. 
We did not succeed in proving this conjecture, apart for $d=1$, in which case $\partial M$ consists of only two points per connected component of $M$.

On the practical side, we conducted our experiments with $r=0$, which yielded satisfactory results. 
Recall that the smaller $r$, the smaller the set of detected boundary points $\mathcal{Y}_{R_0,r,\rho}$, which could lead to potentially too many false negative (i.e. nearby-boundary points missed).

Yet, beyond the present idealized framework where points are not corrupted with noise, this parameter may have an influence.
If so, a possible calibration strategy could consist in investigating the size of the Voronoi cell of $\pi_{\hat{T}_j}(X_i-X_j)$, where $X_j$ ranges among the $k$-nearest neighbors of $X_i$ for growing $k$'s, and to stop when the number of detected points stabilizes. 

\subsubsection{Minimal width $\rho$ of witnessing Voronoi cells}
\label{sec:heuristic_rho}

We finally move to discussion about $\rho$, the last parameter involved in \Cref{def:boundary_points}.
Recall that for $r=0$, the labelling of $X_i$ is only based on the size of the Voronoi cell of $\pi_{\hat{T}_i}(X_i)$. 
Hence, for a given point $X_i$, we are led to compute 
$$
\rho_i = \max_{1 \leq k \leq  K}\|\pi_{\hat{T}_i}(X_i) - v_{i,k}\|
,$$
where $(v_{i,k})_{k \leq K}$ are the vertices of the Voronoi cell (see \BO). 
For interior points, $\rho_i$ is expected to be small, while for points close to the boundary, $\rho_i$ is expected to be larger.
The value of $\rho$ effectively fixes the chosen cutoff between ``small'' and `large'' cells.
To pick $\rho$ wisely, we investigate the distribution of the values $(\rho_i)_{1 \leq i \leq n}$.

That is, we reorder values $\rho_{(1)},\ldots,\rho_{(n)}$, and we plot the graph $(i,\rho_{(i)})_{1 \leq i \leq n}$. 
This graph typically exhibits a sharp jump (see~\Cref{fig:simu1,fig:simu2}).
As the value of this jump corresponds to a phase transition between the two regimes we want to distinguish, we select $\rho$ to be the mid-value of this first jump.
Note that for $r > 0$, a similar strategy based on the $\rho_i^{(j)}$'s (as defined in Section \ref{sec:computational-complexity}) can easily be built.

\subsection{Simulations on low-dimensional examples}
\label{sec:simulations}

We now illustrate the boundary detection, normal vector estimation and parameter tuning heuristics on some examples on four low-dimensional examples.
Namely, the toy distributions that we consider consist of the uniform distributions over the following sets.

\begin{enumerate}
    \item[($d=1$, $D=3$)]
    The spiral given by the parametrization
    $
    [0,5\pi] \ni \theta \mapsto (\cos(\theta), \sin(\theta), \theta/3)
    .
    $

    \item[($d=2$, $D=2$)]
    The annulus $\B(0,1)\setminus \B(0,0.4)$.
    \item[($d=2$, $D=3$)] 
    The unit half sphere $\{(x,y,z) |~x^2+y^2+z^2 = 1 \text{~and~} x \geq 0\}$.
    \item[($d=2$, $D=3$)] 
    The M\"obius strip given by the parametrization 
    $$
    [-1,1]\times[0,2\pi] \ni (u,\theta) \mapsto
    ((u\cos(\theta/2)+3)\cos(\theta),(u\cos(\theta/2)+3)\sin(\theta),u\sin(\theta/2))
    .
    $$ 
\end{enumerate}
For each of these four distributions, \Cref{fig:simu1,fig:simu2} present:
\begin{itemize}
    \item 
    The metric distortion scatterplot used for the calibration of $R_0$ (see \Cref{sec:heuristic_R0});
    \item 
    The order distance histogram plot used for the calibration of $\rho$ (see \Cref{sec:heuristic_rho});
    \item
    Points clouds of $\{500,1000,2000,5000\} \ni n$-samples, with the detected boundary observations and their associated estimated normals displayed in red.
\end{itemize}
In all the cases, bandwidth $h$ is chosen using the rule of \Cref{sec:heuristic_h}, and $r=0$ (see \Cref{sec:heuristic_r}).
These plots are for illustrative purpose only, and are not meant to illustrate minimax convergence rates.
Qualitatively, let us point the following:
\begin{itemize}
    \item 
    For the spiral, exactly two observations are labelled as boundary observations with associated Voronoi cells that are unbounded. This advocates that for possibly setting $\rho=\infty$ in the one-dimensional case $d=1$.
    \item 
    For the annulus, no tangent projection is performed, since we are on a full-dimensional domain ($d=D$).
    This is why the scatterplot of 
    $\bigl(\norm{X_i-X_j}, \norm{\hat{T}_{i}(X_i-X_j))}\bigr)_{1 \leq i,j \leq n}$ coincides with the identity.
    This advocates for setting $R_0=\infty$ in the full-dimensional case $d=D$. 
    Note that if so, only one global Voronoi diagram (that of the complete sample points $\X_n$) needs to be computed, as opposed to one local Voronoi per point.
\end{itemize}

\setlength\mgheight{5em}
\begin{table}
\centering
\begin{tabular}{c m{0.23\textwidth} m{0.23\textwidth}} 
\toprule 
\multicolumn{1}{c}{$M$} & 
\multicolumn{1}{c}{Spiral} & 
\multicolumn{1}{c}{Annulus} \\ 
\midrule
Calibration of $R_0$ &    
\includegraphics[width=\linewidth]{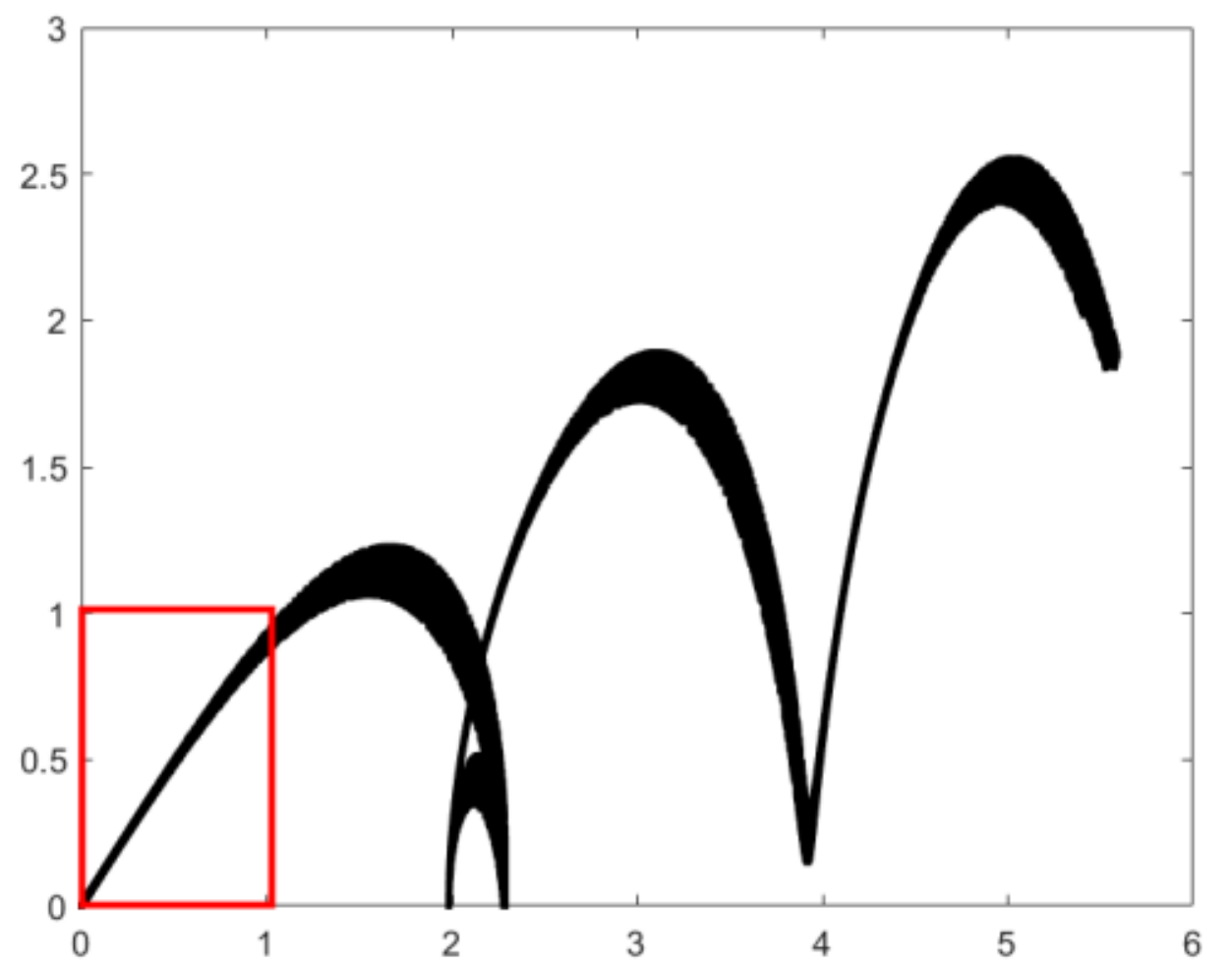}
 & 
\includegraphics[width=\linewidth]{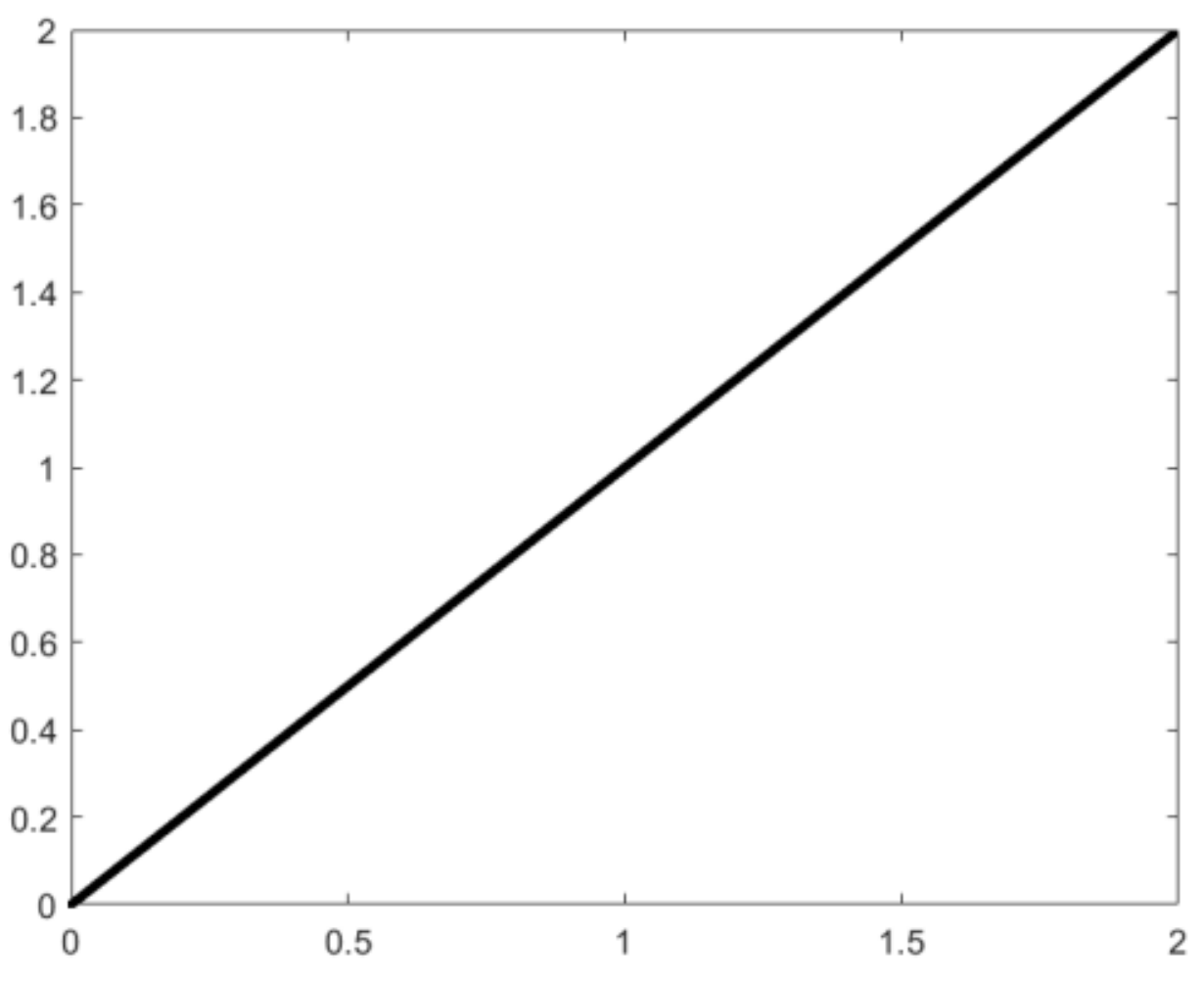}
\\
Calibration of $\rho$ &    
\includegraphics[width=\linewidth]{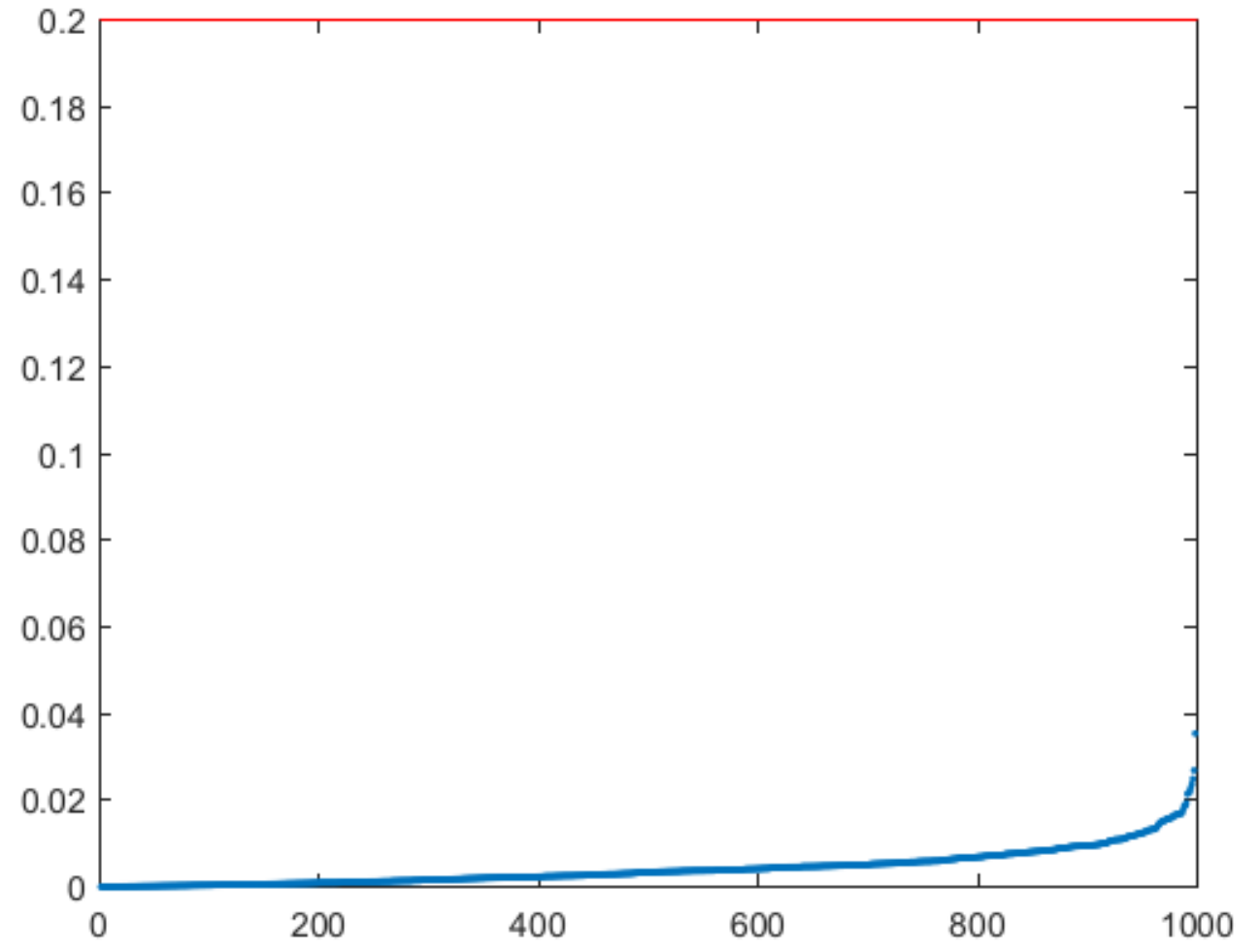}
 & 
\includegraphics[width=\linewidth]{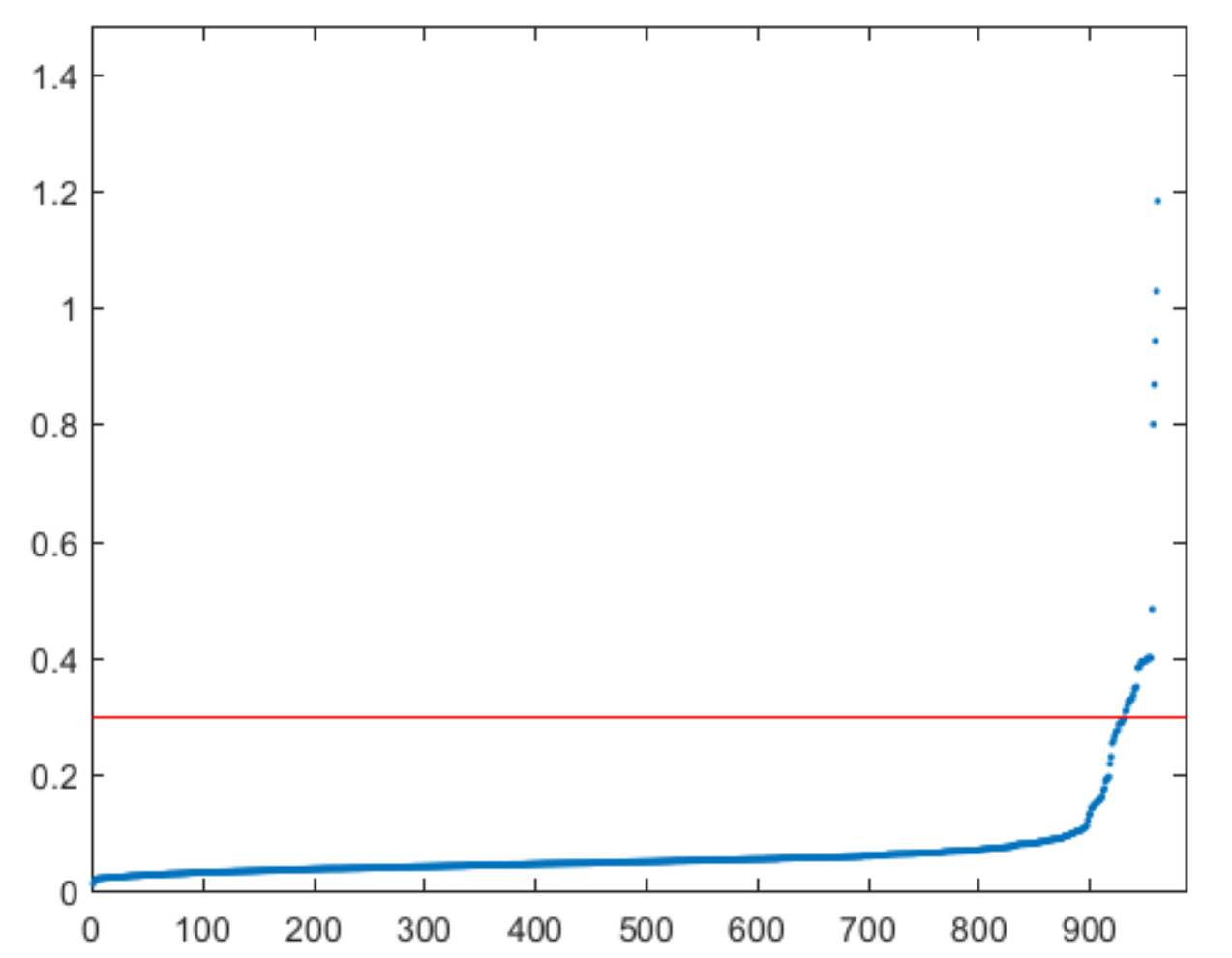}
\\
$n=500$ &   
\includegraphics[width=\linewidth, trim={6.7cm 12.1cm 6cm 11.7cm},clip]{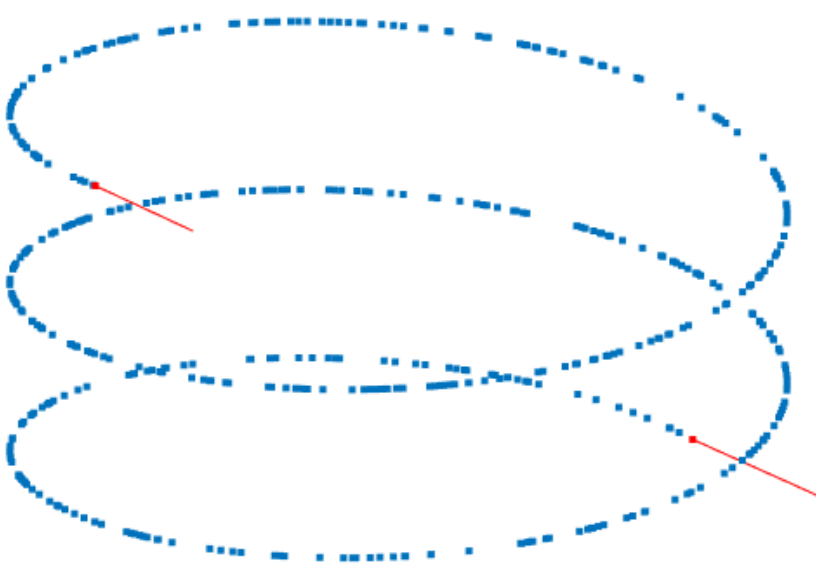}
&
\includegraphics[width=\linewidth, trim={5.7cm 10.8cm 5cm 10.5cm},clip]{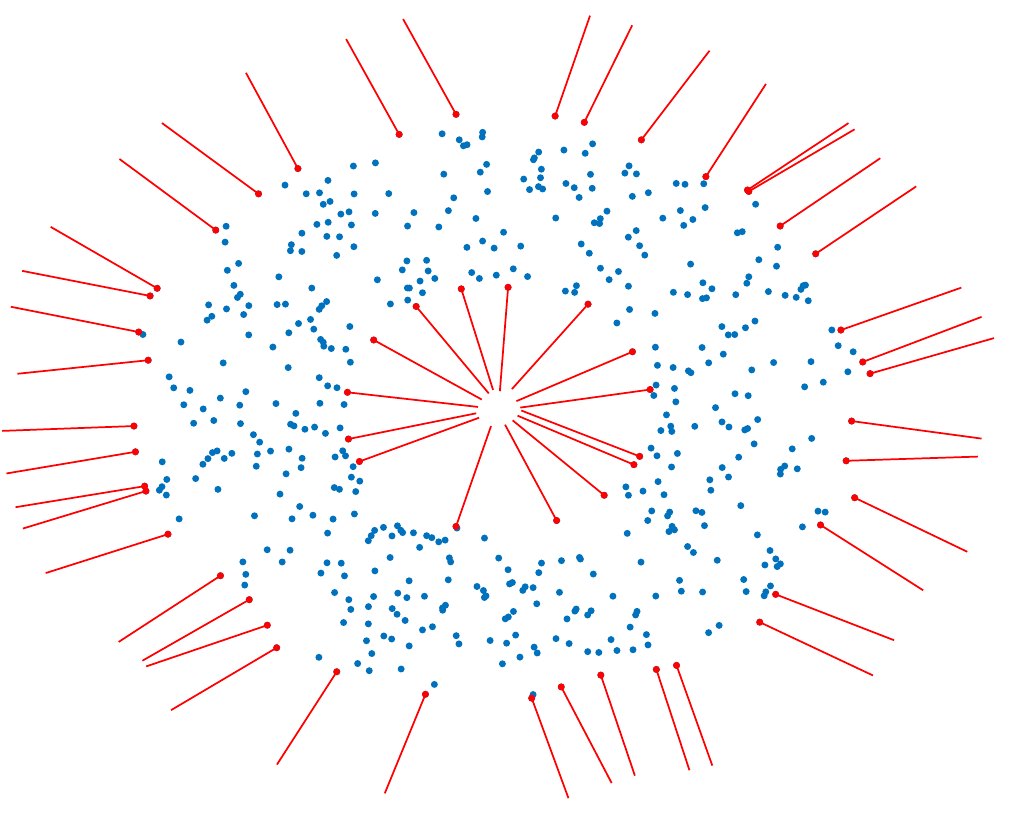}
\\
$n=1000$ &    
\includegraphics[width=\linewidth, trim={6.7cm 12.1cm 6cm 11.7cm},clip]{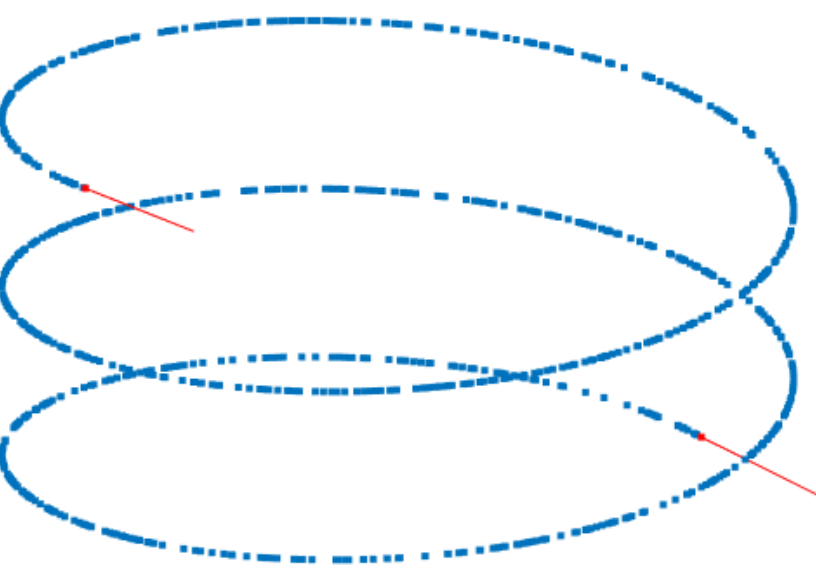}
&
\includegraphics[width=\linewidth, trim={5.7cm 10.8cm 5cm 10.5cm},clip]{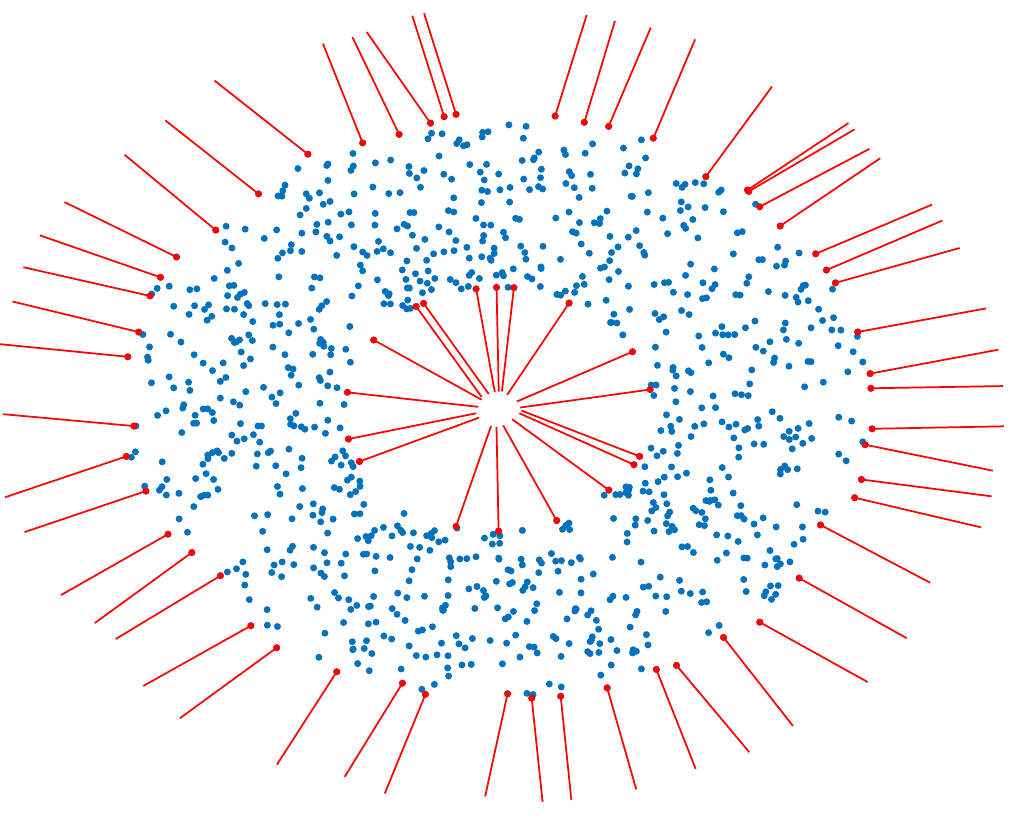}
\\
$n=2000$ &    
\includegraphics[width=\linewidth, trim={6.7cm 12.1cm 6cm 11.7cm},clip]{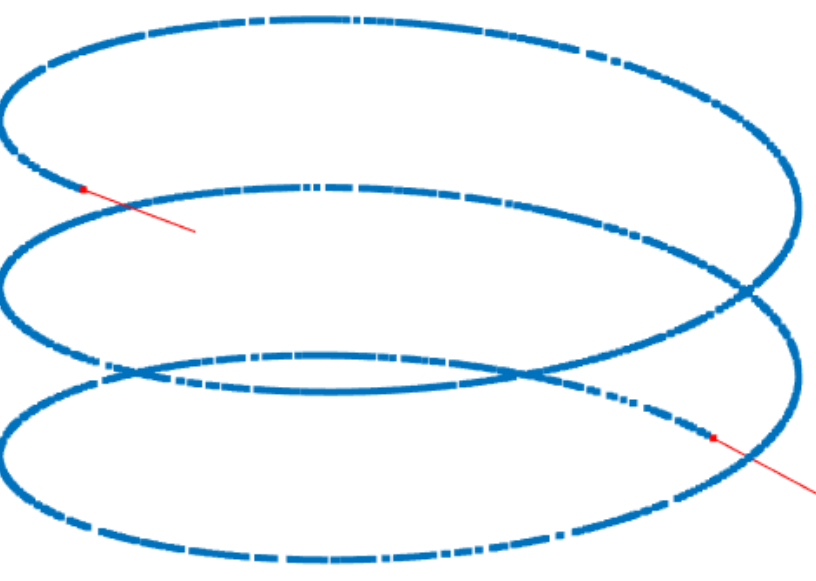}
&
\includegraphics[width=\linewidth, trim={5.7cm 10.8cm 5cm 10.5cm},clip]{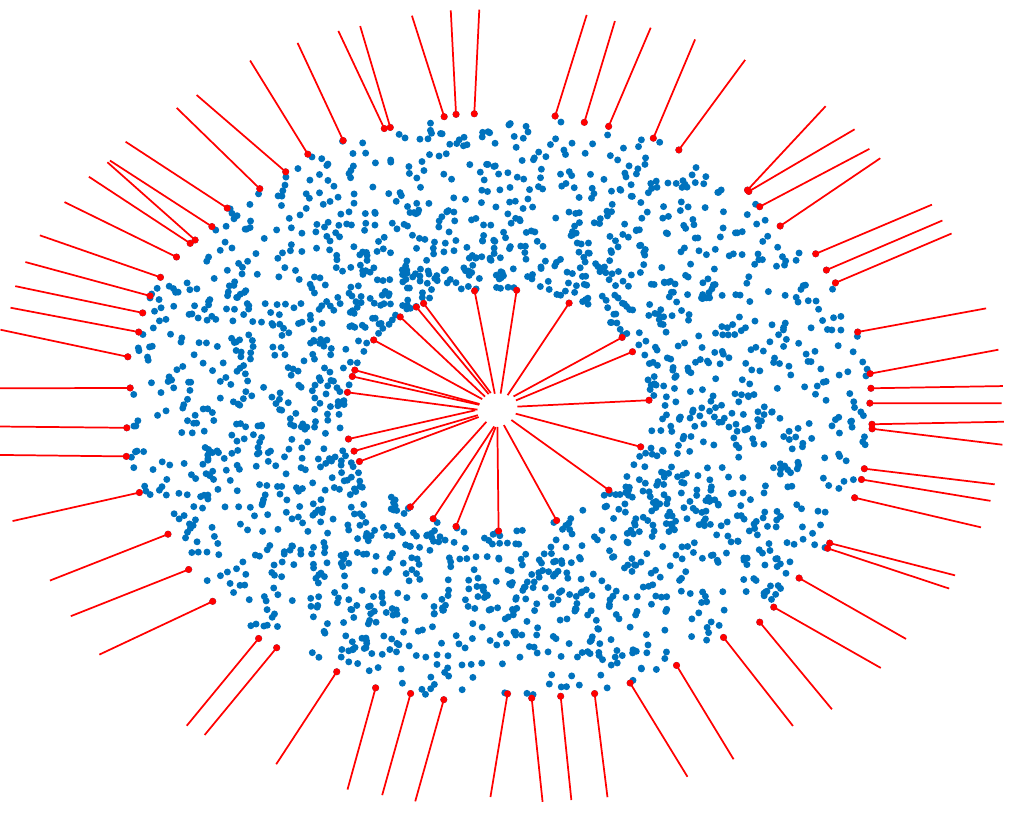}
\\
$n=5000$ &    
\includegraphics[width=\linewidth, trim={6.7cm 12.1cm 6cm 11.7cm},clip]{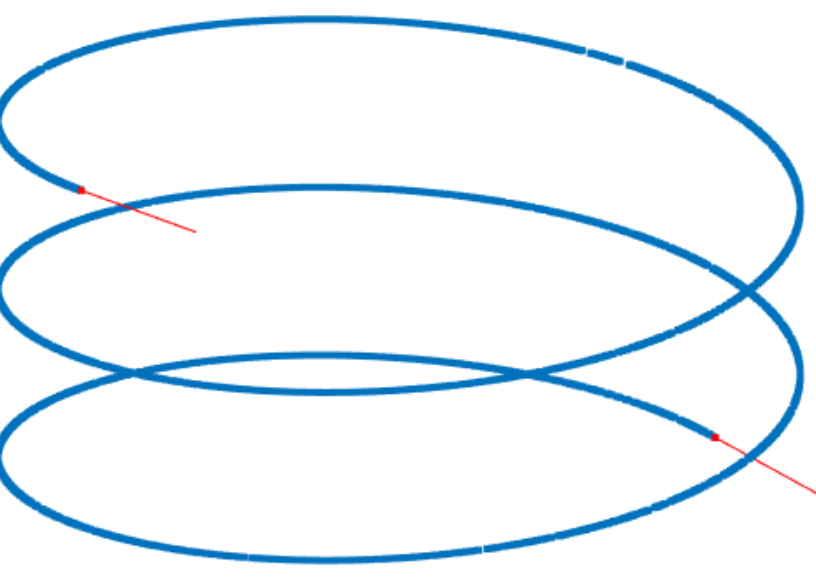}
&
\includegraphics[width=\linewidth, trim={5.7cm 10.8cm 5cm 10.5cm},clip]{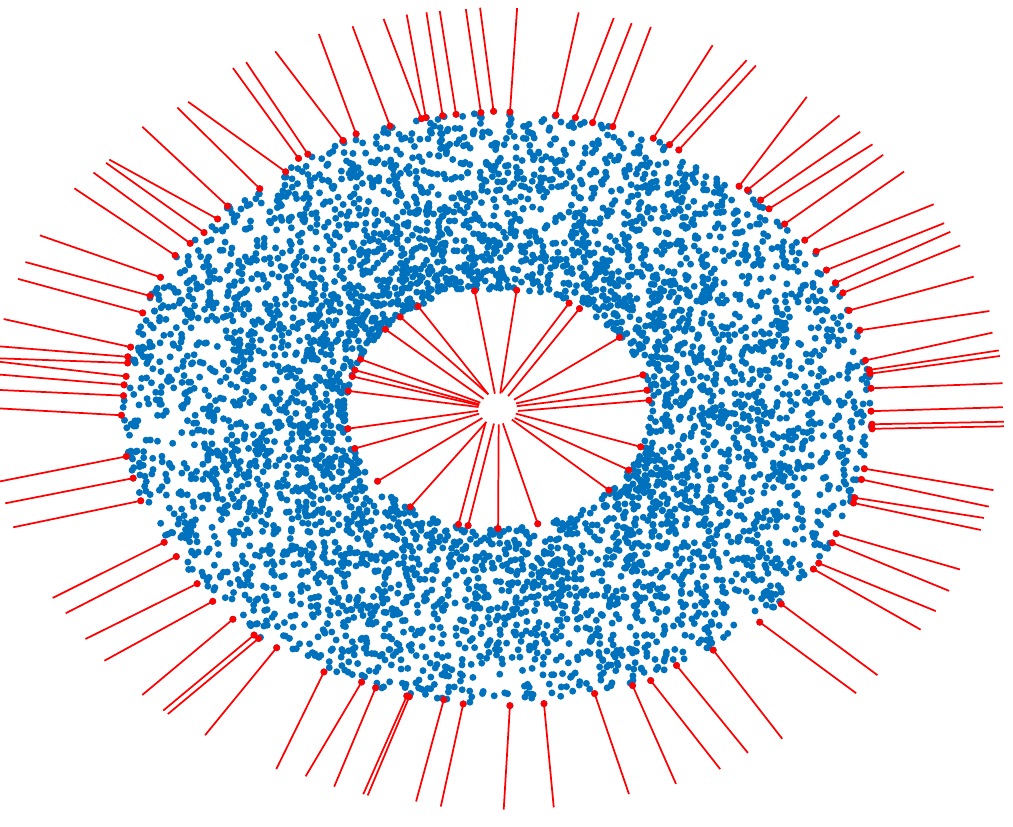}
\\
\bottomrule 
\end{tabular}
    \caption{Simulations results for the spiral and the annulus.}
    \label{fig:simu1}
\end{table}

\begin{table}
\centering
\begin{tabular}{c m{0.23\textwidth} m{0.23\textwidth}} 
\toprule 
\multicolumn{1}{c}{$M$} & 
\multicolumn{1}{c}{Half-sphere} & 
\multicolumn{1}{c}{M\"obius strip} \\ 
\midrule
Calibration of $R_0$ &
\includegraphics[width=\linewidth]{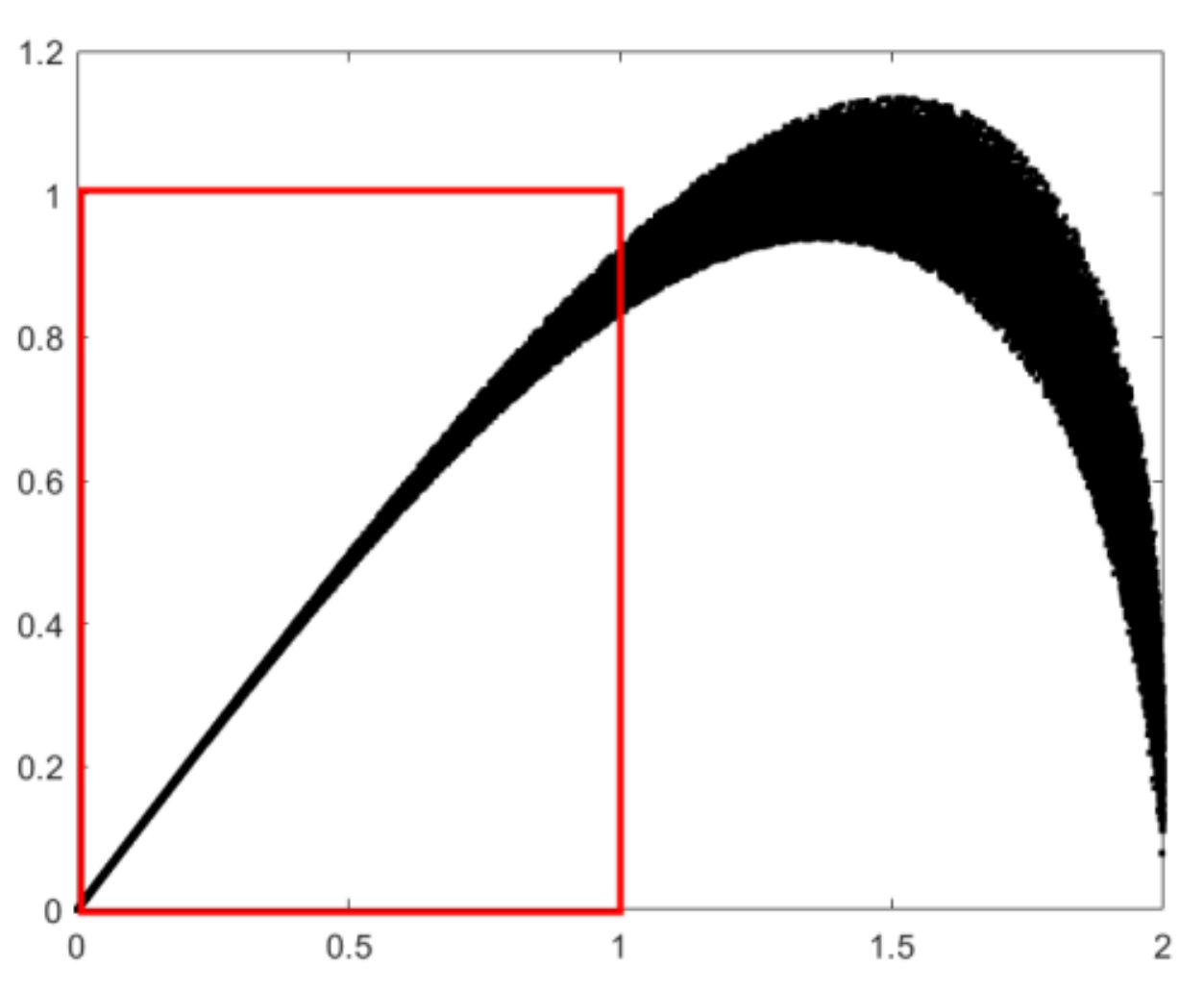}
 & 
\includegraphics[width=\linewidth]{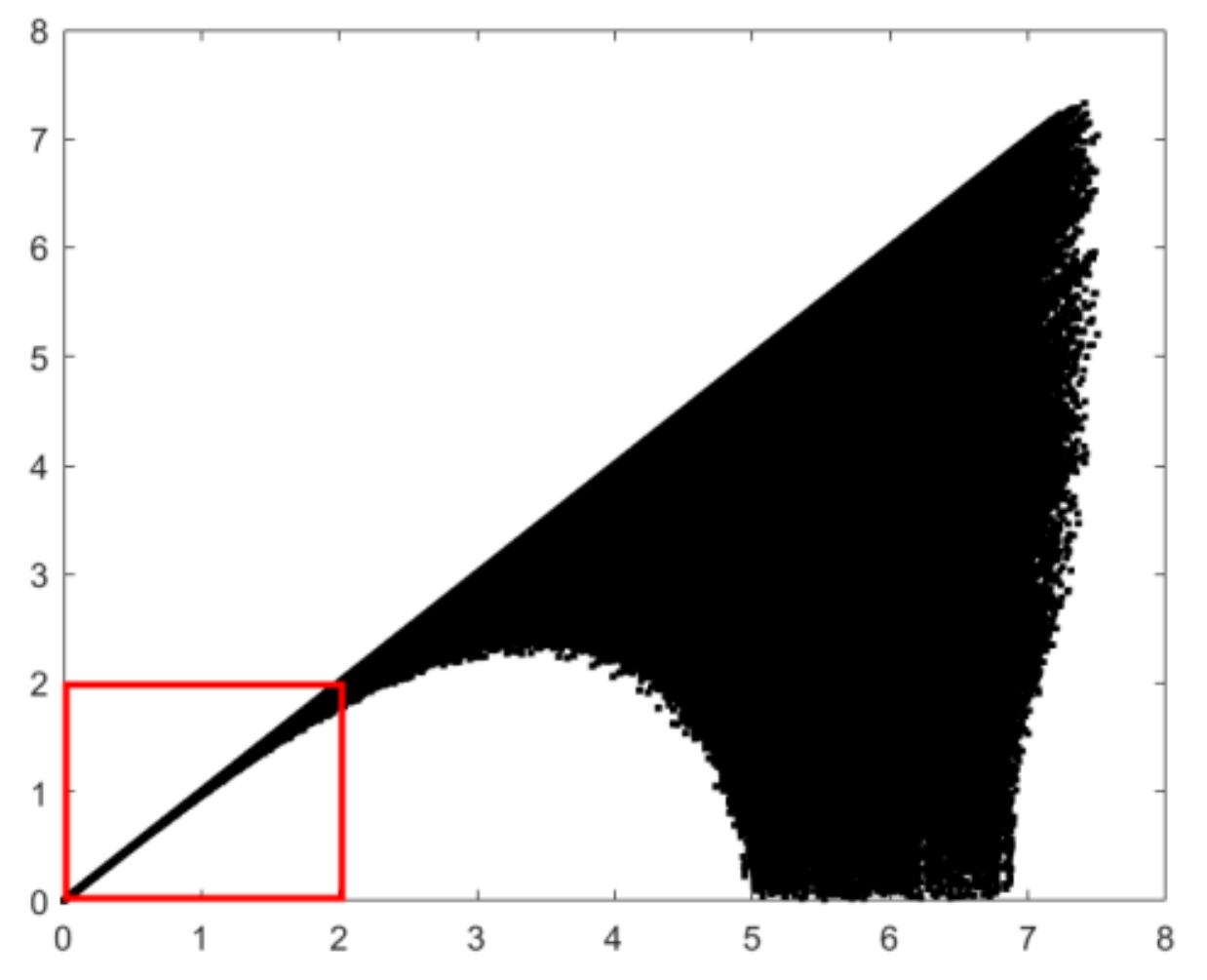}
\\
Calibration of $\rho$ &    
\includegraphics[width=\linewidth]{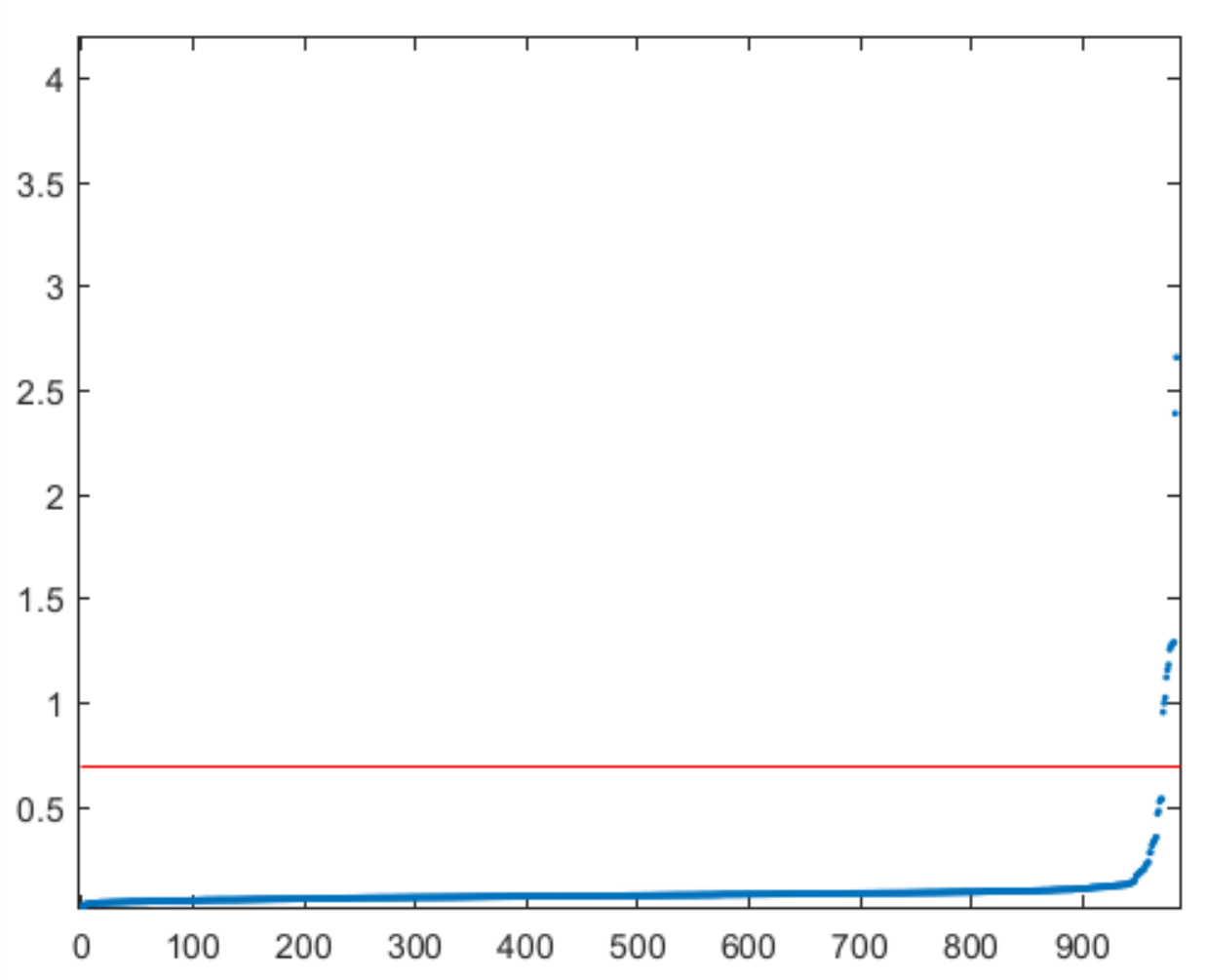}
 & 
\includegraphics[width=\linewidth]{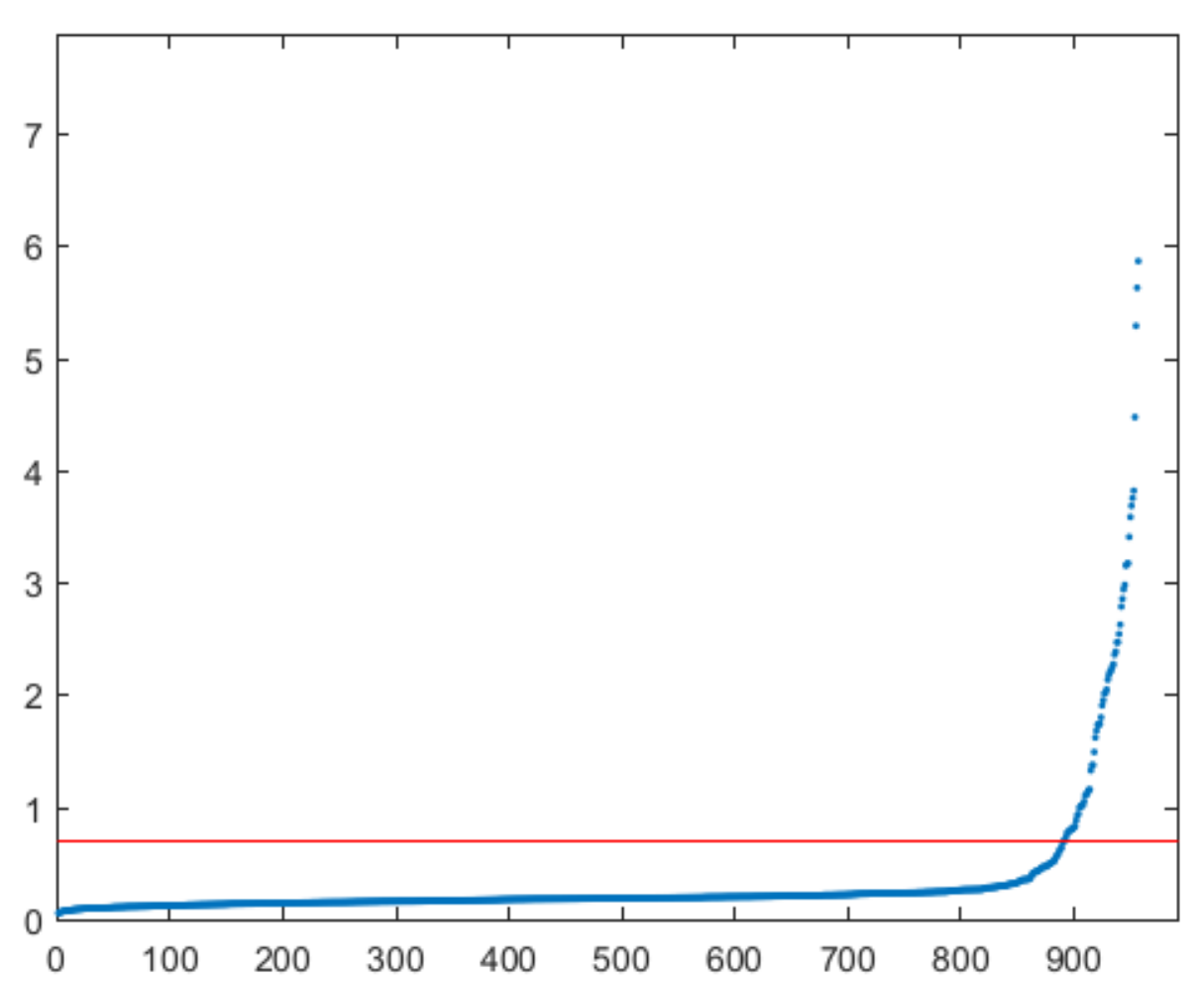}
\\
$n=500$ &    
\includegraphics[width = \linewidth, trim={5cm 11.5cm 6.4cm 12cm},clip]{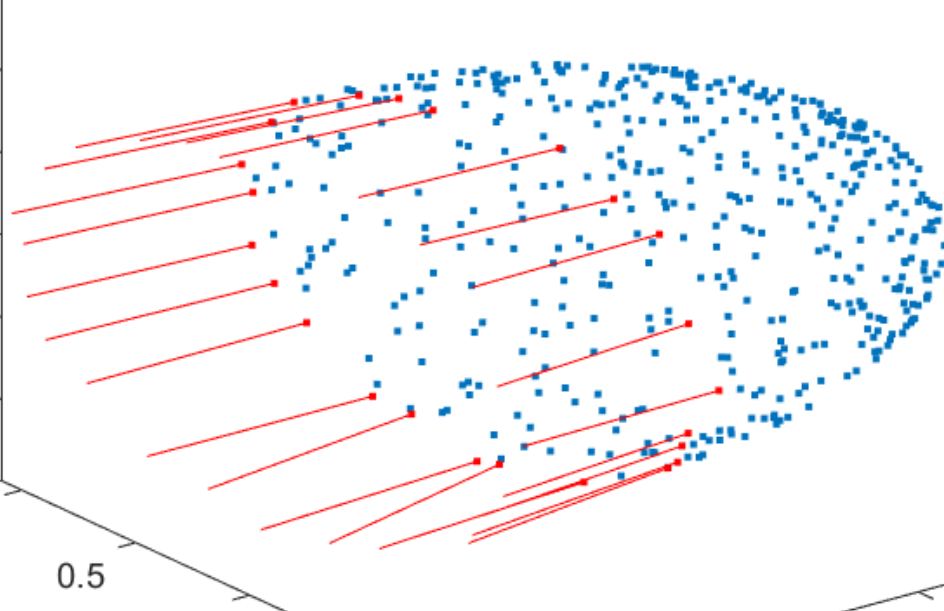}
&
\includegraphics[width = \linewidth, trim={7cm 12.5cm 5.8cm 12cm},clip]{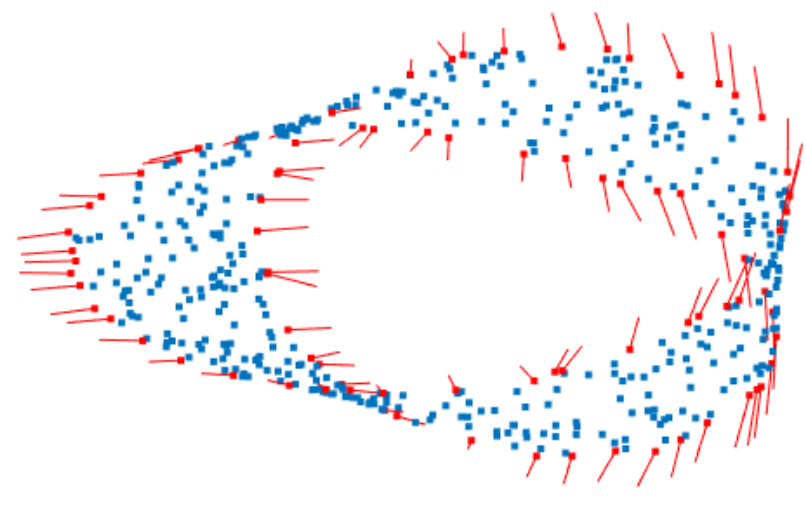}
\\
$n=1000$ &    
\includegraphics[width = \linewidth, trim={5cm 11.5cm 6.4cm 12cm},clip]{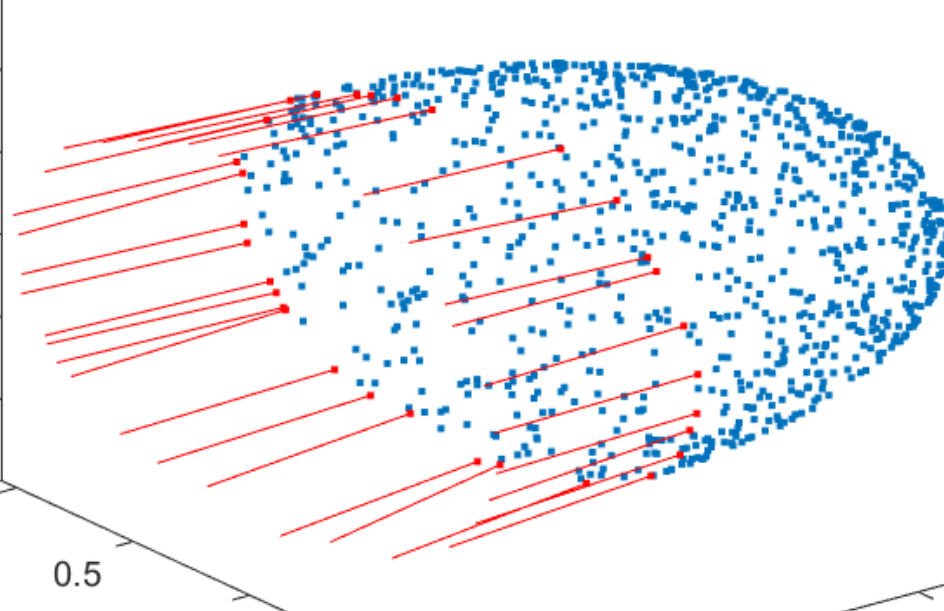}
&
\includegraphics[width = \linewidth, trim={7cm 12.5cm 5.8cm 12cm},clip]{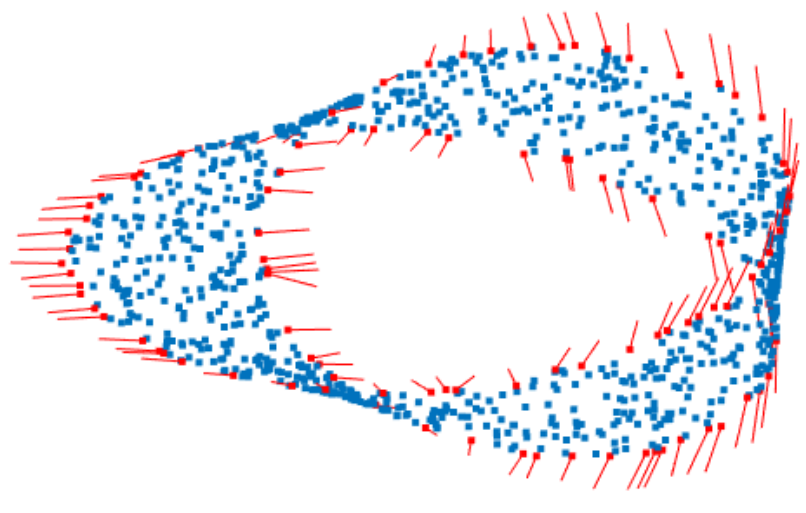}
\\
$n=2000$ &    
\includegraphics[width = \linewidth, trim={5cm 11.5cm 6.4cm 12cm},clip]{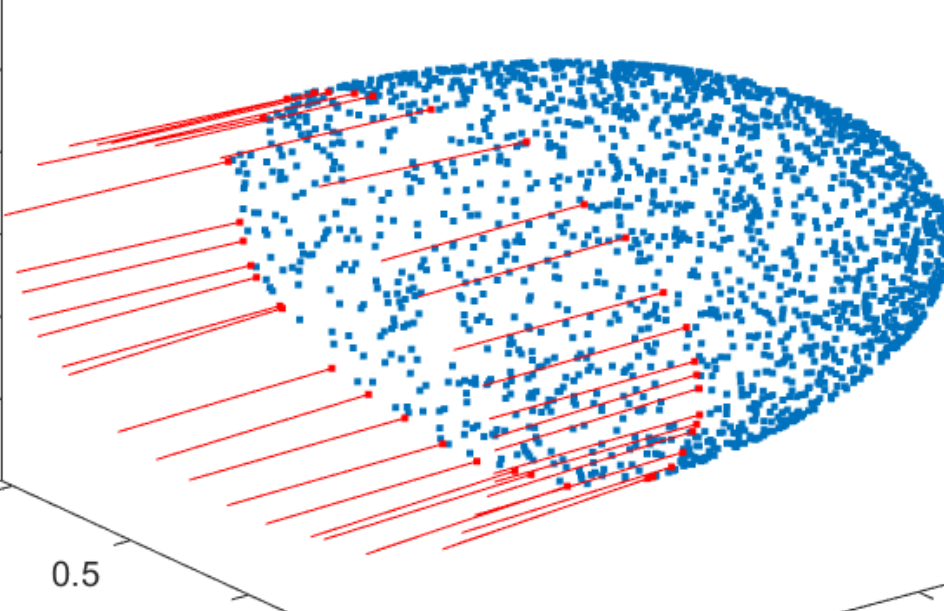}
&
\includegraphics[width = \linewidth, trim={7cm 12.5cm 5.8cm 12cm},clip]{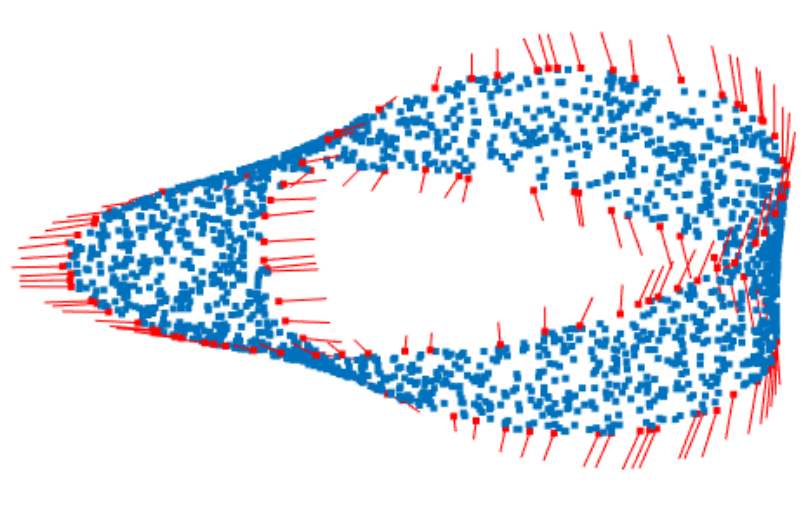}
\\
$n=5000$ &    
\includegraphics[width = \linewidth, trim={5cm 11.5cm 6.4cm 12cm},clip]{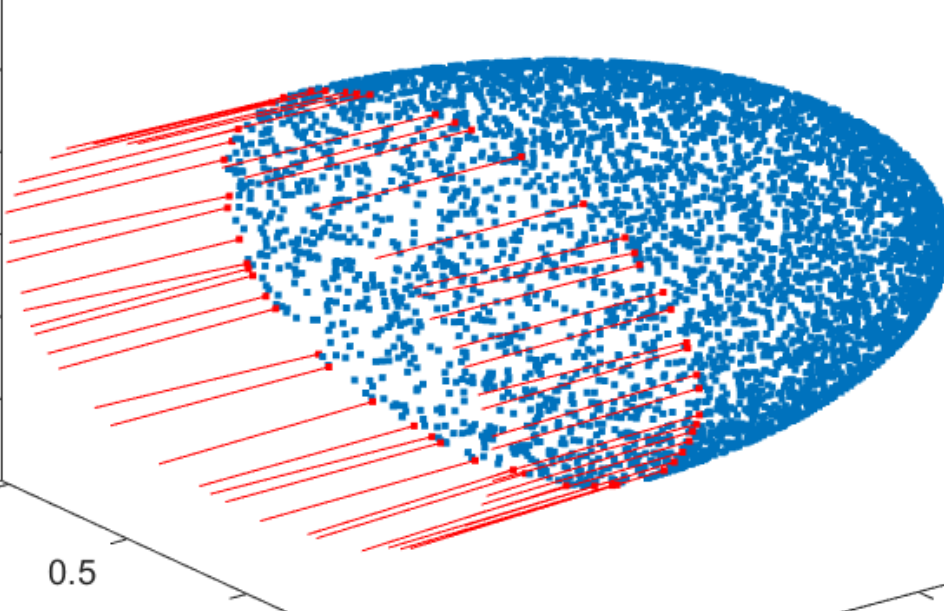}
&
\includegraphics[width = \linewidth, trim={7cm 12.5cm 5.8cm 12cm},clip]{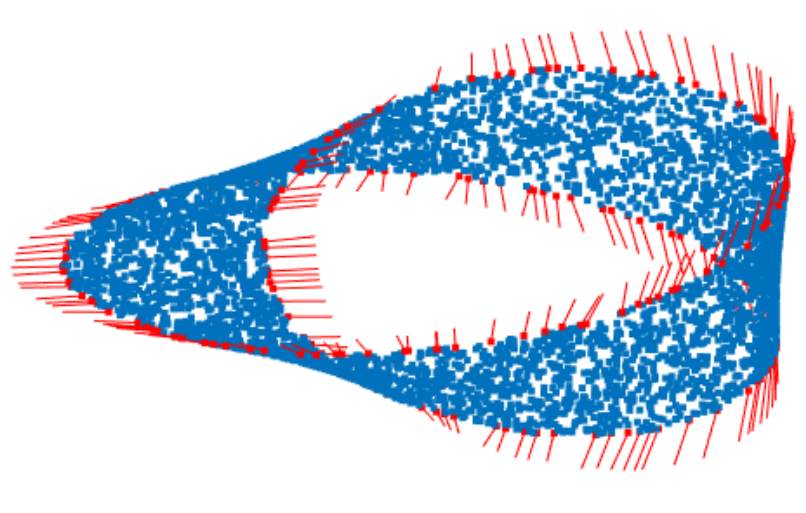}
\\
\bottomrule 
\end{tabular}
    \caption{Simulations results for the half-sphere and the M\"obius strip.}
    \label{fig:simu2}
\end{table}

%% file: geom_results.tex
\section{A shortlist of intermediate geometric results}
	\label{sec:geom_results}
	
This section gathers the main geometric results that are of use in the main derivations (\Cref{sec:proof-outline}).
For the sake of concision, proofs of these results are given in \Cref{sec:lemmas-manifolds-with-boundary}. Throughout, $\Gr{D}{d}$ stands for the Grassmannian --- i.e. the space of $d$-dimensional linear subspaces of $\mathbb{R}^D$ ---, and $\dd_S$ for the geodesic distance of $S \subset \R^D$.

\subsection{Geodesics and tangent spaces}

We begin with a result that connects geodesic and Euclidean distance.
\begin{lemma}[Geodesic Bounds]\label{lem:dist_geod_dist_eucl}
Let $p, q \in M$ such that $\|p-q\| \leq \tau_{\min}$. Then
\[
\|p-q\| \leq \dd_{M}(p,q) \leq 2\|p-q\|.
\]
\end{lemma}
A short proof is given in \Cref{tecsec:geodesics_tangent_spaces}.
This result is well-known in the empty boundary case (see~\cite[Proposition~8.6]{Aamari18}). In the general case,  \Cref{lem:dist_geod_dist_eucl} follows from  {\cite[Lemma~3]{Boissonnat19}}. The last result of this section connects tangent spaces variations with the geodesic distance between their base points.

\begin{restatable}[Tangent Space Stability]{proposition}{proptangentvariationgeodesic}
\label{prop:tangent_variation_geodesic}
Let $M \in \mathcal{M}^{d,D}_{\tau_{\min},\tau_{\partial,\min}}$. Then, for $x,y \in M$,
\begin{align*}
\angle(T_x M, T_y M)
&\leq
\dd_M(x,y)/\tau_{M}
.
\end{align*}
If $\partial M \neq \emptyset$, then for all $p,q \in \partial M$,
\begin{align*}
\angle(T_p \partial M, T_q \partial M)
&\leq
\dd_{\partial M}(p,q)/\tau_{\partial M}
.
\end{align*}
\end{restatable}
A proof of \Cref{prop:tangent_variation_geodesic} is given in \Cref{tecsec:geodesics_tangent_spaces}. Combining the two angle bounds from \Cref{prop:tangent_variation_geodesic} easily yields a bound on the angle between the \emph{linear spaces} $\Span(\eta_p)$ and $\Span(\eta_q)$, for $p,q \in \partial M$. Actually, making use of the structure of normal \emph{cones}, a bound on $\Vert \eta_p - \eta_q \Vert$ can be derived, as presented below.  

\begin{restatable}[Normal Vector Stability]{proposition}{propnormalvectorstability}
\label{prop:normal_vector_stability}
Let $M \in \mathcal{M}^{d,D}_{\tau_{\min},\tau_{\partial,\min}}$. Then for all $p,q \in \partial M$ such that $\|p-q\| \leq (\tau_M \wedge \tau_{\partial M})/32$, we have
\begin{align*}
\Vert \eta_p - \eta_q \Vert 
\leq 
9 \Vert p-q \Vert/(\tau_M \wedge \tau_{\partial M})
.
\end{align*}
\end{restatable}
A proof of \Cref{prop:normal_vector_stability} may be found in \Cref{tecsec:geodesics_tangent_spaces}.

\subsection{Projections}
Projections onto tangent spaces and normal directions play a key role in the estimation schemes on this work. 
First, we adapt \cite[Theorem~4.18]{Federer59} to the case where a small perturbation of the tangent space is allowed.

\begin{restatable}[Tangent and Normal Components of Increments]{proposition}{propprojectionsestimtangentM}
\label{prop:projectionsestimtangentM}
Let $x,y\in M$, and $T \in \Gr{D}{d}$ be such that $\angle(T_x M,T) \leq \theta$.
Write $(x-y)^{T}$ and $(x-y)^{{\perp}}$ for the orthogonal projection of $x-y$ onto $T$ and $T^{\perp}$ respectively. 
Then,
\begin{align*}
\Vert
(y-x)^{{\perp}}
\Vert
&\leq \norm{y-x}\left(\theta +\norm{y-x}/(2\tau_{\min})\right)
,
\\
\Vert
(y-x)^{T}
\Vert
&\geq \norm{y-x}\left(1-\theta-\norm{y-x}/(2\tau_{\min})\right)
.
\end{align*}
\end{restatable}

A proof of \Cref{prop:projectionsestimtangentM} is given in \Cref{tecsec:projections}. The following result ensures that estimates of the normal direction to the boundary may be derived from a suitable tangent space estimator.  

\begin{restatable}[Normals from Tangent Spaces]{proposition}{propexistencesestimnorM}
\label{prop:existencesestimnorM}
Let  $x \in \partial M$, and $T \in \Gr{D}{d}$ such that
$\angle(T_xM,T) < 1$.
Then $T\cap Nor(x,M)$ contains a unique unit vector $\eta$, and it satisfies
\[
\Vert \eta - \eta_x \Vert \leq \sqrt{2}\angle(T_x M,T).
\]
\end{restatable}

A proof of \Cref{prop:existencesestimnorM} can be found in \Cref{tecsec:projections}. The remaining results of this section describe the structure of the projection of balls onto perturbed tangent spaces. We begin by investigating the case where the center of the ball is not on the boundary.

\begin{restatable}[Far-Boundary Balls]{lemma}{lemGeoInterbleupartone}
\label{lem:GeoInter_bleupart1}
Let $x\in M$ and $T \in \Gr{D}{d}$ be such that $\angle (T_x M,T) \leq \theta\leq 1/8$.
If $\dd(x,\partial M)>0$ (with the convention $\dd(x,\emptyset) = + \infty$), and $R\leq \tau_{\min}/16$, then
\begin{align*}
\B_{T}\left(0,\frac{4}{5}\min\set{R, \dd(x,\partial M)} \right) \subset \pi_{T}(\B(x,R)\cap M - x).
\end{align*}
\end{restatable}

A proof of \Cref{lem:GeoInter_bleupart1} is given in \Cref{tecsec:ball-structure}. Next, \Cref{lem:GeoInter_bleupart2} describes $\pi_{T}(\B(x,R)\cap M-x)$ whenever $x$ is a boundary point.

\begin{restatable}[Near-Boundary Balls]{lemma}{lemGeoInterbleuparttwo}
\label{lem:GeoInter_bleupart2}
Assume that $\partial M \neq \emptyset$. 
Let $x\in \partial M$ and $T\in \Gr{D}{d}$ be such that $\angle (T_xM,T) \leq \theta\leq 1/8$.
Denote by $\hat{\eta}$ the unit vector of ${T}\cap Nor(x,M)$, choose $R\leq \tau_{\min}/16$ and $r \leq \min\set{2R/5,7\tau_{\partial,\min}/5}$. 

Then, writing
$O^{\text{in}}:=-r \hat{\eta}$
and $O^{\text{out}}:=r \hat{\eta}$, we have
\begin{equation*}
\B(O^{\text{in}},r)\cap {T}\subset \pi_{{T}}(\B(x,R) \cap M - x)\subset \Bopen(O^{\text{out}},r)^c\cap {T}.
\end{equation*}

\end{restatable}

A proof of \Cref{lem:GeoInter_bleupart2} may be found in \Cref{tecsec:ball-structure}. A consequence of \Cref{lem:GeoInter_bleupart2} is the following \Cref{cor:projorthobord}, that will be useful in the proof of \Cref{thm:deterministic}.

\begin{restatable}[Parallelism of Projected Normals]{corollary}{corprojorthobord}
\label{cor:projorthobord}
Assume that $\partial M \neq \emptyset$. Let $x\in M$ be such that $\dd(x,\partial M)< \tau_{\min}/16$, and $y \in \R^D$. 
For $T \in \Gr{D}{d}$, 
let $x^* \in \pi_{y+T}(\partial M\cap \B(x,\tau_{\min}/16))$ be any point such that
$$
\norm{
x^*-\pi_{y+T}(x)
}
=
\dd(\pi_{y+T}(x),\pi_{y+T}(\partial M\cap \B(x,\tau_{\min}/16))
,
$$
and
$$x'\in \partial M\cap \B(x,\tau_{\min}/16) \text{ such that }\pi_{y+T}(x')=x^*
.
$$
If $\angle (T_{x'}M,T) \leq 1/8$, then  $Nor(x',M)\cap {T}$ contains a unique unit vector $\eta^*(x')$, and
\[
x^*-\pi_{{y+T}}(x)=\norm{x^*-\pi_{y+{T}}(x)}\eta^*(x').
\]
\end{restatable}

A proof of \Cref{cor:projorthobord} is given in \Cref{tecsec:ball-structure}.

\subsection{Covering and volume bounds}
\label{sec:volume_bounds}

This last preliminary section provides probabilistic bounds on the sampling density of $\X_n$ in $M$, and bounds on the volume of intersection of balls. They will drive the convergence rates of \Cref{thm:detection}. 
First, we adapt \cite[Lemma~9.1]{Aamari18} to the non empty boundary case.
\begin{restatable}[Sampling Density Bound]{lemma}{lemcovering}
\label{lem:covering}
Let $\varepsilon_1 = \left ( C_d \frac{\log n}{f_{\min}n} \right )^{\frac{1}{d}}$, for $C_d$ large enough. Then, for $n$ large enough so that $\varepsilon_1 \leq \frac{\tau_{\min}}{16} \wedge \frac{\tau_{\partial,\min}}{2}$, we have, with probability larger than $1-n^{-3}$, 
\[
\dHaus(M,\X_n) \leq \varepsilon_1.
\]  
\end{restatable}

A proof of \Cref{lem:covering} is given in \Cref{tecsec:volumes_and_covering}. It guarantees that the convergence rate of the sample $\X_n$, seen as a Hausdorff estimator of $M$, is the same as in the empty boundary case. 
Next, \Cref{lem:proba_occupation_losange} below provides bounds on the mass of projected intersection of balls.   

\begin{restatable}[Mass of Intersection of Curved Balls]{lemma}{lemprobaoccupationlosange}
\label{lem:proba_occupation_losange}
Let $x \in M$, and $T \in \Gr{D}{d}$. 
Let $O \in T$, and $r,R \geq 0$ be such that $\B_T(O,r) \subset \pi_T(\B(x,R) \cap M -x)$. 
For $A \geq C'_d r^{\frac{1-d}{2}}$, write 
$$
h = \left( \frac{C_d f_{max}^4 \log n}{f_{\min}^{5}(n-1)} \right )^\frac{1}{d}
\text{, and~} 
\varepsilon_2 = \left ( A \frac{f_{max}^4\log n}{f_{\min}^{5}(n-1)} \right )^\frac{2}{d+1}
.
$$
Then for $n$ large enough, for all $\rho \geq r$ and $\Omega \in T$ such that $\| \Omega - O \| \leq r+ \rho - \varepsilon_2$,
\[
\int_{M \cap (\B(x,R) \setminus \B(x,h))} \indicator{\pi_T(u-x) \in \B(O,r) \cap \B(\Omega,\rho)} f(u) \Haus^d(\dd u) 
\geq 
A r^{\frac{d-1}{2}} C''_d  \frac{f_{max}^4 \log n}{f_{\min}^4(n-1)}.
\]
\end{restatable}
A proof of \Cref{lem:proba_occupation_losange} can be found in \Cref{tecsec:volumes_and_covering}. From a sampling point of view, it will ensure that such intersections of (projected) balls will contain at least one sample point with high probability. This point will allow to detect and characterize the boundary observations (see \Cref{thm:deterministic}).

%% file: lemmas-manifolds-with-boundary.tex
\section{Geometric properties of manifolds with boundary}
	\label{sec:lemmas-manifolds-with-boundary}

This Section gathers the proofs for \Cref{sec:geom_results}. 
To ease readability, statements are recalled before their proofs.

\subsection{Geodesics and tangent space variations}\label{tecsec:geodesics_tangent_spaces}
In addition to the Euclidean structure induced by $\R^D$ on $M \subset \R^D$, we can also endow $M$ and $\partial M$ with their intrinsic geodesic distances $\dd_M$ and $\dd_{\partial M}$ respectively.
To cover both cases at once, let $S \in \set{M,\partial M}$.
Given a $\mathcal{C}^1$ curve $c: [a,b] \rightarrow S$, the length of $c$ is defined as $\Length(c)= \int_a^b \norm{c'(t)}dt$. 
Given $p,q \in S$ belonging to the same connected component of $S$, there always exists a path $\gamma_{p \rightarrow q}$ of minimal length joining $p$ and $q$ \cite[Proposition~2.5.19]{Burago01}.
Such a curve $\gamma_{p \rightarrow q}$ is called geodesic, and the geodesic distance between $p$ and $q$ is given by $\dd_S(p,q) = \Length(\gamma_{p \rightarrow q})$. 
If $x$ and $y$ stand in different connected components of $S$, then $\dd_S(x,y) = \infty$.

A geodesic $\gamma$ such that $\norm{\gamma'(t)} = 1$ for all $t$ is called \emph{arc-length parametrized}. Unless stated otherwise, we always assume that geodesics are parametrized by arc-length.
If $S$ has empty boundary, then for all $p \in S$ and all unit vectors $v \in T_p S$, we denote by $\gamma_{p,v}$ the unique arc-length parametrized geodesic of $S$ such that $\gamma_{p,v}(0)=p$ and $\gamma'_{p,v}(0)=v$ \cite[Chap.~7, Theorem~2.8]{DoCarmo92}.
The exponential map is then defined as $\exp_p^S(vt) = \gamma_{p,v}(t)$. Note that if in addition $S$ is compact, $\exp_p^S: T_p S \rightarrow S$ is defined globally on $T_p S$ \cite[Theorem 2.5.28]{Burago01}. We let $\B_S(p,s)$ denote the closed geodesic ball of center $p\in S$ and of radius $s \geq 0$.

Although they might differ drastically at long range, geodesic and Euclidean distances are good approximations of one another when evaluated between close enough points. The following result quantifies this intuition, and implies \Cref{lem:dist_geod_dist_eucl}.

\begin{proposition}
\label{prop:geodesic_vs_euclidean}
Let $S \subset \R^D$ have positive reach $\tau_S > 0$, and $x,y \in S$ be such that $\norm{y-x} \leq \tau_S$.
Then,
\begin{align*}
\norm{y-x}
\leq
\dd_S(x,y)
\leq
\left(
1
+
\frac{\norm{y-x}^2}{20\tau_S^2}
\right)
\norm{y-x}
.
\end{align*}
\end{proposition}

\begin{proof}[Proof of \Cref{prop:geodesic_vs_euclidean}]

We clearly have $\norm{y-x} \leq \dd_S(x,y)$, and on the other hand, \cite[Lemma 3]{Boissonnat19} yields
\begin{align*}
\dd_S(x,y)
&\leq
2\tau_S
\arcsin \left(\frac{\norm{y-x}}{2\tau_S} \right)
\leq
\left(
1
+
\frac{\norm{y-x}^2}{20\tau_S^2}
\right)
\norm{y-x}
,
\end{align*}
where the last inequality follows uses that $\arcsin t \leq t(1+t^2/5)$ for all $0 \leq t \leq 1/2$.
\end{proof}

Next, we ensure that the angle between tangent spaces can be bounded in terms of geodesic distances between base points. In the empty boundary case, this result is well known, and can be shown using via parallel transportation of tangent vectors (see the proof of \cite[Lemma~A.1]{Aamari19b}). In the general case, the tangent space stability property writes as follows.

\proptangentvariationgeodesic*

\begin{proof}[Proof of \Cref{prop:tangent_variation_geodesic}]
If $\partial M = \emptyset$, the first claim follows from \cite[Lemma~6]{Boissonnat19}.

Assume that $\partial M \neq \emptyset$. From \Cref{prop:boundary_is_manifold}, $\partial M$ is a $\mathcal{C}^2$-submanifold without boundary. Then, the second statement also directly follows from \cite[Lemma~6]{Boissonnat19}.
For the first claim, the key technical point is to handle geodesics that would hit the boundary. 

To do this we define a push-inwards operator that will allow to consider path in the interior of $M$ only. First, an elementary results on an atlas of $M$ is needed.

\begin{lemma}\label{lem:atlas_bord}
Let $U_1, \ldots, U_k$ be charts of $M$ that cover $\partial M$. 
Then there exists $r_0>0$ such that 
\[
\forall p \in \partial M \quad \exists j \in \set{1,\ldots,k} \quad \Bopen(p,r_0) \cap M \subset U_j \cap M.
\]
\end{lemma}

We now consider a smooth kernel $K : \R_+ \to [0,1]$ such that
\[
K(x) = 
\begin{cases}
1 & \mbox{if } x \leq \tau_{\partial M}/4 \\
0 & \mbox{if } x \geq \tau_{\partial M}/2
\end{cases}
\]
and we define the vector field $\mathbf{V}$ on $M$ by
\[
\mathbf{V}(p) 
:= 
\begin{cases}
K\left [ \dd(p, \partial M) \right ]\pi_{T_p M} \left (\nabla_p \left ( \dd(\cdot,\partial M)  \right ) \right )
&
\text{if } \dd(p,\partial M) < (r_0 \wedge \tau_{\partial M})/2,
\\
0
&
\text{otherwise.}
\end{cases}
\]
Note that if $q \in \partial M$, $\nabla_q \dd(\cdot, \partial M) = - \eta_q$, where $\eta_q$ is the unit outward-pointing normal vector at $q$. By construction, $\mathbf{V}$ is a $\mathcal{C}^1$ tangent vector field on $M$. We now examine its flow.
\begin{lemma}\label{lem:flow_inward_definition}
For all $p \in M$, the flow of $\mathbf{V}$ starting from $p$ is defined globally on $\R_+$.
\end{lemma}
Equipped with Lemma \ref{lem:flow_inward_definition}, we may define our \emph{push-inwards} operator as follows:
\begin{align*}
g_\varepsilon \colon M &\to M 
\\
p &\mapsto g(p,\varepsilon) 
\end{align*}
where $g(p,t)$ denotes the flow of $\mathbf{V}$ at time $t\geq 0$ starting from $p \in M$.
The following properties of $g_\varepsilon$ will shortly be of technical interest.

\begin{lemma}\label{lem:properties_push_inward}
For all $p \in M$ and $\varepsilon > 0$,
\begin{center}
~\hfill
$\|g_\varepsilon(p) - p\| \leq \varepsilon$,
\hfill
$g_\varepsilon(p) \notin \partial M$,
\hfill
and $\|d_p g_\varepsilon - Id_{T_p M}\|_{op} \leq K \varepsilon e^{K \varepsilon}$,
\hfill
~
\end{center}
where $K = \sup_{p \in M} \|d_p \mathbf{V}\|_{op}$.
\end{lemma}

We can now finish the proof of the first result in Proposition \ref{prop:tangent_variation_geodesic}. We let $p,q \in M$, and $\gamma$ a unit-speed curve joining $p$ and $q$ whith length $\dd_M(p,q)$. We define $\gamma_{\varepsilon}$ as the push-inwards of $\gamma$, that is
\[
\gamma_\varepsilon(t) := g_\varepsilon(\gamma(t)), 
\]
for all $t \in [0,\dd_M(p,q)]$. 
As $g_\varepsilon(p) \notin \partial M$ for all $p \in M$ (\Cref{lem:properties_push_inward}), parallel transportation of tangent vectors in the interior $\Int M$ of $M$ (see for instance the proof of \cite[Lemma~A.1]{Aamari19b}) yields that
\[
\angle (T_{p_\varepsilon}M, T_{q_\varepsilon}M) \leq \frac{L(\gamma_\varepsilon)}{\tau_{M}},
\]
where $p_\varepsilon = g_\varepsilon(p)$, $q_\varepsilon = g_\varepsilon(q)$, and $L(\gamma_\varepsilon)$ denotes the length of $\gamma_\varepsilon$. But from \Cref{lem:properties_push_inward} again,
\begin{align*}
L(\gamma_\varepsilon) = \int_0^{\dd_M(p,q)}\|\gamma'_\varepsilon(t)\|\dd t 
& = \int_0^{\dd_M(p,q)} \|d _{\gamma(t)}g_\varepsilon \left [ \gamma'(t) \right ]\|\dd t \leq (1 + K \varepsilon e^{K \varepsilon}) \dd_M(p,q).
\end{align*}
$\angle (T_{p_\varepsilon M}, T_p M) \leq K \varepsilon e^{K \varepsilon}$, and
$\angle (T_{q_\varepsilon} M, T_q M) \leq K \varepsilon e^{K\varepsilon}$.
As a result, triangle inequality yields 
\[
\angle (T_p M, T_q M)
\leq 
2 K \varepsilon e^{K \varepsilon} 
+ 
(1 + K \varepsilon e^{K \varepsilon}) \frac{\dd_M(p,q)}{\tau_M}
,
\]
so that the result follows after letting $\varepsilon \to 0$.
\end{proof}

We finally prove the intermediate results of \Cref{lem:atlas_bord,lem:flow_inward_definition,lem:properties_push_inward} that we just used to derive \Cref{prop:tangent_variation_geodesic}.

\begin{proof}[Proof of Lemma \ref{lem:atlas_bord}]
For all $p \in \partial M$, set 
$$r(p) := \sup\{r \geq 0 \mid \exists j \in \set{1,\ldots,k}, \Bopen(p,r) \subset U_j\} .$$ 
Note that since $(U_i)_ {1 \leq i \leq k}$ is an open covering of $\partial M$ we have $r(p) >0$. Consider
\[
r_0 := \inf_{p \in \partial M} r(p),
\]
which clearly satisfies the announced statement by definition.
Suppose, for contradiction, that $r_0=0$. 
Then there would exist a sequence $(p_n)_{n \in \mathbb{N}} \in (\partial M)^{\mathbb{N}}$ such that $r(p_n) \to 0$. As $\partial M$ is compact, we may assume (up to extraction) that $p_n \rightarrow p \in \partial M$ as $ n \to +\infty$. 
As a result, for $n$ large enough, we have $\Bopen(p_n,r(p_n)) \subset \Bopen(p,r(p)) \subset U_{j_0}$ for some $j_0$, which is a contradiction. 
\end{proof}

\begin{proof}[Proof of Lemma \ref{lem:flow_inward_definition}]
We distinguish cases according to the value of $\dd(p,\partial M)$ with respect to the chart radius $r_0$ of \Cref{lem:atlas_bord}.
\begin{itemize}[leftmargin=*]
\item
If $\dd(p,\partial M) \geq \tau_{\partial M}/2$, then $\mathbf{V}(p)=0$ and the flow of $\mathbf{V}$ starting from $p$ is $p(t) = p$ for all $t \geq 0$.
\item
If $r_0/2 < \dd(p,\partial M) \leq \tau_{\partial M}/2$, then we may find $r_1 \in (0,r_0/2)$ such that $\Bopen(p,r_1) \cap M$ is diffeomorphic to an open subset of $\R^d$. Using Cauchy-Lipschitz theorem in this chart space, we get that there exists $t_0 >0$ such that the flow of $\mathbf{V}$ starting from $p$ is well-defined at least on $[0, t_0)$.
\item
If  $\dd(p,\partial M) \leq r_0/2$, denote by $q = \pi_{\partial M}(p)$ and let $j \in \set{1,\ldots,k}$ be such that $\Bopen(q,r_0) \subset U_{j}$, where $\psi_j:U_j \cap M \rightarrow (\R^{d-1}\times \R_+) \cap \psi_j(U_j)$ is a chart of $M$.
Without loss of generality we may assume that $d_q (\psi_j) (\eta_q) = -e_d$, where $e_d$ is the $d$-th vector of the canonical basis of $\R^d$.

Let $r_1>0$ be such that $V_1=\Bopen(\psi_j(p),r_1)\cap (\R^{d-1}\times \R_+) \subset (\R^{d-1}\times \R_+) \cap \psi_j(U_j)$, and denote by $\mathbf{V}_2$ the vector field on $V_1$ defined by $d \psi_j \left [ \mathbf{V} \right ]$. Then $\mathbf{V}_2$ can be extended into a Lipschitz vector field $\mathbf{V}_3$ on $\Bopen(\psi_j(p),r_1)$, by choosing $\mathbf{V}_3(x_1, \hdots, x_d) = \mathbf{V}_3(x_1, \hdots, 0)$ if $x_d \leq 0$.  

Then, the Cauchy-Lipschitz theorem ensures that there exists $t_0$ such that the flow of $\mathbf{V}_3$ starting from $\psi_{j}(p)$ is defined on $]-t_0,t_0[$. 
Let $g_2(t,\psi_j(p))$ denote this flow. 
According to \Cref{lem:atlas_bord}, it holds $ \left\langle \mathbf{V}_3(g_2(0,\psi_j(p))), e_d \right\rangle =1$. Thus, there exists $t_1 >0$ such that for all $t \in [0,t_1]$, $g_2(t,p) \in \Bopen(\psi_j(p),r_1) \cap (\R^{d-1}\times \R_+)$, and therefore the flow of $\mathbf{V}_3$ starting from $\psi_j(p)$ stays in $\Bopen(\psi_j(p),r_1) \cap (\R^{d-1}\times \R_+)$.
When pushed back, this means that the flow of $\mathbf{V}$ starting from $p$ stays in the chart $(U_j,\psi_j)$.
\end{itemize}

In summary, we have shown that for all $p \in M$ there exists $t_p >0$ such that the flow $g(t,p)$ of $\mathbf{V}$ starting from $p$ is well-defined for $t \in [0,t_p]$. Since $g(\cdot,p)$ goes to the compact $M$ and satisfies $g(t_1 + t_2,p) = g(t_2, g(t_1,p))$, we deduce that for all $p \in M$, $g(\cdot,p)$ is well-defined on $\R_+$.
\end{proof}

\begin{proof}[Proof of Lemma \ref{lem:properties_push_inward}]
Since $\|\mathbf{V}\| \leq 1$, we directly get that 
$$\|g_\varepsilon(p) - p\| = 
\norm{\int_0^\varepsilon \mathbf{V}(g(p,t))\dd t} 
\leq 
\int_0^\varepsilon \norm{\mathbf{V}(g(p,t))}\dd t
\leq
\varepsilon.$$
To obtain the second point, write $\dd(g_\varepsilon(p),\partial M) -
\dd(p,\partial M)$ as
\begin{align*}
\int_{0}^\varepsilon &
\left\langle \mathbf{V}(g(p,t)),  \nabla_{g(p,t)}\dd(\cdot,\partial M) \right \rangle 
\dd t
\\ 
&= \dd(p,\partial M) + \int_{0}^\varepsilon K \left [ \dd(g(p,t),\partial M) \right ] \left\langle \pi_{T_{g(p,t)}M} (\nabla_{g(p,t)}\dd(\cdot,\partial M)),  \nabla_{g(p,t)}\dd(\cdot,\partial M) \right \rangle
\dd t
.
\end{align*}
Thus,
\begin{itemize}
\item
If $p \notin \partial M$, then $\dd(g_\varepsilon(p),\partial M) \geq \dd(p,\partial M) > 0$. 
\item
If $p \in \partial M$, then $\pi_{T_p M} \left ( \nabla_{g(p,0)} \dd(\cdot,\partial M) \right ) = - \eta_p$. Since $\mathbf{V}$ is continuous, there exists $t_0$ such that for all $t \leq t_0$, we have $$\left\langle \mathbf{V}(g(p,t)), \nabla_{g(p,t)} \dd(\cdot, \partial M) \right\rangle \geq 1/2 >0.$$ 
\end{itemize}
As a result, we also get that
$\dd(g_\varepsilon(p),\partial M) >0$ for all $\varepsilon>0$.

For the third point, we write $K := \sup_{p \in M} \| d_p \mathbf{V} \|_{op}< \infty$, since $\mathbf{V}$ is $\mathcal{C}^1$ and $M$ compact. 
Let $v \in T_pM$ be a unit vector, and $\gamma$ be a path such that $\gamma(0)=p$ and $\gamma'(0) = v$. For a fixed $t$ and $u \leq \varepsilon$, consider  $
f(u) := \| g(\gamma(t),u) - g(p,u) \|^2.
$
Then
\begin{align*}
|f'(u)| & = 2 \left | \left\langle g(\gamma(t),u)-g(p,u),\mathbf{V}(g(\gamma(t),u)) - \mathbf{V}(g(p,u)) \right\rangle \right | \\
        & \leq 2 K f(u).
\end{align*}
Since $f(0) = \| \gamma(t) - p \|^2$, we deduce that $f(u) \leq \| \gamma(t) - p \|^2 e^{2Ku}$, so that 
\[
\| g(\gamma(t),u) - g(p,u) \| \leq \| \gamma(t) - p\| e^{Ku}.
\]
But since 
\[
g_\varepsilon(\gamma(t)) - \gamma(t) = \int_0^\varepsilon \mathbf{V}(g(\gamma(t),u))\dd u, 
\]
we have
\begin{align*}
g_\varepsilon(\gamma(t)) - g_\varepsilon(p) =  tv + o(t) + \int_0^\varepsilon \left (\mathbf{V}(g(\gamma(t),u)) - \mathbf{V}(g(p,u)) \right )\dd u. 
\end{align*}
Thus
\begin{align*}
\left \| \frac{g_\varepsilon(\gamma(t)) - g_\varepsilon(p)}{t} - v \right \| &\leq o(1) + K \varepsilon e^{K\varepsilon}\|{\gamma(t)-p}\|/t.  
\end{align*}
Letting $t \to 0$, we get that $\|d_p g_\varepsilon - Id_{T_pM} \| \leq K \varepsilon e^{K\varepsilon}$, since $\|{\gamma(t)-p}\|/t \to \norm{v} = 1$.
\end{proof}

The two following results guarantee that for all $p \in M$, there exists a ball with large enough radius with center close to $p$ that does not hit $\partial M$.

\begin{lemma}\label{lem:boule_interieure}
Assume that $\partial M \neq \emptyset$. Let $q \in \partial M$ and $0<t \leq \frac{\tau_M}{8}\wedge \frac{\tau_{\partial M}}{2}$. Then there exists $p_t \in \Int(M)$ such that
\begin{itemize}
\item $\|p_t - q\| \in [t -  4t^2/\tau_M,t +  4t^2/\tau_M]$,
\item $\B\bigl(p_t,t - 4t^2/\tau_M\bigr) \cap \partial M = \emptyset$.
\end{itemize}
\end{lemma}
\begin{proof}[Proof of \Cref{lem:boule_interieure}]
Let $\eta_q$ be the outward-pointing unit normal vector of $M$ at $q$. Denote by $q_t := q -t\eta_q$, and $p_t := \pi_M(q_t)$. 
Note that $\dd(q_t,M) \leq t < \tau_M$, so that $p_t$ is well-defined.

Let us first prove that $p_t \notin \partial M$. For this, if we assume that $p_t \in \partial M$, then $p_t = \pi_{\partial M}(q_t)$ and, since $(q_t - q)\in N_q \partial M$ with $\|q_t - q\| < \tau_{\partial M}$, $p_t =  \pi_{\partial M}(q_t) = q$. But as $p_t =q$, we get $\pi_M(q_t) = q$, with $\|q_t - q\| < \tau_M$. 
Thus, we conclude that $q_t-q = -t\eta_q \in Nor(q,M)$, which is a contradiction. Therefore, we do have $p_t \notin \partial M$ for $0<t < \tau_M \wedge \tau_{\partial M}$.

Now, assume that $t \leq \frac{\tau_M}{8}\wedge \frac{\tau_{\partial M}}{2}$. For some unit vector $u_{p_t} \in (T_{p_t} M)^{\perp}$, it holds
\[
\|p_t - q_t \| = \left\langle p_t - q_t ,u_{p_t} \right\rangle.
\]
Since $\|q_t - q\| = t \leq \tau_M /2$, \cite[Theorem~4.8~(8)]{Federer59} entails that $\|p_t - q\| = \| \pi_M(q_t) - \pi_M(q)\| \leq \tau_M t/(\tau_M-t) \leq 2t$. From \Cref{prop:tangent_variation_geodesic}, we deduce that 
$$
\angle( T_{p_t}M^\perp,T_qM^\perp) = \angle( T_{p_t}M,T_qM) \leq 4t/\tau_M
.
$$
Hence, there exists $u_q \in (T_q M)^{\perp}$ such that $\norm{u_q - u_{p_t}} \leq 4t/\tau_M$. It follows that
\begin{align*}
\norm{p_t-q_t} \leq \left\langle p_t - q_t, u_q \right\rangle + \frac{4t}{\tau_M} \norm{p_t - q_t},
\end{align*}
and thus, since $\eta_q \in T_q M$ and $u_q \in (T_q M)^\perp\subset Nor(q,M)$, we can write
\begin{align*}
\frac{1}{2} \norm{p_t - q_t} 
&\leq 
\left(1- \frac{4t}{\tau_M}\right) \norm{p_t - q_t} 
\\
&\leq 
\left\langle p_t - q_t ,u_q \right\rangle 
\\
&= 
\left\langle p_t - q - t \eta_q, u_q \right\rangle 
\\
&= 
\left\langle p_t - q,u_q \right\rangle
\\
&\leq
\frac{\norm{p_t-q}^2}{2 \tau_M} 
\\
&\leq 
\frac{2t^2}{\tau_M}
,
\end{align*}
where the last but one inequality follows from \cite[Theorem~4.18]{Federer59}.
As $\norm{q_t - q} = t$, triangle inequality then yields $\|p_t - q\| \in [t - \frac{4t^2}{\tau_M},t + \frac{4t^2}{\tau_M}]$. At last, since $\eta_q \in (T_q \partial M)^\perp$ and $t < \tau_{\partial M}$, $\Bopen(q_t,t) \cap \partial M = \emptyset$. Noting that $\B(p_t, t - \frac{4t^2}{\tau_M}) \subset \Bopen(q_t,t)$ concludes the proof.
\end{proof}

\begin{corollary}\label{cor:sub_balls}
For all $r \leq \frac{\tau_M}{32} \wedge \frac{\tau_{\partial M}}{3}$ and $x \in M$, there exists $x' \in \B(x,3r/4) \cap M$ such that $\B(x', r/4) \cap \partial M = \emptyset$.
\end{corollary}

\begin{proof}[Proof of \Cref{cor:sub_balls}]
Let us write $\Delta := \dd(x,\partial M)$, with the convention $\dd(x,\emptyset) = + \infty$.
If $\Delta >r/2$, then taking $x' := x$ gives the result directly.
We shall now assume that $\Delta \leq r/2$.
Denote by $q := \pi_{\partial M}(x)$ and $q_t := q - t \eta_q$, where $t>0$ and $\eta_q$ is the unit outward-pointing vector of $M$ at $q$. 

Write $v:= \pi_{Tan(q, M)}(x-q)$. 
Since $x-q \in (T_q \partial M)^\perp$ and that $\pi_{(T_q \partial M)^\perp}(Tan(q, M)) = \R_- \eta_q$ (see \Cref{prop:cones_of_manifold}), we can write $v = - \ell \eta_q$ for some $\ell \geq 0$.
Thus, we may decompose
\[
x-q = - \ell \eta_q + u, 
\]
with $u \in Nor(q, M)$ and $\|u\| = \dd(x-q,Tan(q,M)) \leq \Delta^2/(2\tau_M)$, from \cite[Theorem~4.18]{Federer59}. 
From this decomposition, reverse triangle inequality yields
\begin{align*}
|\ell - \Delta | 
&=
| \norm{-\ell \eta_q} - \norm{x-q} |
\\
&\leq 
\norm{u}
\\
&\leq
\Delta^2/(2\tau_M)
.
\end{align*}
We hence deduce that 
$
\norm{x-q_\Delta} 
\leq
|\ell - \Delta| + \norm{u}
\leq 
\Delta^2/\tau_M
.
$

Now, pick $x' := \pi_{M}(q_{\Delta + r/2})$. It is immediate that $\| q_{\Delta + r/2} - q_{\Delta}\| = r/2$. Then, following the proof of \Cref{lem:boule_interieure}, since $\Delta + \frac{r}{2} \leq \frac{3r}{2} < \frac{\tau_M}{2}$, it holds
\[
\| q_{\Delta + r/2} - x' \| \leq \frac{4(\Delta + r/2)^2}{\tau_M}.
\]
These bounds altogether lead to
\begin{align*}
\|x' - x\| 
&\leq
\norm{x'-q_{\Delta+r/2}}
+
\norm{q_{\Delta+r/2}-q_\Delta}
+
\norm{q_\Delta-x}
\\
& \leq 
\frac{4 ( \Delta + \frac{r}{2})^2}{\tau_M} 
+
\frac{r}{2}
+ 
\frac{\Delta^2}{\tau_M} 
\\
& \leq r \left (\frac{1}{8} +  \frac{1}{2} + \frac{1}{128} \right ) \leq \frac{3r}{4}.
\end{align*}
At last, since $\Delta + \frac{r}{2} \leq \tau_{\partial M}/2$ and $(\Delta + \frac{r}{2}) - \frac{4 ( \Delta + \frac{r}{2})^2}{\tau_M} \geq \frac{r}{2}(1-1/6) > r/4$, we have
\begin{align*}
\B\left(x',\frac{r}{4}\right) \cap \partial M 
&\subset 
\Bopen\left(x', \left(\Delta + \frac{r}{2}\right) - \frac{4 ( \Delta + \frac{r}{2})^2}{\tau_M}\right) \cap \partial M 
\\
&\subset 
\Bopen\left(q_{\Delta + \frac{r}{2}}, \Delta + \frac{r}{2}\right) \cap \partial M 
\\
&= 
\emptyset
,
\end{align*}
which concludes the proof.
\end{proof}

\propnormalvectorstability*

\begin{proof}[Proof of \Cref{prop:normal_vector_stability}]
Let $p,q \in \partial M$, with $\Vert p - q \Vert = \kappa (\tau_M \wedge \tau_{\partial M})$, where $\kappa \leq 1/32$. According to \Cref{prop:tangent_variation_geodesic} and \Cref{prop:geodesic_vs_euclidean} (applied with $M$), there exists $u \in T_q M$  such that $\Vert \eta_p - u \Vert \leq 2 \Vert p-q \Vert /\tau_M \leq  2 \kappa$.
Decompose $u$ as
\begin{align*}
u = \alpha \eta_q + v_q,
\end{align*}
where $v_q \in T_q \partial M$. We may bound $\Vert v_q \Vert$ as follows. Let $w_q \in T_q \partial M$ with $\Vert w_q \Vert = 1$ be fixed. Using \Cref{prop:tangent_variation_geodesic} and \Cref{prop:geodesic_vs_euclidean} again (but applied with $\partial M$), let $w_p \in T_p \partial M$ be such that $\Vert  w _ p - w _q \Vert \leq 2 \|p-q\|/\tau_{\partial M} \leq 2 \kappa$. We may write
\begin{align*}
\left \langle w_q , v_q \right \rangle & = \left \langle w_q,u \right\rangle \\
& = \left \langle w_p + (w_q - w_p), \eta_p + (u-\eta_p) \right\rangle \\
& \leq \frac{4(1+\kappa)\|p-q\|}{\tau_M \wedge \tau_{\partial M}}, 
\end{align*}
so that $\Vert v_q \Vert \leq 4(1+\kappa)\|p-q\|/(\tau_M \wedge \tau_{\partial M})$. 

Next, let us prove that $\alpha \geq 0$ by contradiction. 
For this, assume that $\alpha < 0 $, and let $\Delta_0 = (\tau_M \wedge \tau_ {\partial M} )/8$. Proceeding as in the proof of \Cref{lem:boule_interieure} yields that 
\begin{align*}
\dd( q + \alpha \Delta_0 \eta_q,M) \leq \frac{4 \alpha ^2 \Delta_0^2}{\tau_M} \leq \frac{\Delta_0}{2}.
\end{align*} 
On the other hand, since $\eta_p \in Nor(p,M)$, \cite[Theorem~4.8~(12)]{Federer59} asserts that $\Bopen(p + \Delta_0 \eta_p, \Delta_0) \cap M = \emptyset$.
But triangle inequality allows to write
\begin{align*}
& \Bopen\bigl( q + \alpha \Delta_0 \eta_q, \Delta_0 ( 1- 10 \kappa - 4 \kappa(\kappa+1)\bigr) \cap M
\\
& \qquad 
\subset \Bopen\bigl( q + (p-q) + \Delta_0(\eta_p-u) + \Delta_0 \alpha \eta_q + \Delta_0 v_q, \Delta_0\bigr) 
\cap M
\\
& \qquad
= \Bopen(p + \Delta_0 \eta_p,\Delta_0)
\cap M
\\
& \qquad
=
\emptyset
,
\end{align*}
so that we get to
\begin{align*}
\dd( q + \alpha \Delta_0 \eta_q,M)
\geq
(1-10\kappa - 4 \kappa(\kappa+1))\Delta_0
> \Delta_0/2,
\end{align*}
which is the desired contradiction. 
Thus, we have proven that $\alpha \geq 0$. Next, note that
\begin{align*}
1 = \Vert \eta_p \Vert & \leq \Vert \eta_p - u \Vert + \Vert u \Vert \\
& \leq \alpha + \Vert v_q \Vert + 2 \kappa, 
\end{align*}
 so that $\alpha \geq 1- 2 \kappa - 4\kappa(1+\kappa) \geq 1/2$. Further, we may write
 \begin{align*}
(1- \alpha)^2 + 2 \alpha \left ( 1- \left\langle\eta_p, \eta_q \right\rangle \right )  & = \Vert \eta_p - \alpha \eta_q \Vert^2 \\ & \leq \left ( \Vert \eta_p - u \Vert + \Vert v_q \Vert \right )^2 \\
 & \leq \left ( \frac{2+4(1+\kappa)}{\tau_M \wedge \tau_{\partial M}} \right )^2 \Vert p-q\Vert^2,
 \end{align*}
 that leads to
 \begin{align*}
 \Vert \eta_p - \eta_q \Vert^2 = 2 \left ( 1- \left\langle\eta_p, \eta_q \right\rangle \right ) & \leq \left ( \frac{2+4(1+\kappa)}{\tau_M \wedge \tau_{\partial M}} \right )^2 \frac{\Vert p-q\Vert^2}{\alpha} \\
 & \leq 2 \left ( \frac{2+4(1+\kappa)}{\tau_M \wedge \tau_{\partial M}} \right )^2 \Vert p-q\Vert^2,
 \end{align*}
 hence the result.
\end{proof}

\subsection{Projections and normals}\label{tecsec:projections}

\propprojectionsestimtangentM*

\begin{proof}[Proof of \Cref{prop:projectionsestimtangentM}]
Let $(y-x)^{T_x}$ and $(y-x)^{\perp_x}$ be the orthogonal projections of $y-x$ onto $T_xM$ and $(T_x M)^{\perp}$ respectively.
Since $\angle(T_x M,T) \leq \theta$, we have 
\begin{align*}
\norm{(y-x)^{{\perp}}}
&\leq 
\norm{((y-x)^{\perp_x})^{{\perp}}}+\norm{((y-x)^{T_x})^{{\perp}}}
\\
&\leq
\norm{(y-x)^{\perp_x}}+\theta \norm{(y-x)^{T_x}}
\\
&\leq
\frac{\norm{y-x}^2}{2\tau_{\min}}
+
\theta \norm{y-x}
,
\end{align*}
where the last line comes from \cite[Theorem 4.18]{Federer59}. This proves the first inequality.
The second one follows from the first one and triangle inequality.
\end{proof}

We now move to the proof of \Cref{prop:existencesestimnorM}, which we split into two intermediate results.

\propexistencesestimnorM*

\begin{proof}[Proof of \Cref{prop:existencesestimnorM}]
This is a straightforward consequence of \Cref{prop:projectionsestimnorM} and \Cref{prop:estimeta}.
\end{proof}

The following two results imply \Cref{prop:existencesestimnorM}. 
First, \Cref{prop:projectionsestimnorM} ensures that estimates of tangent spaces at boundary points contain a normal vector to $\partial M$. Second, \Cref{prop:estimeta} ensures that this normal vector is close to the unit outward-pointing vector at the considered boundary point.

\begin{proposition}\label{prop:projectionsestimnorM}
Assume that $\partial M \neq \emptyset$. Let  $x \in \partial M$ and ${T} \in \Gr{D}{d}$ be such that
$\angle(T_xM,T) < 1$.
Then $T\cap Nor(x,M)$ is a half-line: it contains a unique unit vector $\eta$.

Furthermore, if $y \in \partial M$ and 
$(y-x)^{\eta}$ denotes the orthogonal projection of $(y-x)$ onto $\Span(\eta)$, we have
\[
\|(y-x)^{\eta}\|\leq \frac{\norm{y-x}^2}{2\tau_{\partial M}}
.
\]
\end{proposition}
\begin{proof}[Proof of \Cref{prop:projectionsestimnorM}]
Since $\angle(T_x M,T) < 1$, for all $z \in \R^D\setminus \set{0}$, 
$$\norm{(\pi_{{T}} + \pi_{{T_x M}^\perp})(z)} =\norm{z - (\pi_{{T}} - \pi_{T_x M} )(z)} \geq (1-\angle(T_x M,T)) \norm{z} > 0
.
$$ 
Hence, $\pi_{{T}} + \pi_{{T_x M}^\perp}$ has full rank, which means that $\R^D = T + T_x M^\perp \subset T + N_x M$.  Furthermore, $\dim({T}) + \dim(N_x M) = D+1$ entails that ${T} \cap N_x M = \R u$
for some $u \neq 0$. 
We may thus decompose $u$ as $u=u^{t_x} + u^{\eta_x} + u^{\perp_x}$, where $u^{t_x} = \pi_{{N_xM}^\perp}(u)$, $u^{\perp_x} = \pi_{{T_x M}^\perp}(u)$, and $u^{\eta_x} = \pi_{N_x M \cap T_x M}(u)$. 
Since $u \in N_x M$, we have $u^{t_x}=0$,
and the angle bound $\angle(T_x M,{T}) < 1$ yields that $\|u^{\eta_x}\| \geq \|u\|(1-\angle(T_x M,{T}))>0$.
As a result, $\eta := \sign(\left\langle u, \eta_x \right\rangle) u$ provides us with the announced unique unit $\eta \in T \cap Nor(x,M)$.

Now, the fact that $\eta \in Nor(x,M) \subset (T_x \partial M)^\perp$ allows to write
\begin{align*}
\norm{(y-x)^\eta}
&=
|\inner{y-x}{\eta}|
\\
&=
|\inner{\pi_{(T_x \partial M)^\perp}(y-x)}{\eta}|
\\
&\leq
\norm{\pi_{(T_x \partial M)^\perp}(y-x)}
\\
&\leq
\frac{\norm{y-x}^2}{2\tau_{\partial M}}
,
\end{align*} 
where the last inequality follows from the reach condition on $\partial M$ and \cite[Theorem~4.18]{Federer59}.
\end{proof}

\begin{proposition}\label{prop:estimeta}
Assume that $\partial M \neq \emptyset$. 
Let $x\in \partial M$ and $T \in \Gr{D}{d}$ be such that  $\angle(T_xM,T) \leq \theta<1$. 
Write $\eta$ for \emph{the} unit vector of $Nor(x,M)\cap T$ (\Cref{prop:projectionsestimnorM}). 
Then,
\begin{equation*}
\norm{\eta-\eta_x}
\leq
\sqrt{2} \theta
.
\end{equation*}
\end{proposition}

\begin{proof}[Proof of \Cref{prop:estimeta}]
Since $\eta \in  Nor(x,M)$, $\eta^{t_x} = 0$. 
Furthermore, the angle condition yields that 
$\|\eta^{\perp_x}\| \leq \theta \|\eta\|$. 
We may thus decompose $\eta= \inner{\eta}{\eta_x} \eta_x + \beta u$ for some unit $u\in (\eta_x)^{\perp}$ and $|\beta|\leq \theta$. 
In particular, $|\inner{\eta}{\eta_x}|\geq \sqrt{1-\theta^2}$.
But since $\eta \in Nor(x,M)$, $\inner{\eta}{\eta_x} \geq 0$, so that in fact, $\inner{\eta}{\eta_x}\geq \sqrt{1-\theta^2}$. Finally, as $\eta$ and $\eta_x$ are both unit vectors, we get
\begin{align*}
\norm{\eta-\eta_x}
&=
\sqrt{2}\sqrt{1-\inner{\eta}{\eta_x}}
\leq
\sqrt{2}\sqrt{1-\sqrt{1-\theta^2}}
\leq
\sqrt{2} \theta
.
\qedhere
\end{align*}
\end{proof}

Next, we state a simple lemma that will be useful for describing boundary balls.

\begin{lemma}\label{lem:boulesvidesautourdubord}
Assume that $\partial M \neq \emptyset$. Let $r <\tau_{\min}$, $x\in \partial M$ and $u\in N_x\partial M$ be such that $\langle \eta_x,u \rangle \geq 0$. Then $\B(x+ru,r)\cap M=\{x\}$
\end{lemma}
\begin{proof}[Proof of Lemma \ref{lem:boulesvidesautourdubord}]
As $u \in N_x \partial M$ and $\langle \eta_x,u \rangle \geq 0$, \Cref{prop:cones_of_manifold} yields that $u \in Nor(x,M)$, so that \cite[Theorem~4.8~(12)]{Federer59} asserts that $x$ is the unique projection of $x + ru$ onto $M$.
\end{proof}

The following result provides a quantitative bound on the metric distortion induced by projecting $M$ locally onto (approximate) tangent spaces.

\begin{proposition}\label{lem:piTbijetcontract}
Let $x\in M$ and $T \in \Gr{D}{d}$ be such that $\angle(T_xM,T)\leq \theta$.
Then, for all $y,z \in M\cap \B(x,\tau_{\min}/4)$, we have
$$\left(6/10-\theta \right)\norm{y-z}\leq \norm{\pi_T(y)-\pi_{T}(z)}\leq \norm{y-z}. $$
In particular, if $\theta \leq 1/2$, then $\pi_T : M\cap \B(x,\tau_{\min}/4) \to \pi_{T}(M\cap \B(x,\tau_{\min}/4))$ is a homeomorphism.
\end{proposition}
\begin{proof}[Proof of \Cref{lem:piTbijetcontract}]
The right hand side inequality is straightforward, since $\pi_T$ is an orthogonal projection.
For the other inequality, combine \Cref{prop:tangent_variation_geodesic} and \Cref{prop:geodesic_vs_euclidean} to get
\begin{align*}
\angle (T,T_yM) 
&\leq 
\angle (T,T_xM) + \angle (T_xM,T_yM)
\\
&\leq 
\theta
+
\frac{\dd_M(x,y)}{\tau_{\min}}
\\
&\leq 
\theta 
+
\left(
1
+
\frac{\norm{y-x}^2}{20\tau_{\min}^2}
\right)
\frac{\norm{y-x}}{\tau_{\min}}
\\
&\leq 
\theta + (1+ 1/320)\frac{\norm{y-x}}{\tau_{\min}}.
\end{align*}
Thus,  
\Cref{prop:projectionsestimtangentM} applied at $y$ and $z$ entails
\begin{align*}
\norm{\pi_{T}(y)-\pi_{T}(z)}
&\geq
\left(1-\{\theta + (1+ 1/320)\norm{y-x}/\tau_{\min}\}-\frac{\norm{y-z}}{2\tau_{\min}}\right)
\norm{y-z}
\\
&\geq
\left(6/10-\theta \right)
\norm{y-z}
,
\end{align*}
which concludes the proof.
\end{proof}

For $q \in M$, the following result characterizes the boundary of $\pi_{T}( M \cap \B(q,r)-q)$, when seen as a subset of $T \cong \R^d$.

\begin{lemma}\label{lem:projection_boundary_commutation}
Let $0\leq r\leq \tau_{\min}/16$. Then for all $q\in M$ and $T \in \Gr{D}{d}$ such that $\angle (T_qM,T) \leq \theta \leq 1/8$,
$$
\partial \pi_{q+T}\bigl( M \cap \B(q,r)\bigr)=\pi_{q+T}\bigl(\partial M \cap \B(q,r)\bigr)\cup \pi_{q+T}\bigl( M \cap \partial \B(q,r)\bigr)
.
$$
\end{lemma}

\begin{proof}[Proof of \Cref{lem:projection_boundary_commutation}]
As preliminary remarks, first note that since $M \cap \B(q,r)$ is compact and $\pi_{q+T}$ is continuous, we have 
\[
\overline{\pi_{q+T}(M \cap \B(q,r))} = \pi_{q+T}(M \cap \B(q,r)).
\]
Furthermore, for all $p \in \B(q,r)$, \Cref{prop:tangent_variation_geodesic} and \Cref{lem:dist_geod_dist_eucl} yield that $\angle(T_pM, T) \leq   1/4$. 
We recall that $\Int(M)=M\setminus\partial M$.

\begin{itemize}[leftmargin=*]
\item[\textbf{Step 1:}]
First, we prove that $\pi_{q+T} \bigl(\Int(M) \cap \Bopen(q,r) \bigr ) \subset \bigl ( \pi_{q+T}(\B(q,r) \cap M) \bigr)^{\mathrm{o}}$.
\\
\noindent
For this, let $p \in \Int(M) \cap \Bopen(q,r)$ be fixed. 
Let $\rho_M \in (0 ,\min\set{r-\norm{p-q},\dd(p,\partial M)})$ (with the convention $\dd(p, \emptyset) = + \infty$), so that in particular, $M \cap \Bopen(p,\rho_M) \subset \Int(M) \cap \Bopen(q,r)$. 
According to~\cite[Lemma 1]{Aamari19b}, there exists $0<r_2 \leq \tau_{M}/8$ such that 
\[
\exp_p: \Bopen_{T_p M}(0,r_2) \longrightarrow \Bopen(p,\rho_M) \cap \Int(M)
\]
is a diffeomorphism onto its image, and can be decomposed as $\exp_p(v) = p + v + N_p(v)$, with $N_p(0)=0$, $d_0 N_p=0$, $\|d_v N_p \|_{op} \leq 5/(4\tau_M)$. 
We now consider the map $g$ defined as
\begin{align*}
g \colon \Bopen_T(0,r_2) &\to \Bopen(p,\rho_M) \cap \Int(M)
\\
u &\mapsto \exp_p(\pi_{T_p M}(u))
\end{align*}
Note that, since $\angle(T_p M,T) \leq 1/4$, $\pi_{T_p M} : \Bopen_T(0,r_2) \rightarrow \Bopen_{T_p M}(0,r_2)$ is a diffeomorphism  onto its image that satisfies $\|u - \pi_{T_p M}(u)\| \leq \|u\|/4$ for all $u \in \Bopen_{T}(0,r_2)$.
In particular, $\pi_{T_p M}$ is injective on $T$, and hence so is $g$ on its domain.
As a result, for all $u_1, u_2 \in \Bopen_T(0,r_2)$, 
\begin{align*}
g(u_1) - g(u_2) 
&= 
(u_1 - u_2) + (\pi_{T_p M}(u_1 - u_2)-(u_1-u_2)) 
\\
&~~~~~~+ N_p(\pi_{T_p M}(u_1)) - N_p(\pi_{T_p M}(u_2)).
\end{align*}
We may thus bound
\begin{align*}
\| g(u_1) - g(u_2) - (u_1 - u_2) \| 
&\leq 
\frac{1}{4} \|u_1 - u_2\| + 5r_2/(4\tau_{\min}) \|u_1 - u_2\| 
\\
& \leq 
\frac{1}{2}\|u_1 - u_2\|
.
\end{align*}
Let now $f: \Bopen_T(0,r_2) \rightarrow \Bopen_T(0,\rho_M)$ be defined as $f(\cdot) := \pi_{q+T} \circ (g(\cdot) - p)$. 
By composition and \Cref{lem:piTbijetcontract}, $f$ is clearly injective. Moreover, for all $u_1, u_2 \in \Bopen_T(0,r_2)$,
\[
\frac{1}{2} \| u_1 - u_2\| \leq \|f(u_1) - f(u_2)\| \leq \frac{3}{2} \|u_1 - u_2\|,
\]
since $\pi_T(u_1-u_2) = u_1-u_2$ and $\| \pi_{T} \left (g(u_1) - g(u_2) - (u_1 - u_2) \right ) \| \leq \|u_1 - u_2\|/2$.
Thus, $f: \Bopen_T(0,r_2) \rightarrow f(\Bopen_T(0,r_2))$ is a homeomorphism, which ensures that $f(\Bopen_T(0,r_2))$ is an open subset of $T$ that contains $0 = f(0)$.
But by construction,
$$\pi_{q+T}(p) + f(\Bopen_T(0,r_2)) \subset \pi_{q+T}(\Bopen(p,\rho_M) \cap \Int(M)),$$ 
which shows that $\pi_{q+T}(p) \in \left ( \pi_{q+T}(\B(q,r) \cap M) \right )^{\mathrm{o}}$, and concludes the first step.

\item[\textbf{Step 2:}] 
Next, we show that no element of $\pi_{q+T}((\partial M \cap \B(q,r)) \cup (M \cap \Sphere(q,r)))$ 
can belong to the interior set $\left ( \pi_{q+T}(\B(q,r) \cap M) \right )^{\mathrm{o}}$. 
\begin{itemize}[leftmargin=*]
\item
If $\partial M \neq \emptyset$, let $p \in \partial M \cap \B(q,r)$ be fixed. Striving for a contradiction, assume that $\pi_{q+T}(p) \in \pi_{q+T}(M \cap \B(q,r))^\mathrm{o}$. 
In particular, for $\delta > 0$ small enough, $\pi_{q+T}(p + \delta \eta_p) \in \pi_{q+T}(\B(q,r) \cap M)$.
Without loss of generality, we shall pick $\delta \in (0, \tau_{\min}/16)$ small enough so that $p + \delta \eta_p \in \B(q,r)$.

Then there exists $p'\in \B(q,r) \cap M$ such that $\pi_{q+T}(p') = \pi_{q+T}(p+\delta \eta_p)$, or equivalently, $\pi_T(p'-p)=\delta \pi_T(\eta_p)$. 
Consider $v :=p'-p - \delta \eta_p$.
By construction, $\pi_T(v) = 0$, so that $v\in T^{\perp}$, and its norm is at most
\[
\norm{v}
\leq
\norm{p'-p} + \norm{\delta \eta_p}
\leq
2r + \delta
\leq
3\tau_{\min}/8.
\]
Furthermore, $v \neq 0$, as otherwise this would mean that $p + \delta \eta_p = p' \in \B(q,r) \cap M \subset M$, which is impossible since $\dd(p+\delta \eta_p,M) = \delta$ from \cite[Theorem~4.8~(12)]{Federer59}.
We may now decompose $v$ as $v = v_1 + v_2$, with $v_1 \in T_p M$ and $v_2 \in T_p M^\perp$.
\begin{itemize}[leftmargin=*]
\item
On one hand, the angle bound $\angle(T,T_p M) \leq 1/4$ and $v \in T_p M^\perp$ yield $\norm{v_1} \leq \|v\|/4$.
\item
Furthermore, $\delta\leq \tau_{\min}/16$ ensures that $\|v_2\| \leq \|v\| \leq 3 \tau_{\min}/8 < \tau_{M} - \delta$. 
Let us now consider $s:= p + \delta \eta_p + v_2$. 
As $\delta \eta_p + v_2 \in Nor(p,M)$ and $\|\delta \eta_p + v_2\| < \tau_M$, \cite[Theorem~4.8~(12)]{Federer59} asserts that $\pi_M(s) = p$ and $\dd(s,M) = \|\delta \eta_p + v_2\|$. 
But on the other hand, $s + v_1= p' \in M$, so clearly $\|v_1\| \geq \dd(s,M)$. Therefore,
\begin{align*}
\|v_1\|^2 
&\geq 
\|\delta \eta_p + v_2\|^2
\\
&=
\delta^2 + \|v_2\|^2 
\\
&= 
\delta^2 + \|v\|^2 - \|v_1\|^2
\\
&\geq
\|v\|^2 - \|v_1\|^2
,
\end{align*}
and thus $\|v_1\| \geq \|v\|/\sqrt{2}$.
\end{itemize}
The last two items contradicting each other, we finally obtain that $p \notin \pi_{q+T}(M \cap \B(q,r))^\mathrm{o}$.

\item
Let now $p \in \partial \B(q,r) \cap M$ be fixed.
Striving for a contradiction, let us assume that $\pi_{q+T}(p) \in \pi_{q+T}(M \cap \B(q,r))^\mathrm{o}$. 
This implies in particular that for all $\delta <1$ small enough, $\pi_{q+T}(p + \delta(p-q)) \in \pi_{q+T}(\B(q,r) \cap M)$. 
Then there exists $v \in T^\perp$ such that $p + \delta(p-q) + v \in M \cap \B(q,r)$.
Denote by $v_2 = \pi_{T_p M^\perp}(v)$. 
Since $\angle(T_p M,T) \leq 1/4$, we have $\|v\| \geq 3 \|v_2\|/4$. 
On the other hand, since $p + \delta(q-p) + v \in M$, we have 
\begin{align*}
\| \pi_{T_p M^\perp}( \delta(p-q) + v) \|
&=
\dd\bigl((p + \delta(p-q) + v) - p,T_p M\bigr)
\\
&\leq 
\frac{\| \delta (p-q) + v \|^2}{2 \tau_M}
\\
&\leq
\frac{\delta^2 r^2 +\| v \|^2}{ \tau_M}
,
\end{align*}
from \cite[Theorem~4.18]{Federer59}.
And noting that
\begin{align*}
\norm{
\pi_{T_p M^\perp}( \delta(p-q) + v)
}
&=
\norm{
\delta\pi_{T_p M^\perp}(p-q)
+
v_2
}
\\
&\geq 
\norm{v_2}
-
\delta \dd(q-p,T_p M)
\\
&\geq
\frac{3\|v\|}{4} - \frac{\delta r^2}{2 \tau_M}
,
\end{align*}
we obtain
\begin{align}
\|v\|
& \leq 
\frac{4}{3}
\left(
\frac{\delta r^2}{2 \tau_M}
+
\frac{\delta^2 r^2 +\| v \|^2}{ \tau_M}
\right)
\leq
2\left(
\frac{\delta r^2}{2 \tau_M}
+
\frac{\delta^2 r^2 +\| v \|^2}{ \tau_M}
\right)
.
\label{eq:majoration_v_bord_projete}
\end{align}
On the other hand, since $p + \delta(p-q) + v \in \B(q,r)$, we have $\|(1+\delta)(p-q) + v\|^2 \leq r^2$, and therefore
\[
(2 \delta + \delta^2)r^2 + \|v\|^2 - 2 (1+ \delta) r \|v\| \leq 0,
\]
But according to \eqref{eq:majoration_v_bord_projete}, this last inequality yields
\begin{align*}
(2 \delta + \delta^2)r^2
&
+ \|v\|^2 - 2 (1+ \delta) r \|v\|
\\
&\geq 
(2 \delta + \delta^2)r^2 + \|v\|^2  -{4(1+ \delta)r} 
\left(
\frac{\delta r^2}{2 \tau_M}
+
\frac{\delta^2 r^2 +\| v \|^2}{ \tau_M}
\right)
\\
&=
\norm{v}^2
\left(
1- 4(1+\delta) \frac{r}{\tau_M}
\right)
+
r^2
\left(
(2 \delta + \delta^2)
-
4(1+\delta)\left\{
\frac{\delta r}{2 \tau_M}
+
\frac{r \delta^2}{\tau_M}
\right\}
\right)
,
\end{align*}
and since $r \leq \tau_M / 16$ and $\delta \in (0, 1]$, we finally get
\begin{align*}
(2 \delta + \delta^2)r^2
+ \|v\|^2 - 2 (1+ \delta) r \|v\|
&
\geq
\frac{\norm{v}^2}{2}
+
r^2
\left(
(2 \delta + \delta^2)
-
\delta(1+\delta)
\left\{\frac{1}{8} + \frac{1}{4}\right\}
\right)
\\
&\geq
\frac{\norm{v}^2}{2}
+
r^2\delta
\\
&>
0
\end{align*}
which is the desired contradiction.
That is, we have $\pi_{q+T}(p) \notin \pi_{q+T}(M \cap \B(q,r))^\mathrm{o}$, as announced.
\end{itemize}

\item[\textbf{Conclusion:}]
Putting everything together, we deduce that 
\begin{align*}
\pi_{q+T}((\partial M \cap \B(q,r)) \cup (M \cap \partial \B(q,r))) 
&= 
\overline{\pi_{q+T}(M \cap \B(q,r)} \setminus \pi_{q+T}(M \cap \B(q,r))^\mathrm{o} 
\\
&= \partial \pi_{q+T}(M \cap \B(q,r)
,
\end{align*}
which is the announced result.
\end{itemize}
\end{proof}

\subsection{Structure of balls on manifolds with boundary}
\label{tecsec:ball-structure}

Using Lemma \ref{lem:projection_boundary_commutation}, we are now able to derive the two key results on the structure of $\pi_{{T}}(\B(x,R_0)-x)$. This structure depends on whether $x$ is either near or far from $\partial M$.
We start with the case where $x$ is an interior point.

\lemGeoInterbleupartone*

\begin{proof}[Proof of \Cref{lem:GeoInter_bleupart1}]
Let $z'$ be in $\Bopen\left(x,4\min\set{R, \dd(x,\partial M)}/5 \right) \cap (x+T)$, and assume for contradiction that 
$z'\notin  \pi_{x+T}(\B(x,R)\cap M)$. 
Then by connectedness, there exists $z\in[x,z']$ such that $z \in \partial \pi_{x+T}(\B(x,R)\cap M)$. 
\begin{itemize}[leftmargin=*]
\item
Note that, since $\Bopen\left(x,4\min\set{R, \dd(x,\partial M)}/5 \right) \cap (x+T)$ is convex and contains $\{x,z'\}$, we have $z \in \Bopen\left(x,4\min\set{R, \dd(x,\partial M)}/5 \right) \cap x+T$. 
\item
According to \Cref{lem:projection_boundary_commutation}, we can write $z=\pi_{x+T}(y)$ with
$ y \in \partial \B(x,R)\cap M$ or  $y\in \B(x,R)\cap \partial M $. 
Therefore, we have either $\norm{y-x}= R$ or $\norm{y-x}\geq \dd(x,\partial M)$, which entails 
$\norm{y-x}\geq \min\set{R,\dd(x,\partial M)}$. Applying \Cref{prop:projectionsestimtangentM}  gives that
\begin{align*}
\norm{x-z} 
&=
\norm{\pi_T(x) - \pi_T(z)}
\\
& \geq 
\min\set{R,\dd(x,\partial M)}\left ( 1-\theta-\frac{\norm{x-y}}{2\tau_{\min}} \right ) 
\\
& \geq  
\frac{27}{32}\min\set{R,\dd(x,\partial M)}
\\
&\geq
\frac{4}{5}\min\set{R,\dd(x,\partial M)},
\end{align*}
leading to $z \notin \Bopen\left(x,4\min\set{R,\dd(x,\partial M)}/5 \right)$, and hence a contradiction.
\end{itemize}
It follows that
$
\Bopen\left(x,4\set{R, \dd(x,\partial M)}/5 \right) \cap (x+T)
\subset 
\pi_{x+T}(\B(x,R)\cap M)
$
. 
Finally, the closedness of $\pi_{x+T}(\B(x,R)\cap M)$ concludes the proof.
\end{proof}

Next we turn to the case where $x$ is a boundary point.

\lemGeoInterbleuparttwo*

\begin{proof}[Proof of \Cref{lem:GeoInter_bleupart2}]

Take $O=x+\alpha \hat{\eta}$ with $|\alpha|=r$. 

We first prove that
$(\B(O,r)\cap (x+T))\cap \partial \pi_{x+T}(\B(x,R)\cap M)=\{x\}$.
For this, consider $z\in \pi_{x+T}( M \cap \B(x,R))\setminus \{x\}$ 
and $y\in  M \cap \B(x,R)$ such that $z=x+(y-x)^{T}=x+(y-x)^{\hat{t}}+(y-x)^{\hat{\eta}}$. Recall that
$(y-x)^{\hat{t}}$ denotes the orthogonal projection of $y-x$ onto  $\hat{\eta}^{\perp}\cap {T}$. We have that
\begin{align*}
\norm{O-z}^2 & = \left (\norm{(y-x)^{\hat{\eta}}}\pm|\alpha| \right )^2+\norm{(y-x)^{\hat{t}}}^2
\\ 
            &\geq
             \left (\norm{(y-x)^{\hat{\eta}}}-|\alpha| \right )^2+\norm{(y-x)^{\hat{t}}}^2
             \\ 
            &=
            r^2+\norm{(y-x)^{{T}}}^2-2r\norm{(y-x)^{\hat{\eta}}}.
\end{align*}
According to \Cref{lem:projection_boundary_commutation}, if $z\in \partial \pi_{x+{T}}( M \cap  \B(x,R))$,  we have either
$z\in \pi_{x+{T}}(M\cap \partial \B(x,R))$, or $z\in \pi_{x+{T}}(\partial M\cap \B(x,R))$. In the first case,  \Cref{prop:projectionsestimtangentM} gives 
\begin{align*}
\norm{O-z}^2& \geq  r^2+\norm{(y-x)^{{T}}}^2-2r\norm{(y-x)^{{T}}} \\
& \geq r^2+\norm{(y-x)^{{T}}}\left(\norm{(y-x)^{{T}}}-2r\right) \\
&\geq r^2+\norm{(y-x)^{{T}}}\left(\frac{27}{32}R-2r\right).
\end{align*}
In the second case, using \Cref{prop:projectionsestimtangentM} and \Cref{prop:projectionsestimnorM} leads to
\begin{equation*}
\norm{O-z}^2 \geq r^2+\norm{y-x}^2\left( \left(\frac{27}{32}\right)^2-\frac{r}{2\tau_{\partial, \min}}\right).
\end{equation*}
In both cases, since $z \neq x$ by assumption, we have $(y-x)^T \neq  0$ and hence $y-x \neq 0$, so that if 
$
r \leq \min \set{2R/5,7\tau_{\partial, \min}/5}
$, 
we have $\norm{O-z}> r$, which entails $z\notin \B(O,r)$. In other words, we have proved that $\B(O,r) \cap \partial \pi_{x+T}(\B(x,R_{0})\cap M)=\{x\}$.  

By connectedness, it follows that if $O\in \{x+O^{\text{in}}, x+ O^{\text{out}}\}$, we have either
$$\B(O,r)\cap (x+{T})\subset \pi_{x+{T}}(\B(x,R_{0})\cap M),$$
or 
$$
    \B(O,r)\cap (x+{T})
    \subset 
    \left(\pi_{x+{T}}(\B(x,R_{0})\cap M\right)^c \cup \{x\}.
$$

Let us now focus on $\B(x+O^{\text{out}},r)\cap (x+ T)$. 
Consider a sequence $x^*_n=x+\eps_n \hat{\eta}$ with $\eps_n>0$ converging to $0$. Suppose that $x^*_n\in \pi_{x+{T}}(\B(x,R)\cap M)$ 
i.e. there exits $x_n\in M$ such that 
$x^*_n-x=(x_n-x)^{{T}}$. By \Cref{prop:projectionsestimtangentM}, we have $\norm{(x_n-x)^{{\perp}}}\leq \eps_n(\theta+1/4)$. 
Let $\Omega=x+r'\hat{\eta}$ with $r'<\min(\tau_{\min},\tau_{\partial,\min})$. On one hand \Cref{lem:boulesvidesautourdubord} 
ensures that $\norm{\Omega-x_n}\geq r'$ and, on the other hand 
$$
    \norm{
        \Omega-x_n
    }^2
    =
    (r'-\eps_n)^2
    +
    \norm{
        (x_n-x)^{{\perp}}
    }^2
    \leq 
    r'^2-2\eps_n r' 
    +
    \eps_n^2
    \left(
        1+(\theta+1/4)^2
    \right)
.
$$
Thus, for $n$ large enough $\norm{\Omega-x_n}^2<r'^2$, which is impossible. Hence,  for $n$ large enough $x_n^*\notin  \pi_{x+ {T} }(\B(x,R)\cap M)$, which proves the right hand side inclusion
\begin{equation*}
\pi_{x+{T}}(\B(x,R)\cap M)\subset \bigl(\B(x+O^{\text{out}},r)^c\cap (x+{T})\bigr)\cup\{x\}.
\end{equation*}

Next, we prove that if $\theta \leq 1/8$, then there exists $x^*\in x+{T} \cap \B(x+O^{\text{in}},r)$ such that
$x^*\in \pi_{x+{T}}(\B(x,R)\cap M)$, and thus $\B(x+O^{\text{in}},r)\cap x+{T}\subset \pi_{x+{T}}(\B(x,R_{0})\cap M)$.
For this, introduce $\eta=\pi_{T_xM}(\hat{\eta})$ and $\eta'=\pi_{{T}}(\eta)$. We clearly have $\norm{\eta}\leq 1$, $\norm{\eta'}\leq 1$,
$\norm{\hat{\eta}-\eta}\leq \theta$ and $\norm{\eta-\eta'}\leq \theta$. In particular, this implies that $\norm{\eta' - \hat{\eta}} \leq 2 \theta< 1$ and $\norm{\eta'} \geq 1-2\theta$.
Hence, decomposing $\eta'=\lambda \hat{\eta} + \mu {v}$, with ${v}\in {T}\cap(\hat{\eta})^{\perp}$ and $\norm{{v}}=1$, we have $\lambda >0$, with
\begin{equation*}
(1-2\theta)^2 \leq \lambda^2+\mu^2\leq 1 \text{ and }
\lambda \geq 1-2\theta.
\end{equation*}
Furthermore, since $\eta\in T_xM$ and that 
$$
\inner{\eta}{\eta_x} \geq 1- \norm{\eta-\eta_x} \geq 1 - \norm{\eta_x - \hat{\eta}} - \norm{\hat{\eta}-\eta}
\geq
1 - \sqrt{2} \theta - \theta
>
0
$$
from \Cref{prop:estimeta}, we get that $\eta \in Nor(x,M)$ from \Cref{prop:cones_of_manifold}, or equivalently that $-\eta \in Tan(x,M)$.
Hence, \cite[Definition~4.3]{Federer59} asserts that there exists a sequence $(x_n)_n\in M\setminus \{x\}$ converging to $x$ such that 
$\norm{\frac{x_n-x}{\norm{x_n-x}}-\frac{-\eta}{\norm{\eta}}}\leq \frac{1}{n}$,
that is
$$x_n=x-\norm{x-x_n}\left( \frac{\eta}{\norm{\eta}}+\frac{1}{n}w_n\right)\text{ with }\norm{w_n}\leq 1. $$
Considering $x_n^*=\pi_{x+{T}}(x_n)$, $w_n^*=\pi_{{T}}(w_n)$, and $\eps_n=\frac{\norm{x-x_n}}{\norm{\eta}}$, we may hence write
\begin{equation*}
x_n^*=x-\eps_n\left(\lambda \hat{\eta}+\mu {v} + \frac{\norm{\eta}}{n}w_n^*\right),
\end{equation*}
so that
\begin{align*}
\norm{x+O^{\text{in}}-x_n^*} & \leq  \norm{(r-\lambda \eps_n)\hat{\eta}+\eps_n \mu {v}}+\frac{\eps_n}{n}\\
	 &\leq \sqrt{r^2-2r\lambda\eps_n+\eps_n^2}+\frac{\eps_n}{n}\\
	 &\leq \sqrt{(r-\lambda \eps_n)^2+\eps_n^2(1-\lambda^2)} +\frac{\eps_n}{n}\\
	 &\leq (r-\lambda \eps_n)+\eps_n\sqrt{1-\lambda^2} +\frac{\eps_n}{n}.
\end{align*}
Since $\lambda \geq 1 - 2 \theta \geq 3/4$, this yields
$$\norm{x+O^{\text{in}}-x_n^*}\leq r-\eps_n
\left(\frac{3}{4}-\frac{\sqrt{7}}{4}+\frac{1}{n} \right).$$
On the other hand, we have
$$\norm{x_n^*-x}\geq \frac{\norm{x_n-x}}{\norm{\eta}} \left(\sqrt{\lambda^2+\mu^2}-\frac{1}{n} \right)\geq\frac{\norm{x_n-x}}{\norm{\eta}} \left(\frac{3}{4}-\frac{1}{n} \right) >0, $$
for $n$ large enough. Thus, for $n$ large enough, $x_n^*\in (x +{T})\cap \B(x+O^{\text{in}},r)$ with $x_n^*\in \pi_{x+{T}}(\B(x,R)\cap M)$ and $x_n^*\neq x$, ensuring that
\begin{equation*}
\B(x+O^{\text{in}},r)\cap (x +{T}) \subset \pi_{x+T}(\B(x,R_{0})\cap M),
\end{equation*}
which is the left hand side inclusion.
\end{proof}

At last, the following consequence of \Cref{lem:GeoInter_bleupart2} will be of particular interest in the proof of \Cref{thm:deterministic}.

\corprojorthobord*

\begin{proof}[Proof of \Cref{cor:projorthobord}]
According to \Cref{prop:existencesestimnorM}, $Nor(x',M) \cap T$ contains a unique unit vector $\eta^*(x')$. By definition of $x^*$ we have 
  \begin{equation}\label{corobleu:inclu1}
  \Bopen_{y+T}(\pi_{y+T}(x),\norm{x^*-\pi_{y+T}(x)})\cap \pi_{y+T}(\partial M\cap \B(x,\tau_{\min}/16))=\emptyset.
 \end{equation}
 Since $\pi_{y+T} = \pi_{x'+T} + \pi_{T^\perp}(y-x')$, \Cref{lem:GeoInter_bleupart2} applied at $x'$ with $R_0=\tau_{\min}/16$ yields
\begin{align*}
  \Bopen_{x'+T}(x'+r_0\eta^*(x'),r_0)\cap \pi_{x'+T}(M\cap \B(x',R_0)) =\emptyset.
 \end{align*}
 Since $\pi_{x'+T} = \pi_{y+T} + \pi_{T^\perp}(x'-y)$, and that for all $p \in x'+T$ and $r>0$, 
 $$
 \Bopen_{x'+T}(p,r) = \pi_{T^\perp}(x'-y) + \Bopen_{y+T}(\pi_{y+T}(p),r)
 ,
 $$
 we deduce that 
  \begin{equation}\label{corobleu:inclu2}
  \Bopen_{y+T}(x^*+r_0\eta^*(x'),r_0)\cap \pi_{y+T}(M\cap \B(x',R_0)) =\emptyset.
 \end{equation}
Now, decompose 
$$
x^*-\pi_{y+T}(x)=\cos \varphi \norm{x^*-\pi_{y+T}(x)} \eta^*(x') + \sin \varphi \norm{x^*-\pi_{y+T}(x)} v
$$ 
with
$v\in \eta^*(x')^{\perp}$ and $\varphi \in [0,2 \pi)$, and consider 
$$
x_t
:=
x^* + t \sin (\pi-\varphi/2) \eta^*(x')+ t\cos (\pi-\varphi/2)  v
,
$$
for $t\geq 0$.
Straightforward calculus yields
\begin{align*}
\begin{cases}
 \norm{x^*+r_0\eta^*(x')-x_t}^2=r_0^2+t^2-2r_0t\sin (\pi-\varphi/2), \\
 \norm{\pi_{y+T}(x)-x_t}^2=\norm{x^*-\pi_{y+T}(x)}^2+t^2+2t\norm{x^*-\pi_{y+T}(x)}\sin (\pi+\varphi/2),\\
 \norm{x-x_t}\leq \norm{x-x^*}+t \text { with }\norm{x-x^*}\leq \dd(x,\partial M)<\tau_{\min}/16.
 \end{cases}
\end{align*}
Suppose, to derive a contradiction, that $\varphi\neq 0$. 
Then for small enough $t$, we have 
\[
x_t \in \Bopen(x,\tau_{\min}/16) \cap \Bopen_{y+T}(x^*+r_0\eta^*(x'),r_0) \cap \Bopen_{y+T}(\pi_{y+T}(x), \norm{x^*-\pi_{y+T}(x)}).
\]
Then, \Cref{corobleu:inclu2} provides $z\in(x_t,\pi_{y+T}(x))$ such that $z\in \pi_{y+T}(\partial M\cap \B(x,\tau_{\min}/16))$.
But since $\norm{z-\pi_{y+T}(x)} < \norm{x^*-\pi_{y+T}(x)}$ by construction,  \Cref{corobleu:inclu1} leads to the desired contradiction.
Hence, $\varphi = 0$, which yields the announced result.
\end{proof}

\subsection{Volume bounds and covering numbers}\label{tecsec:volumes_and_covering}

\lemcovering*

\begin{proof}[Proof of \Cref{lem:covering}]
Let $\varepsilon_1 \leq \frac{\tau_{\min}}{16} \wedge \frac{\tau_{\partial,\min}}{2}$, and $x \in M$. 
As $\X_n \subset M$, the Hausdorff distance between $M$ and $\X_n$ writes as $\dHaus(M,\X_n) = \max_{x \in M} \dd(x,\X_n)$. Furthermore, according to \Cref{cor:sub_balls}, 
\begin{align*}
\mathbb{P} \left ( \max_{x \in M} \dd(x,\X_n) \geq \varepsilon_1 \right ) 
&\leq 
\mathbb{P}  \left ( \max_{\substack{x' \in M \\ \dd(x',\partial M)\geq \varepsilon_1/4}} \dd(x',\X_n) \geq \varepsilon_1/4 \right ) \\
& \leq \frac{16^d}{c_d f_{\min} \varepsilon_1^d} \exp \left ( -n \frac{c_d f_{\min}}{8^d}\varepsilon_1^d \right ),
\end{align*} 
where the second inequality follows as \cite[Lemma~9.1]{Aamari18}. Thus, choosing $\varepsilon_1 = \left ( C_d \frac{\log n }{f_{\min}n} \right )^\frac{1}{d}$, for $C_d$ large enough, yields that 
$\dHaus(M,\X_n) \leq \varepsilon_1$, with probability larger than $1-n^{-3}$. 
\end{proof}

\begin{lemma}[Volume of Intersection of Balls]
\label{teclem:volume_intersection_balls}
Let $0 \leq r' \leq r$, and $O, O' \in \R^d$ that satisfy 
\[
\| O - O' \| = r +r' - h, 
\]
for some $0 \leq h \leq r'$. Then
\[
\Haus^d\bigl(
\B(O,r) \cap \B(O',r')
\bigr)
\geq 
\frac{\omega_{d-1}}{d 2^{\frac{d-1}{2}}}h^{\frac{d+1}{2}}(r')^\frac{d-1}{2}.
\]
\end{lemma}

\begin{proof}[Proof of \Cref{teclem:volume_intersection_balls}]
Let $A:=\partial \B(O,r)\cap [O,O']$,  $B:=\partial \B(O',r')\cap [O,O']$, 
and $\Omega$ be the orthogonal projection of any point of $\partial \B(O,r) \cap \partial \B(O',r')$ onto $[O,O']$.
Also define $a:=\norm{A-\Omega}$, $b:=\norm{B-\Omega}$ and $\ell:=\dd(\Omega,\partial \B(O,r) \cap \partial \B(O',r')) $ (see \Cref{fig:voldessin}).
Let $\mathcal{C}$ (resp. $\mathcal{C}'$) denote the section of cone of apex $B$ (resp. $A$), direction $O-O'$ (resp. $O'-O$), and basis $\B(\Omega,\ell)\cap \left(\Omega + \Span(O'-O)^\perp \right)$. 
\begin{figure}[!htp]
	\centering
	\includegraphics[width=0.6\textwidth]{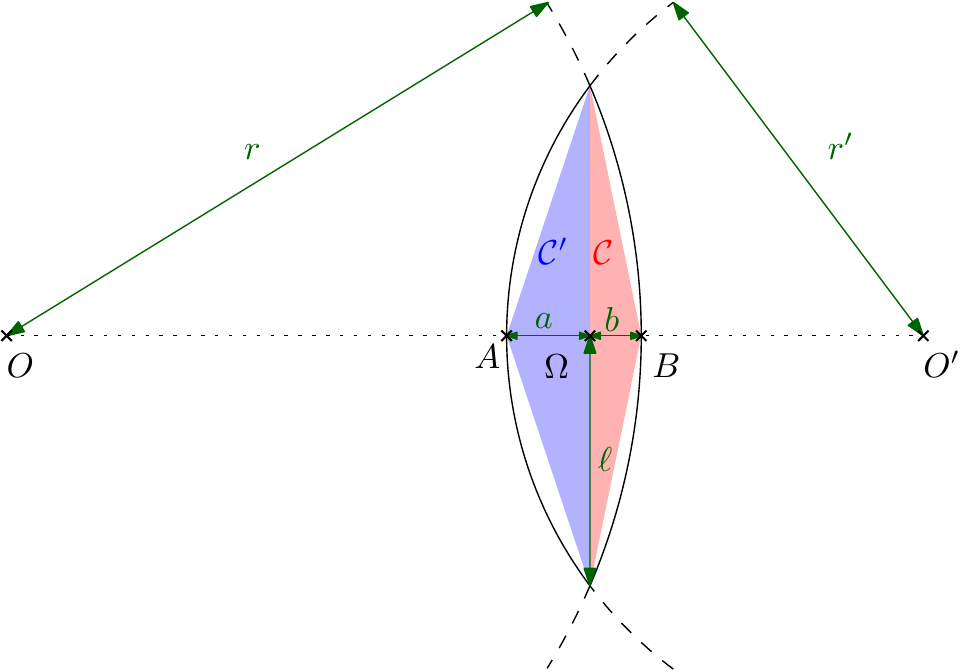}	
	\caption{Layout for \Cref{teclem:volume_intersection_balls}.}
	\label{fig:voldessin}
\end{figure}

By convexity, we have $\mathcal{C},\mathcal{C}' \subset \B(O,r) \cap \B(O',r')$, and since $\mathcal{C} \cap \mathcal{C}'$ is included in a hyperplane, we get
\begin{align}\label{eq:minovol1}
\Haus^d\bigl(
\B(O,r)\cap \B(O',r')
\bigr) 
& \geq 
\Haus^d(\mathcal{C} \cup \mathcal{C}')
\notag
\\
&=
\Haus^d(\mathcal{C})
+
\Haus^d(
\mathcal{C'}
)
\notag 
\\  
& = 
\frac{\omega_{d-1}}{d}\ell^{d-1} (a+b) \notag
\\
&=\frac{\omega_{d-1}}{d}\ell^{d-1}h
.
\end{align}
Furthermore, since $a+b=h$, Pythagoras theorem gives
$$
(r-b)^2+\ell		^2=r^2
\text{~~and~~}
(r'-a)^2+\ell^2=r'^2,
$$
leading to 
$$a=\frac{2rh-h^2}{2(r+r'-h)}=\frac{rh}{r+r'}+\frac{h^2}{r+r'-h}\left( \frac{r}{r+r'}-\frac{1}{2}\right).$$ 
Recalling that $r' \leq r$, we may write 
\begin{equation*}
\frac{rh}{r+r'}\leq a\leq \frac{rh}{r+r'-h}.
\end{equation*}
Finally, since $\ell^2=2r'a-a^2$, we hence obtain
\begin{align*}
\ell^2 & \geq a \left ( 2r' - \frac{rh}{r+r'-h} \right ) \geq \frac{2r'rh}{r+r'} \left( 1 - \frac{rh}{2r'(r+r'-h)} \right )
\geq \frac{r'rh}{r+r'}.
\end{align*}
Combining the equation above with \eqref{eq:minovol1} concludes the proof.\end{proof}

\lemprobaoccupationlosange*

\begin{proof}[Proof of \Cref{lem:proba_occupation_losange}]
As $\|\pi_T^t \circ \pi_T\|_{op} = \normop{\pi_T} \leq 1$, we have $\sqrt{|\det (\pi_T^t \circ \pi_T)|)} \leq 1$, so that the co-area formula \cite[Theorem~3.1]{Federer59} entails that
\begin{align*}
\int_{M \cap (\B(x,R) \setminus \B(x,h))} & \indicator{\pi_T(u-x) \in \B(O,r) \cap \B(\Omega,\rho)} f(u) \Haus^d(\dd u) 
\\
& \geq 
f_{\min} \int_{\pi_T(M \cap \B(x,R)-x)} \indicator{\pi_{T}^{-1}(v) \notin \B(0,h)} \indicator{\B(O,r) \cap \B(\Omega, \rho)}(v) \dd v.
\end{align*}
Since $\indicator{\pi_{T}^{-1}(v) \notin \B(0,h)} \geq \indicator{v \notin \B(0,h)}$, we get, provided $A$ is large enough,
\begin{align*}
& \int_{M \cap (\B(x,R) \setminus \B(x,h))}  \indicator{\pi_T(u-x) \in \B(O,r) \cap \B(\Omega,\rho)} f(u) \Haus^d(\dd u) 
\\
& \qquad \geq f_{\min} \int_{\pi_T(M \cap \B(x,R)-x)} \indicator{v \notin \B(0,h)} \indicator{\B(O,r) \cap \B(\Omega, \rho)}(v)\dd v 
\\
& \qquad \geq f_{\min} \int_{\B_T(0,r)} \indicator{\B(O,r) \cap \B(\Omega,\rho)}(v) \dd v - f_{\min} \omega_d h^d 
\\
&  \qquad \geq f_{\min} \left (\tilde{C}_d r^{\frac{d-1}{2}} A \frac{f_{\max}^4 \log n}{f_{\min}^5(n-1)} - C_d \omega_d \frac{f_{\max}^4 \log n}{f_{\min}^5(n-1)} \right ) 
\\
& \qquad  \geq 
A r^{\frac{d-1}{2}} \tilde{C}'_d  \frac{f_{\max}^4 \log n}{f_{\min}^4(n-1)},
\end{align*}
where the second to last inequality comes from \Cref{teclem:volume_intersection_balls}.
\end{proof}

%% file: proofs-main-detection.tex
\section{Proof of \Cref{lem:deterministic_protection}}
\label{sec:proof-of-lemdeterministicprotection}

\lemdeterministicprotection*

\begin{proof}[Proof of \Cref{lem:deterministic_protection}]

Suppose that  $X_i$ is detected in the tangent space $T_j$. Then $\norm{X_i-X_j}\leq r$, and there exists $\Omega\in T_j$ such that  
$\norm{\Omega-\pi_{T_j}(X_i-X_j)}\geq \rho \geq \rho_-$ and $\mathcal{Y}_{j} \cap \B(\Omega,\norm{\Omega-\pi_{T_j}(X_i-X_j)}) =\emptyset$. 
Since $\norm{\pi_{T_j}(X_i-X_j)}\leq r$, it follows that  
$
\norm{\Omega}\geq \rho_- -r\geq 2r >r + \eps_1
.
$
Hence, 
define $u:=\Omega/\norm{\Omega}$, and $x:=X_j + (r+\eps_1)u$.
As $\norm{\Omega- (x-X_j)} = \norm{\Omega} - r - \eps_1 \leq \norm{\Omega-\pi_{T_j}(X_i-X_j)} - \eps_1$, we get
\begin{equation}\label{deter:protec:boulvide}
(\B(x,\eps_1)-X_j)\cap \mathcal{Y}_{j}
\subset 
\B\bigl(\Omega,\norm{\Omega-\pi_{T_j}(X_i-X_j)}\bigr)
\cap 
\mathcal{Y}_{j}
=
\emptyset
.
\end{equation}
From \Cref{deter:protec:boulvide}, we now deduce that $x-X_j\notin \pi_{T_j}(\B(0,R_0)\cap (M-X_j))$. 
Indeed, if that was not the case,
there would exist $y\in M \cap \B(X_j,R_0)$ such that $\pi_{T_j}(y-X_j)=x-X_j$.
As $\dHaus(M,\mathcal{X}_n) \leq \eps_1$, there exists $X_k\in \B(y,\eps_1)\cap \mathcal{X}_n$. 
Since $\norm{y-X_j}\leq R_0\leq \tau_{\min}/32$ and $\theta \leq 1/24$,  \Cref{prop:projectionsestimtangentM} yields that 
$$
\norm{X_k-X_j}\leq \eps_1 + \norm{y-X_j}\leq \eps_1+\frac{\norm{x-X_j}}{1-\theta-\frac{\norm{y-X_j}}{2\tau_{\min}}}\leq \eps_1+\frac{192}{181}(r+\eps_1)
\leq 
R_0,
$$
and thus $X_k \in \B(X_j,R_0)$.
By definition of $\mathcal{Y}_j$, this leads to $\pi_{T_j}(X_k-X_j) \in \mathcal{Y}_{j}$, and since 
$$
\norm{\pi_{T_j}(X_k-X_j) - (x-X_j)} = \norm{\pi_{T_j}(X_k-y)} \leq \norm{X_k-y} \leq \varepsilon_1
,
$$
we get
$
\pi_{T_j}(X_k-X_j) \in \mathcal{Y}_{j}\cap (\B(x, \varepsilon_1)-X_j)
,
$
which contradicts \Cref{deter:protec:boulvide}.
As a result, $x-X_j\notin \pi_{T_j}(\B(0,R_0)\cap (M-X_j))$, so that \Cref{lem:GeoInter_bleupart1} asserts that
$$\frac{4}{5}\min\set{R_0-2 \eps_1,\dd(X_j,\partial M)} < \norm{x-X_j}=r+\eps_1
. 
$$
As $4(R_0-2 \varepsilon_1)/5 \geq  R_0/4 \geq r + \varepsilon_1$ by assumption, the above inequality yields that $\dd(X_j,\partial M) \leq 5(r+\eps_1)/4 \leq 2r < \infty$, and in particular that $\partial M \neq \emptyset$, hence the result.
 \end{proof}

%% file: tangent-space-estimation.tex
\section{Tangent space estimation}
	\label{sec:appendix:tangent-space-estimation}
	
\subsection{Tangent space of the manifold}
\label{sec:proof_tangent_space_estimation}

\proptangentspaceestimation*

\begin{proof}[Proof of \Cref{prop:tangent_space_estimation}] 
We let $h = \left ( \frac{\kappa}{f_{\min}} \frac{\log n}{n-1} \right )^{\frac{1}{d}}$, where $\kappa>1$ will be fixed later, and assume that $n$ is large enough so that $ h \leq \frac{\tau_{M}}{32}\wedge \frac{\tau_{\partial M}}{3} \wedge \frac{\tau_M}{\sqrt{d}}$. Without loss of generality we consider the case where $i=1$ and $X_1=0$. We let $x \in \B(0,h) \cap M$ be such that $\B(x,h/4) \cap \partial M = \emptyset$, according to \Cref{cor:sub_balls}. Slightly differing from the notation in \Cref{prop:projectionsestimtangentM}, for any vector $u \in \mathbb{R}^d$, we denote by $u_T = \pi_{T_x M}(u)$ and $u_\perp = \pi_{(T_x M)^\perp}(u)$. 
For short, we also write $p(h) := P(\B(0,h))$ and $p_n(h) := P_n(\B(0,h))$, where $P_n = n^{-1} \sum_{i = 1}^n \delta_{X_i}$ stands for the empirical measure. 
The proof of \Cref{prop:tangent_space_estimation} will make use of the following concentration result, borrowed from \cite{Aamari18}.
\begin{lemma}[{\cite[Lemma~9.5]{Aamari18}}]
\label{lem:concentration_TSE}
Write 
$$
\Sigma(h) 
:= 
\mathbb{E} \left ( X_T (X_T) ^t \indicator{\B(0,h)}(X) \right )
.
$$
Then for $n$ large enough, with probability larger than $1- 2 \left ( \frac{1}{n} \right )^{1 + \frac{2}{d}}$, we have,
\[
 p_n(h)  \leq  2 p(h) + \frac{10(2 + \frac{2}{d}) \log n}{n-1},
 \]
and
 \[
 \left \| \frac{1}{n-1} \sum_{i=2}^n  (X_i)_T (X_i)_T^t\indicator{\B(0,h)}(X_i) - \Sigma(h) \right \|_{F}  \leq  C_d \frac{f_{\max}}{f_{\min}\sqrt{\kappa}} p(h) h^2. 
\]
\end{lemma}
We now assume that the event described by \Cref{lem:concentration_TSE} occurs. We may decompose the local covariance matrix as
\begin{align*}
\frac{1}{n-1} \sum_{i=2}^n (X_i)(X_i)^t = \sum_{i=2}^n (X_i)_T (X_i)_T^t + R_1,
\end{align*}
where
$$R_1 := \frac{1}{n-1}\sum_{i=2}^n \left [ (X_i)_T (X_i)_\perp ^t  + (X_i)_{\perp}(X_i)_T^t + (X_i)_\perp (X_i)_\perp^t \right ].$$ 
Since $\B(0,h) \subset \B(x,2h)$, we have $\|(X_i)_T\| \leq h$ and, according to \cite[Theorem~4.18]{Federer59}, $\|(X_i)_\perp\| \leq \|(X_i-x)_\perp\| + \|(x-0)_\perp\| \leq \frac{3h^2}{\tau_M}$. Thus, $\|R_1\|_F \leq \frac{9h^3}{\tau_M}p_n(h) \leq C_d \frac{f_{\max} h^{d+3}}{\tau_M}$, according to \Cref{lem:concentration_TSE}.

Next, using \Cref{lem:concentration_TSE} again, we have 
\[
\lambda_{\min} \left ( \frac{1}{n-1} \sum_{i=2}^n  (X_i)_T (X_i)_T^t \mathbbm{1}_{\B(0,h)}(X_i) \right ) \geq \lambda_{\min} \left(\Sigma(h) \right ) -  C_d \frac{f_{\max}^2}{f_{\min}\sqrt{\kappa}} h^{d+2}.
\]
On the other hand, for $u \in T_x M$, we have
\begin{align*}
u^t \Sigma(h) u & = \int_{\B(0,h) \cap M} \left\langle u ,y_T \right\rangle^2 f(y) \Haus^d(\dd y) \\ 
                & \geq f_{\min} \int_{\B(x,h/4) \cap M} \left\langle u ,y_T \right\rangle^2 f(y) \Haus^d(\dd y) \\
                & \geq f_{\min} \int_{\B_d(0,h/4)} \left\langle u, \exp_x(v)_T - x_T   +x_T \right\rangle^2 \left | \det \left (d_v(\exp_x) \right ) \right | \dd v,
\end{align*}
according to \cite[Propositions~8.5 and~8.6]{Aamari18}. Moreover, \cite[Proposition~8.7]{Aamari18} ensures that we have $\left | \det \left (d_v(\exp_x) \right ) \right | \geq c_d$ provided that $\norm{v} \leq h/4 \leq \tau_M/{4}$, and \cite[Proposition~8.6]{Aamari18} gives $\exp_x(v) = x + v + R(v)$, with $\|R(v)\| \leq \frac{5\|v\|^2}{8\tau_M}$, under the same condition. Thus, 
\begin{align*}
u^t \Sigma(h) u & \geq c_d f_{\min} \int_{\B_d(0,h/4)} \left\langle u, v + R(v)_T + x_T \right\rangle^2 \dd v \\
                & \geq \frac{1}{2}c_d f_{\min} \int_{\B_d(0,h/4)} \left\langle u, v + x_T \right\rangle^2 \dd v \\
                & \qquad - 3c_d f_{\min}\int_{\B_d(0,h/4)} \left ( \frac{5 \|v\|^2}{8 \tau_{\min}} \right )^2 \dd v.
\end{align*}
Denoting by $\sigma_{d-1}$ the surface of the $(d-1)$-dimensional unit sphere and using polar coordinates yields
\begin{align*}
\int_{\B_d(0,h/4)} \left\langle u, v +  x_T \right\rangle ^2 \dd v & \geq \int_{\B_d(0,h/4)} \left\langle u, v\right\rangle ^2 \dd v  \geq \left (\frac{h}{4} \right )^{d+2} \frac{1}{d(d+2)}\sigma_{d-1},
\end{align*}
and
\begin{align*}
\int_{\B_d(0,h/4)} \left ( \frac{5 \|v\|^2}{8 \tau_{\min}} \right )^2 \dd v \leq \left ( \frac{5}{8} \right )^2  \frac{\sigma_{d-1}}{(d+4)\tau_M^2} \left ( \frac{h}{4} \right )^{d+4}.
\end{align*}
Since $h \leq \tau_M/\sqrt{d}$, it follows that 
\[
\lambda_{\min}(\Sigma(h)) \geq c_d f_{\min} h^{d+2},
\]
for some positive constant $c_d$. Gathering all pieces and using \cite[Theorem~10.1]{Aamari18} leads to 
\[
\angle (T_x M,\hat{T}_i) \leq C_d  \frac{f_{\max} h}{\tau_M(c_d f_{\min} - C_d (f_{\max}^2/(f_{\min}\sqrt{\kappa})))}.
\]
Thus, choosing $\kappa = C_d \left ( \frac{f_{\max}}{f_{\min}} \right )^4$, for $C_d$ large enough, gives 
\[
\angle(T_x M,\hat{T_i}) \leq C_d \frac{f_{\max}h}{f_{\min}\tau_M}.
\]
Noting that $ \angle(T_0 M, T_x M) \leq 2h/\tau_M$ from \Cref{prop:tangent_variation_geodesic,lem:dist_geod_dist_eucl}, the result of \Cref{prop:tangent_space_estimation} follows after using a union bound.
\end{proof}

\subsection{Tangent space of the boundary}\label{sec:proof_cor_boundary_tangent_space_estimation}

\corboundarytangentspaceestimation*

\begin{proof}[Proof of \Cref{cor:boundary_tangent_space_estimation}]

Under the assumptions of \Cref{thm:detection} and \Cref{prop:tangent_space_estimation}, we let $X_i \in \mathcal{Y}_{R_0,r,\rho}$,  $\varepsilon_{\partial M} = \left ( C_d R_0 \frac{f_{\max}}{f_{\min}^2} \frac{\log n}{n} \right )^{\frac{1}{d+1}}$, and $h = \left ( C_d \frac{f_{\max}^4}{f_{\min}^5} \frac{\log n}{n-1} \right )^{\frac{1}{d}}$ so that with probability larger than $1- 4n^{-2/d}$, we have
\[
\angle(\eta_{\pi_{\partial M}(X_i)},\tilde{\eta}_i) \leq \frac{\varepsilon_{\partial M}}{R_0}
\text{~~~and~~~}
\angle(T_{X_i}M,\hat{T}_i) \leq C_d \frac{f_{\max}}{f_{\min}} \frac{h}{\tau_{\min}}.
\]
Combining \Cref{thm:detection} \ref{item:thm_detection_emptyboundary} with \Cref{lem:dist_geod_dist_eucl} and \Cref{prop:tangent_variation_geodesic} entails
\begin{align*}
\angle(T_{\pi_{\partial M}(X_i)}M,\hat{T}_i)
&\leq  
\angle(T_{\pi_{\partial M}(X_i)}M,T_{X_i}M) 
+
\angle(T_{X_i}M,\hat{T}_i)
\\
&\leq 
2 \frac{\varepsilon_{\partial M}^2}{R_0} 
+
C_d \frac{f_{\max}}{f_{\min}} \frac{h}{\tau_{\min}}
\leq 
C_d \varepsilon_{\partial M}
,
\end{align*} 
for $n$ large enough. Finally, since
\[
\angle(T_{\pi_{\partial M}(X_i)} \partial M,\hat{T}_{\partial,i})
\leq 
\angle(T_{\pi_{\partial M}(X_i)}M,\hat{T}_i)
+ 
\angle(\eta_{\pi_{\partial M}(X_i)},\tilde{\eta}_i)
,
\]
the bound follows.
\end{proof}

%% file: linear-patches.tex
\section{Local linear patches}
	\label{sec:appendix:linear-patches}

\thmlocallinearpatchdeterministic*

\begin{proof}[Proof of \Cref{thm:local_linear_patch_deterministic}]

First, note that the choice $r_0 = (\tau_{\min} \wedge \tau_{\partial, \min})/40$ satisfies the requirements of \Cref{lem:GeoInter_bleupart2}, 
for  a radius $R_0=\tau_{\min}/16$.
For short, let $\mathbb{M} := \mathbb{M}(\mathcal{X}_n,\mathcal{X}_\partial,T,\eta)$. 
\begin{itemize}[leftmargin=*]
\item
Let $x\in M$ be fixed. 
We bound $\dd(x,\mathbb{M})$ depending on its closeness to $\partial M$.
\begin{itemize}[leftmargin=*]
\item
First assume that $\dd(x,\partial M) \leq \varepsilon_{\partial M}-\delta$.
Then $\dd(x,\mathcal{X}_{\partial})\leq \varepsilon_{\partial M}$, and we let  $X_{i_0}\in \mathcal{X}_{\partial}$ be such that $\norm{x-X_{i_0}}\leq \varepsilon_{\partial M}$.  Without loss of generality we may assume that $i_0=1$.
Let $\mathbb{P}_1:=X_1 + \B_{T_1}(0,\varepsilon_{\partial M})\cap\{z,\langle z-X_1,{\eta}_1\rangle \leq 0 \} \subset \mathbb{M}_\partial$ denote the half-patch at $X_1$.

From \Cref{prop:projectionsestimtangentM}, we have
\begin{equation}\label{cas11}
 \norm{\pi_{X_1+T_1}(x)-x}\leq \varepsilon_{\partial M}\left(\theta+\frac{\varepsilon_{\partial M}}{2\tau_{\min}} \right).
\end{equation}
As a result, if $\pi_{X_1+T_1}(x)\in \mathbb{P}_1$, then $\dd(x,\mathbb{M})\leq \norm{\pi_{X_1+T_1}(x)-x}$ yields the desired bound.
Otherwise, if $\pi_{X_1+T_1}(x)\notin \mathbb{P}_1$, since $\dd(x,\mathbb{P}_1) \leq \norm{x-X_1} \leq \eps_{\partial M}$, we can decompose $\pi_{X_1+T_1}(x)$ as $\pi_{X_1+T_1}(x)=X_1+\alpha {\eta}_1 + \beta v$,  with unit $v \in T_1 \cap \Span(\eta_1)^\perp$, and $\alpha = \dd(\pi_{X_1+T_1}(x),\mathbb{P}_1) > 0$ such that $\alpha^2+\beta^2 \leq \eps_{\partial M}^2$.
Writing $x_1:=\pi_{\partial M}(X_1)$, triangle inequality ensures  that 
$$
\norm{x-x_1}
\leq
\norm{x-X_1} + \norm{X_1-x_1}
\leq
\eps_{\partial M} + a\delta^2
\leq
\tau_{\min}/32
.
$$
From \Cref{lem:dist_geod_dist_eucl,prop:tangent_variation_geodesic}, we also have 
$$
\angle(T_{x_1}M,T_1)
\leq 
\angle(T_{x_1}M,T_{X_1} M)
+
\angle(T_{X_1}M,T_1)
\leq
2\eps_{\partial M}/\tau_{\min}
+
\theta 
\leq 
1/8
.
$$
As a result, \Cref{lem:GeoInter_bleupart2} applies and gives 
$$\pi_{T_1}(x-X_1)\in \B_{T_1}(0,\varepsilon_{\partial M})\cap (\B_{T_1}(\pi_{T_1}(x_1-X_1)+r_0 \eta_{x_1},r_0))^c.$$
Thus, we have 
	\begin{align*}
	r_0
	&\leq
	\norm{\alpha {\eta}_1+\beta v - r_0 \eta_{x_1} + \pi_{T_1}(X_1-x_1)}
	\\
	&\leq
	\norm{\alpha {\eta}_1 + \beta v - r_0 \eta_{x_1} }+a\delta^2
	,
	\end{align*}
which, since $a\delta^2 \leq r_0$ and $\alpha^2 + \beta^2 \leq \eps_{\partial M}$, leads to
\begin{align*}
(r_0-a\delta^2)^2
&\leq 
\norm{(\alpha {\eta}_1+\beta v) - r_0 \eta_{x_1}}^2
\\
&\leq
\varepsilon_{\partial M}^2+r_0^2-2r_0\alpha \langle {\eta}_1, \eta_{x_1}\rangle -2r_0 \beta \langle v,\eta_{x_1}\rangle
.
\end{align*}
  As $\langle {\eta}_1, \eta_{x_1}\rangle = 1-\norm{\eta_1-\eta_{x_1}}^2/2 \geq 1- \theta'^2/2>0$ 
  and
  $|\langle v,\eta_{x_1}\rangle| = |\langle v,\eta_1-\eta_{x_1}\rangle| \leq \theta'
  $, 
  we deduce that
 \begin{equation*}
  \alpha 
  = 
  \dd(\pi_{X_1+T_1}(x),\mathbb{P}_1) 
  \leq 
  \frac{\varepsilon_{\partial M}^2+2r_0a\delta^2+ r_0\varepsilon_{\partial M} \theta'}{2r_0(1-\theta'^2/2)}
  \leq 
  \frac{\varepsilon_{\partial M}^2}{r_0}
  +
  2a\delta^2
  +
  \varepsilon_{\partial M}\theta'
  .
 \end{equation*}
At the end of the day, combining the above inequality with \eqref{cas11} yields the bound
 \begin{equation}\label{eq:patchesineg1}  \dd(x,\mathbb{M}) \leq 2a\delta^2 + \varepsilon_{\partial M}\left(\theta+\theta'+ \frac{2 \varepsilon_{\partial M}}{r_0}
 \right)
 ,
 \end{equation}
 which also holds if $\pi_{X_1+T_1}(x) \in \mathbb{P}_1$.

\item 
Now, assume that $\dd(x,\partial M)> \varepsilon_{\partial M}-\delta$. Let $X_{i_0}$ denote the closest point to 
$x$ in $\mathcal{X}_n$, with $i_0 = 1$ without loss of generality. Since 
$\norm{x-X_1}\leq\eps_0$, we deduce that 
\begin{align*}
\dd(X_1,\mathcal{X}_{\partial})
&\geq 
\dd(x,\mathcal{X}_{\partial})
-
\norm{x-X_1}
\\
&\geq
\dd(x,\partial M)
-
a\delta^2
-
\eps_0
\\
&\geq
\varepsilon_{\partial M}-\delta-\eps_0-a \delta^2
\\
&\geq 
\varepsilon_{\partial M}/2
.
\end{align*}
Thus $X_1 \in \interior{\mathcal{X}}_{\varepsilon_{\partial M}}$, and therefore
$\mathbb{P}_1:=X_1+\B_{T_1}(0,\varepsilon_{\mathring{M}})$ is a patch of 
$\mathbb{M}_{\Int} \subset \mathbb{M}$. 
Because $\eps_{\mathring{M}} \geq \eps_0$, the point $\pi_{X_1+T_1}(x)$ belongs to $\mathbb{P}_1$,
so that $\dd(x,\mathbb{M})\leq \norm{\pi_{X_1+T_1}(x)-x}$. Using \Cref{prop:projectionsestimtangentM} again, we get
 \begin{equation}\label{eq:patchesineg2}
   \dd(x,\mathbb{M}) \leq  \eps_0\left(\theta + \frac{\eps_0}{2\tau_{\min}} \right). 
 \end{equation}
 
\end{itemize}

\item
 Let now $x\in \mathbb{M}$ be fixed. We bound $\dd(x,M)$ depending on whether $x$ belongs to a ``boundary patch'' (i.e. to $\mathbb{M}_\partial$) or an ``interior patch'' (i.e. to $\mathbb{M}_{\Int}$).

\begin{itemize}[leftmargin=*]
\item
Assume that $x\in \mathbb{M}_\partial$ belongs to ``boundary patch''. That is, without loss of generality, $x\in X_1 + \B_{T_1}(0,\varepsilon_{\partial M})\cap\{z,\langle z-X_1,{\eta}_1 \rangle \leq 0\}$ with $X_1 \in \mathcal{X}_{\partial}$.
Define 
$x_1 := \pi_{\partial M}(X_1)$, 
$$
x_1^*:=\pi_{\pi_{X_1+T_1}(\partial M\cap \B(X_1,\tau_{\min}/16))}(X_1)
,
$$
and let
$x'_1\in \partial M\cap \B(X_1,\tau_{\min}/16)$ be such that $\pi_{X_1 + {T}_1}(x_1')=x_1^*$. 
According to  \Cref{cor:projorthobord}, we have  $x_1^*-X_1= \norm{x_1^*-X_1}\eta_1^*$, where
$\eta_1^*$ is the unit vector of $Nor(x_1',M)\cap {T}_1$.
Furthermore, \Cref{prop:tangent_variation_geodesic}, \Cref{prop:geodesic_vs_euclidean} and \Cref{prop:estimeta} combined
yield the bound 
\begin{align*}
\Vert \eta_1^* - \eta_{x'_1} \Vert 
\leq
\sqrt{2}
\angle(T_{x_1'} M , T_1)
\leq \sqrt{2}(\theta + {2\|x'_1-X_1\|}/{\tau_{\min}})
.
\end{align*}
Furthermore, by definition of $x_1^*$ and the fact that $x_1 \in \partial M \cap \B(X_1,\tau_{\min}/16)$, we also have
$$\norm{X_1-x_1^*} \leq  \norm{\pi_{{T}_1}(X_1-x_1)}
\leq 
\norm{X_1-x_1}
=
\dd(X_1,\partial M)
\leq
a\delta^2
.
$$
As $\norm{X_1-x_1'} \leq \tau_{\min}/16$, \Cref{prop:projectionsestimtangentM} ensures that
$\norm{X_1-x_1^*}\geq \norm{X_1-x_1'}(1-\theta-1/32)$, which leads to $\norm{X_1-x_1'}\leq 2 a \delta^2$ and hence to $\norm{x_1-x_1'}\leq 3a\delta^2 \leq (\tau_{\min} \wedge \tau_{\partial, \min})/32$. 
As a result, \Cref{prop:normal_vector_stability} applies and asserts that
\begin{align*}
\Vert \eta_{x_1} - \eta_{x_1'} \Vert  
&\leq 
\frac{9 \Vert x_1 - x'_1 \Vert}{\tau_{\min} \wedge \tau_{\partial, \min}}
\leq
\frac{27 a\delta^2}{\tau_{\min} \wedge \tau_{\partial, \min}}
.
\end{align*}
Gathering all the pieces together, we obtain
\begin{align*}
\Vert \eta^*_1 - {\eta}_1 \Vert  
& \leq 
 \Vert \eta_1^* - \eta_{x'_1} \Vert  + \Vert \eta_{x'_1} - \eta_{x_1} \Vert  + \Vert \eta_{x_1} - {\eta}_1 \Vert  
\\
& \leq 
\sqrt{2} \theta + \theta' + \frac{(27+4\sqrt{2}) a\delta^2}{\tau_{\min} \wedge \tau_{\partial, \min}} 
\\
&\leq
\sqrt{2}\theta + \theta' + \frac{a\delta^2}{r_0}
\\
&:= \theta''.
\end{align*}
Now, if $x\in \B_{X_1+T_1}(x_1^*-r_0 \eta_{1}^*,r_0)$, we have 
 $
 \dd(x,\B_{X_1+T_1}(x_1^*-r_0 \eta_{1}^*,r_0))=0
 .
 $
 Otherwise, if $x\notin \B_{X_1+T_1}(x_1^*-r_0 \eta_{1}^*,r_0)$,
we have 
 $$
 \dd(x,\B_{X_1+T_1}(x_1^*-r_0 \eta_{1}^*,r_0))= \norm{x-(x^*_1-r_0\eta_1^*)}-r_0>0
.
 $$
We may hence write
 \begin{equation*}
 \begin{cases}
x-X_1=-\alpha {\eta}_1+ \beta v \text{ with }\alpha\geq 0, \alpha^2+\beta^2\leq\varepsilon_{\partial M}^2, \text{ and unit }v\in T_1\cap\Span(\eta_1)^{\perp},\\
x_1^*-X_1= t \eta_1^* \text{ with } 0\leq t \leq a \delta^2 \text{ and } \Vert {\eta}_1 - \eta_1^* \Vert \leq \theta''.\\
 \end{cases}
\end{equation*}
Since $\langle{\eta_1},{\eta_1^*}\rangle \geq 0$ and $|\langle{v},{\eta_1^*}\rangle|=|\langle{v},{\eta_1^*-\eta_1}\rangle| \leq \theta''$, it follows that 
\begin{align*}
\norm{x-(x^*_1-r_0\eta_1^*) }^2
&=
\norm{(x-X_1)+(r_0-t)\eta_1^*}^2
\\
&\leq
\varepsilon_{\partial M}^2
+
2(r_0-t)
(
\inner{-\alpha \eta_1}{\eta_1^*}
+
\inner{\beta v}{\eta_1^*}
)
+
(r_0-t)^2
\\
&\leq  
\varepsilon_{\partial M}^2
+
2\varepsilon_{\partial M} (r_0-t)\theta''
+
(r_0-t)^2
.
\end{align*}
Therefore, no matter whether or not $x$ belongs to $\B_{X_1+T_1}(x_1^*-r_0 \eta_{1}^*,r_0)$, we have
 \begin{align*}
\dd(x, \B_{X_1+T_1}(x_1^*-r_0 \eta_{1}^*,r_0))
&\leq 
\varepsilon_{\partial M} \theta''+\frac{\varepsilon_{\partial M}^2}{2(r_0-a\delta^2)}
\\
&\leq 
\varepsilon_{\partial M} \theta''+\frac{\varepsilon_{\partial M}^2}{r_0}
.
 \end{align*}
From the left-hand side inclusion of \Cref{lem:GeoInter_bleupart2}, we hence get the existence of some $y\in \B(x_1',\tau_{\min}/16)\cap M$ such that
  $$
  \norm{x-\pi_{X_1+T_1}(y)}
  \leq  
\varepsilon_{\partial M} \theta''+\frac{\varepsilon_{\partial M}^2}{r_0}
  .
  $$
We will now show that this point $y \in M$ is close to $x$.

For this, a first (rough) bound on $\norm{y-X_1}$ may be derived, using $\norm{y-X_1}\leq \norm{y-x_1'}+\norm{x_1'-X_1}\leq \tau_{\min}/16+2a\delta^2\leq \tau_{\min}/8$.
According to \Cref{prop:projectionsestimtangentM}, we have  
$$
\norm{y-X_1}
\leq 
\frac{\norm{\pi_{T_1}(y-X_1)}}{1-\theta-\norm{y-X_1}/(2\tau_{\min})}
\leq
2\norm{\pi_{T_1}(y-X_1)},
$$
which, by using the other bound of \Cref{prop:projectionsestimtangentM}, leads to 
\begin{align*}
\norm{y-\pi_{X_1+T_1}(y)}
&\leq 
2\norm{\pi_{T_1}(y-X_1)}
\left(\theta+\frac{\norm{\pi_{T_1}(y-X_1)}}{\tau_{\min}} \right)
,
\end{align*}
Hence, further bounding
\begin{align*}
\norm{\pi_{T_1}(y-X_1)}& \leq \norm{x-X_1} + \norm{x - \pi_{X_1+T_1}(y)}\\ 
&\leq 
\varepsilon_{\partial M}
+ 
\varepsilon_{\partial M} \theta''+\frac{\varepsilon_{\partial M}^2}{r_0}
\\
&\leq
2\eps_{\partial M}
\end{align*} 
since $\theta''\leq 1/2$ and $\eps_{\partial M} \leq r_0/2$, we finally obtain
\begin{align*}
\norm{x-y}
&\leq 
\norm{x - \pi_{X_1+T_1}(y)} + \norm{y - \pi_{X_1+T_1}(y)} 
\\
& \leq  
\varepsilon_{\partial M} \theta''
+\frac{\varepsilon_{\partial M}^2}{r_0}
+
4\eps_{\partial M}
\left(\theta+\frac{2\eps_{\partial M}}{\tau_{\min}} \right)
\\
&\leq
8\eps_{\partial M}
\left(
\theta + \theta' + \frac{\eps_{\partial M}}{r_0}
\right)
,
\end{align*} 
where we used that $a\delta^2 \leq \eps_{\partial M}$.
In particular, we have
 \begin{equation}\label{patchesineg3}
 \dd(x,M) \leq
 8\eps_{\partial M}
\left(
\theta + \theta' + \frac{\eps_{\partial M}}{r_0}
\right)
.
 \end{equation}
\item
Assume that $x \in \mathbb{M}_{\Int}$ belongs to an ``interior patch''. That is, without loss of generality, $x\in X_1+ \B_{T_1}(0,\varepsilon_{\mathring{M}})$ with 
 $\dd(X_1,\mathcal{X}_{\partial})\geq \varepsilon_{\partial M}/2$. We have 
 $\dd(X_1,\partial M)\geq \varepsilon_{\partial M}/2-\delta \geq {3\varepsilon_{\mathring{M}}}/{2}$, so that an applying
 \Cref{lem:GeoInter_bleupart1} at $X_1$ provides the existence of some $y\in M\cap \B(X_1,\eps_{\mathring{M}})$ such that
$x=\pi_{X_1+T_1}(y)$.  Thus,  \Cref{prop:projectionsestimtangentM} entails
  \begin{align}
   \dd(x,M) 
   &\leq 
   \norm{y-x}
   =
   \Vert(y-X_1)^\perp\Vert 
   \leq  
   {\varepsilon_{\mathring{M}}} \left( \theta+ \frac{\varepsilon_{\mathring{M}}}{2\tau_{\min}}\right). \label{patchesineg4}
 \end{align}

\end{itemize}

\end{itemize} 
 
To conclude the proof of \Cref{thm:local_linear_patch_deterministic}, we combine the above results as follows.
\begin{enumerate}[label=(\roman*),leftmargin=*]
 \item If $\partial M = \emptyset$, then		 $\dd(x,\partial M) = \infty$ for all $x\in \R^D$, so that $\mathcal{X}_\partial = \emptyset$ and hence $\mathbb{M}_\partial = \emptyset$. As a result, $\dHaus(M,\mathbb{M})$ is bounded by the maximum of \Cref{eq:patchesineg2,patchesineg4}. The requirement
 $\eps_0\leq \varepsilon_{\mathring{M}}$ ensures that
 \[
 \dHaus(M,\mathbb{M})\leq  {\varepsilon_{\mathring{M}}} \left( \theta+ \frac{\varepsilon_{\mathring{M}}}{2\tau_{\min}}\right).
 \]
 \item If $\partial M\neq \emptyset$, then $\dHaus(M,\mathbb{M})$ is bounded by the maximum of \Cref{eq:patchesineg1,eq:patchesineg2,patchesineg3,patchesineg4}. This boils down to
  $$\dHaus(M,\mathbb{M}) 
  \leq  
  2a\delta^2
  +
   8\eps_{\partial M}
\left(
\theta + \theta' + \frac{\eps_{\partial M}}{r_0}
\right)
.
 $$
\end{enumerate}

\end{proof}

%% file: lower-bounds-proofs.tex
\section{Proofs of the minimax lower bounds}
\label{sec:proof-lower-bounds-main}
The minimax lower bounds (\Cref{thm:boundary-main-lower-bound,thm:manifold-main-lower-bound}) will be proven using the standard Bayesian arguments relying on hypotheses comparison method. This is usually referred to as Le Cam's method. It involves the total variation distance, for which we recall a definition.

\begin{definition}[Total Variation]
	\label{def:total_variation}
	For any two Borel probability distributions $P_0,P_1$ over $\R^D$, the total variation between them is defined as
	\begin{align*}
	\TV(P_0,P_1)
	:=
	\frac{1}{2}
	\int_{\R^D}
	|f_1 - f_0| d\mu
	,
	\end{align*}
where $\mu$ is a $\sigma$-finite measure dominating $P_0$ and $P_1$, with respective densities $f_0$ and $f_1$.
\end{definition}

In the context of manifold and boundary estimation for the Hausdorff distance $\dHaus$, Le Cam's lemma~\cite{Yu97} writes as follows.
\begin{lemma}
	\label{lem:le_cam_applied}
	Fix an integer $n\geq 1$ and write $\mathcal{P} = \mathcal{P}^{d,D}_{\tau_{\min},\tau_{\partial,\min}}(f_{\min},f_{\max})$. 
	\begin{enumerate}[label=(\roman*)]
	\item 
	\label{item:le_cam_applied-manifold}
	Then for all $P_0,P_1 \in \mathcal{P}$ with respective supports $M_0$ and $M_1$,
	\begin{align*}
	\inf_{\hat{M}}
	\sup_{P \in \mathcal{P}}
	\E_{P^n}
	\left[
	\dHaus\bigl(M,\hat{M}\bigr)
	\right]
	&\geq
	\frac{1}{2}
	\dHaus(M_0,M_1)
	\left(1-\TV(P_0,P_1)\right)^n,
	\end{align*}
	where the infimum ranges among all the estimators $\hat{M} = \hat{M}(X_1,\ldots,X_n)$.
	\item
		\label{item:le_cam_applied-boundary}
	If in addition, $\partial M_0$ and $\partial M_1$ are non-empty,
	\begin{align*}
	\inf_{\hat{B}}
	\sup_{P \in \mathcal{P}}
	\E_{P^n}
	\left[
	\dHaus\bigl(\partial M,\hat{B}\bigr)
	\mathbbm{1}_{\partial M \neq \emptyset}
	\right]
	&\geq
	\frac{1}{2}
	\dHaus(\partial M_0, \partial M_1)
	\left(1-\TV(P_0,P_1)\right)^n,
	\end{align*}
	where the infimum ranges among all the estimators $\hat{B} = \hat{B}(X_1,\ldots,X_n)$.
	\end{enumerate}
	\end{lemma}

\begin{proof}[Proof of \Cref{lem:le_cam_applied}]
Apply \cite[Lemma 1]{Yu97} with loss function $\dHaus$, model $\mathcal{P}$, parameters of interest $\theta(P) = \supp(P)$ and $\theta(P) = \partial \bigl(\supp(P)\bigr)$ respectively, and conclude with the bound $(1-\TV(P_0^n,P_1^n)) \geq (1-\TV(P_0,P_1))^n$.
\end{proof}

Aiming at applying \Cref{lem:le_cam_applied}, we shall first describe how to construct hypotheses $P_0$ and $P_1$ that belong to the models, close in total variation distance but with supports (or boundary) far away in Hausdorff distance.

\subsection{Hypotheses with empty boundary}

To do so in the boundariless case $\tau_{\partial,\min} = \infty$, we will use a structural stability result of the family of model.
We recall that $\normop{\cdot}$ denotes the operator norm, that is $\normop{A} = \max_{\norm{v} = 1} \norm{Av}$ for all $A \in \R^{D\times D}$.

\begin{restatable}[Reach Stability]{proposition}{propdiffeomorphismstability}
\label{prop:diffeomorphism_stability}
Let $M \in  \mathcal{M}^{d,D}_{\tau_{\min},\tau_{\partial,\min}}$ and $\Phi : \R^D \rightarrow \R^D$ be a $\mathcal{C}^2$ map such that $\lim_{\norm{x} \to \infty} \norm{\Phi(x)} = \infty$.
Assume that 
$
\sup_{x \in \R^D} \normop{I_D - d_x \Phi} \leq 1/10
$
.
Then $\Phi$ is a global diffeomorphism, and the image $\Phi(M)$ of $M$ by $\Phi$ satisfies:
\begin{itemize}
\item
$\partial \Phi(M) = \Phi(\partial M)$,
\item
If $\sup_{x \in \R^D} \normop{d^2_x \Phi} \leq 1/\left(2 \tau_{\min} \right)$, then $\tau_{\Phi(M)} \geq \tau_{\min}/2$,
\item
If $\sup_{x \in \R^D} \normop{d^2_x \Phi} \leq 1/\left(2 \tau_{\partial \min} \right)$, then $\tau_{\partial \Phi(M)} \geq \tau_{\partial, \min}/2$.
\end{itemize}
\end{restatable}

\begin{figure}[!htbp]
\centering
\includegraphics[width = 0.7\linewidth,page = 1]{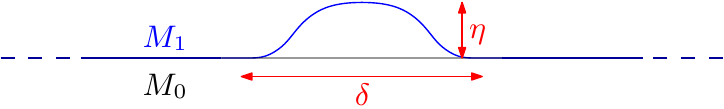}
\caption{Boundariless supports $M_0$ and $M_1$ of \Cref{prop:hypotheses_empty_boundary} for $d=1$ and $D=2$.
Here, the total variation between the associated uniform distributions is of order $\TV(P_0,P_1) \asymp f_{\min} \Haus^d(M_0 \triangle M_1) \asymp f_{\min} \delta^d$ and Hausdorff distance $\dHaus(M_0,M_1) = \eta$.
The reach bound forces the bump to have height $\eta \lesssim \delta^2/\tau_{\min}$, so that optimal parameter choices yield:\vspace{0.5em}
\\
\vspace{0.5em}
\centerline{
$\displaystyle \delta \asymp \left(\frac{1}{f_{\min} n}\right)^{1/d}$
\text{ and }
$\displaystyle  \eta \asymp 
\frac{\delta^2}{\tau_{\min}}
\asymp
\tau_{\min}
\left(\frac{1}{f_{\min} \tau_{\min}^d n}\right)^{2/d}
$.
}
As $\TV(P_0,P_1) \leq 1$, this can only be done when $f_{\min} \delta^d \lesssim 1$, i.e. $n \gtrsim 1/(f_{\min}\tau_{\min}^d)$.
}
\label{fig:hypotheses_le_cam_boundariless}
\end{figure}

The proof is to be found in \Cref{sec:model_stability}.
Essentially, the class $\bigl\{\mathcal{M}^{d,D}_{\tau_{\min},\tau_{\partial,\min}}\bigr\}_{\tau_{\min},\tau_{\partial,\min}}$ is stable up to $\mathcal{C}^2$-diffeomorphism, with explicit bounds on the parameters.
From there, we consider $P_0$ over a boundariless manifold $M_0 \in \mathcal{M}^{d,D}_{2 \tau_{\min},\infty}$, and $P_1$ over $M_1$ that is obtained by bumping $M_1$ locally (see \Cref{fig:hypotheses_le_cam_boundariless}). The method is similar to that of \cite[Lemma~5]{Aamari19b}, with an explicit dependency in the parameters of the model.

\begin{restatable}[Hypotheses with Empty Boundary]{proposition}{prophypothesesemptyboundary}
	\label{prop:hypotheses_empty_boundary}
	Assume that $f_{\min} \leq c_d/\tau_{\min}^d$ and $c'_d/\tau_{\min}^d \leq f_{\max}$, for some small enough $c_d,(c'_d)^{-1}>0$. 
	
	If $d \leq D-1$, then for all $n \geq C_d/(f_{\min}\tau_{\min}^d)$, there exist $P_0,P_1 \in \mathcal{P}^{d,D}_{\tau_{\min},\infty}(f_{\min},f_{\max})$ with boundariless supports $M_0$ and $M_1$ such that
	\begin{align*}
	\TV(P_0,P_1) \leq \frac{1}{n}
	\text{~~~~and~~~~}
	\dHaus(M_0,M_1) 
	\geq 
	C'_d \tau_{\min}
	\left(\frac{1}{f_{\min} \tau_{\min}^d n}\right)^{2/d}
	.
	\end{align*}
\end{restatable}

	See \Cref{sec:construction-of-hypotheses} for the construction of these hypotheses. 
	We are now in position to prove \Cref{thm:manifold-main-lower-bound}~\ref{item:minimax_lower_bound_boundaryless}.

    \begin{proof}[Proof of \Cref{thm:manifold-main-lower-bound}~\ref{item:minimax_lower_bound_boundaryless}]
    Let $\mathcal{P}$ denote the model $\mathcal{P}^{d,D}_{\tau_{\min},\infty}(f_{\min},f_{\max})$, and write $n_0 :=  \ceil{C_d/(f_{\min} \tau_{\min}^d)}$, where $C_d>0$ is the constant of \Cref{prop:hypotheses_empty_boundary}.
    \begin{itemize}
	\item    
    If $n \geq n_0$, applying \Cref{lem:le_cam_applied} \ref{item:le_cam_applied-manifold} with hypotheses $P_0$ and $P_1$ of \Cref{prop:hypotheses_empty_boundary}, yields
	    	    \begin{align*}
	    	\inf_{\hat{M}}
			\sup_{P \in \mathcal{P}}
			\E_{P^n}
			\left[
			\dHaus\bigl(M,\hat{M}\bigr)
			\right]
			&\geq
			\frac{1}{2} C'_d \tau_{\min}
			\left(\frac{1}{f_{\min} \tau_{\min}^d n}\right)^{2/d}	
			\left(1 - \frac{1}{n} \right)^{n}
			\\
			&\geq
			C''_d \tau_{\min}
			\left\{
				1
				\wedge
				\left(
					\frac{1}{f_{\min} \tau_{\min}^d n}
				\right)^{2/d}
			\right\}
			.
	    \end{align*}
	    \item
	    Otherwise, if $n < n_0$, note that since
	    	$
	    	\inf_{\hat{M}}
			\sup_{P \in \mathcal{P}}
			\E_{P^n}
			\left[
			\dHaus\bigl(M,\hat{M}\bigr)
			\right]
			$
	    is a non-increasing sequence, the previous point yields
	    \begin{align*}
	    	\inf_{\hat{M}}
			\sup_{P \in \mathcal{P}}
			\E_{P^n}
			\left[
			\dHaus\bigl(M,\hat{M}\bigr)
			\right]
			&\geq
	    	\inf_{\hat{M}}
			\sup_{P \in \mathcal{P}}
			\E_{P^{n_0}}
			\left[
			\dHaus\bigl(M,\hat{M}\bigr)
			\right]
			\\
			&\geq
			C''_d \tau_{\min}
			\left(\frac{1}{f_{\min} \tau_{\min}^d n_0}\right)^{2/d}
			\\
			&\geq
			\tilde{C}'_d
			\tau_{\min}
			\geq
			\tilde{C}'_d
			\tau_{\min}
			\left\{
				1
				\wedge
				\left(
					\frac{1}{f_{\min} \tau_{\min}^d n}
				\right)^{2/d}
			\right\}
			,
	    \end{align*}
	    which concludes the proof.
	    \qedhere
    \end{itemize}
    \end{proof}

\subsection{Convex hypotheses (with boundary)}
Similarly to the previous section, we shall use a stability result under diffeomorphisms in the convex case $\tau_{\min} = \infty$. 
Unfortunately, \Cref{prop:diffeomorphism_stability} only provides convexity of $\Phi(M)$ (i.e. $\tau_{\Phi (M)} = \infty$) for diffeomorphisms $\Phi$ that are affine maps, which does not allow enough flexibility.
Beyond affine maps, the following result allows to quantify how much one may bump a strictly convex full dimensional domain while keeping it convex.

\begin{restatable}[Stability of Strict Convexity]{proposition}{propdiffeomorphismconvexitystability}
\label{prop:diffeomorphism_convexity_stability}
	Let $C \subset \R^d$ be a compact domain with $\interior{C} \neq \emptyset$, that has a $\mathcal{C}^2$ boundary $\bar{\partial} C$. 
	Assume that:
	\begin{itemize}
	\item
	for all $x \in \bar{\partial} C$,
	$\bar{\partial} C \setminus \{x\}$ is connected;
	\item
	 for all $x, y \in \bar{\partial} C$,
	$
	\dd(y-x,T_x \bar{\partial} C) 
	\geq
	A \norm{y-x}^2
	,
	$
	for some $A > 0$.
	\end{itemize}
	Let $\Phi : \R^d \rightarrow \R^d$ be a $\mathcal{C}^2$ map such that $\lim_{\norm{x} \to \infty} \norm{\Phi(x)} = \infty$, $\normop{I_d - d \Phi} \leq 1/10$ and $\normop{d^2 \Phi} \leq A$, then $C$ and $\Phi(C)$ are convex.
\end{restatable}

\begin{figure}[!htbp]
\centering
\includegraphics[width = 0.7\linewidth,page = 2]{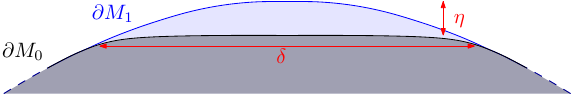}
\caption{
Convex supports $M_0$ and $M_1$ of \Cref{prop:hypotheses_with_boundary} for $d=D=2$.
Here, the total variation between the associated uniform distributions is of order $\TV(P_0,P_1) \asymp f_{\min} \Haus^d(M_0 \triangle M_1) \asymp f_{\min} \delta^{d-1}\eta$ and Hausdorff distance $\dHaus(M_0,M_1) = \dHaus(\partial M_0,\partial M_1) = \eta$.
The reach bound forces the bump to have height $\eta \lesssim \delta^2/\tau_{\partial, \min}$, so that optimal parameter choices yield:\vspace{0.5em}
\\
\vspace{0.5em}
\centerline{
$\displaystyle \delta \asymp \left(\frac{1}{\tau_{\partial, \min} f_{\min} n}\right)^{1/(d+1)}$
\text{ and }
$\displaystyle  \eta \asymp 
\frac{\delta^2}{\tau_{\partial,\min}}
\asymp
\tau_{\partial,\min}
\left(\frac{1}{f_{\min} \tau_{\partial, \min}^d n}\right)^{2/(d+1)}
$.
}
As $\TV(P_0,P_1) \leq 1$, this can only be done when $f_{\min} \delta^{d-1} \eta \lesssim 1$, i.e. $n \gtrsim 1/(f_{\min}\tau_{\partial,\min}^d)$.
}
\label{fig:hypotheses_le_cam_boundary}
\end{figure}

See \Cref{sec:model_stability} for the proof.
Equipped with \Cref{prop:diffeomorphism_stability,prop:diffeomorphism_convexity_stability}, we build hypotheses as shown in \Cref{fig:hypotheses_le_cam_boundary}. The formal statement goes as follows.

\begin{restatable}[Convex Hypotheses]{proposition}{prophypotheseswithboundary}
	\label{prop:hypotheses_with_boundary}
	Assume that $f_{\min} \leq c_d/\tau_{\partial, \min}^d$ and $c'_d/\tau_{\partial, \min}^d \leq f_{\max}$ for some small enough $c_d,(c'_d)^{-1}>0$. 
	
	Then for all $n \geq C_d/(f_{\min}\tau_{\partial, \min}^d)$, 
	there exist $P_0,P_1 \in \mathcal{P}^{d,D}_{\infty,\tau_{\partial,\min}}(f_{\min},f_{\max})$ with convex supports $M_0$ and $M_1$ such that
	\begin{align*}
	\TV(P_0,P_1) \leq \frac{1}{n}
	\text{~~~~and~~~~}
	\dHaus(\partial M_0, \partial M_1)
	&=
	\dHaus(M_0,M_1)
	\geq 
	C'_d \tau_{\partial,\min}
	\left(\frac{1}{f_{\min} \tau_{\partial,\min}^d n}\right)^{2/(d+1)}.
	\end{align*}
\end{restatable}

See \Cref{sec:construction-of-hypotheses} for the construction of these hypotheses. 
We are finally in position to prove \Cref{thm:boundary-main-lower-bound} and \Cref{thm:manifold-main-lower-bound}~\ref{item:minimax_lower_bound_with_boundary}.
    
    \begin{proof}[Proofs of \Cref{thm:boundary-main-lower-bound} and \Cref{thm:manifold-main-lower-bound}~\ref{item:minimax_lower_bound_with_boundary}]
    The proof follows the lines of that of \Cref{thm:manifold-main-lower-bound}~\ref{item:minimax_lower_bound_boundaryless} mutatis mutandis. That is, by setting $\mathcal{P}:=\mathcal{P}^{d,D}_{\infty,\tau_{\partial,\min}}(f_{\min},f_{\max})$, $n_0 :=  \lceil C_d/(f_{\min} \tau_{\partial, \min}^d) \rceil$ where $C_d>0$ is the constant of \Cref{prop:hypotheses_with_boundary}, and applying \Cref{lem:le_cam_applied}~\ref{item:le_cam_applied-manifold} and~\ref{item:le_cam_applied-boundary} with the hypotheses $P_0$ and $P_1$  of \Cref{prop:hypotheses_with_boundary}.
    \end{proof}

%% file: lower-bounds-main-tools.tex
\section{Main tools for the minimax lower bounds}
	\label{sec:appendix:lower-bounds}
	
\subsection{Stability of the model}
\label{sec:model_stability}

\subsubsection{Reach bounds}
To prove \Cref{prop:diffeomorphism_stability}, we will use the following general reach stability result.

\begin{lemma}[{\cite[Theorem 4.19]{Federer59}}]
\label{lem:reach_stability}
Let $S \subset \mathbb{R}^D$ with $\tau_S \geq \tau_0 > 0$, and $ \Phi  : \mathbb{R}^D \longrightarrow \mathbb{R}^D$ be a $\mathcal{C}^1$-diffeomorphism such that $ \Phi $,$ \Phi ^{-1}$, and $d  \Phi $ are Lipschitz, with Lipschitz constants $K$,$N$ and $R$ respectively, then
\[ 
\tau_{\Phi(S)} 
\geq 
\dfrac{\tau_0}{(K + R\tau_0)N^2}. 
\]
\end{lemma}

\propdiffeomorphismstability*

\begin{proof}[Proof of \Cref{prop:diffeomorphism_stability}]
First note that since $\sup_x \normop{d_x \Phi - I_D} < 1$, $d_x \Phi$ is invertible for all $x \in \R^D$, so that $\Phi$ is a local diffeomorphism in the neighborhood of $x$. 
In addition, $\lim_{\norm{x} \rightarrow \infty}\norm{\Phi(x)} = \infty$, so that the Hadamard-Cacciopoli theorem \cite{DeMarco94} asserts that $\Phi$ is a global diffeomorphism of $\R^D$.	

Now, for short, let us write $M' = \Phi\left(M\right)$. As $\Phi$ is a global diffeomorphism of $\R^D$, $M'$ is a $d$-dimensional submanifold: indeed, using notation of \Cref{def:manifold_with_boundary}, any local $\mathcal{C}^2$ parametrization $\Psi_p$ of $M$ at $p \in M$ lifts to the local $\mathcal{C}^2$ parametrization $\tilde{\Psi}_{\Phi(p)} = \Phi \circ \Psi_p$ of $M'$ at $\Phi(p) \in M'$.
In particular, $\partial M' = \Phi(\partial M\bigr)$.
Moreover, $\Phi$ is $\norm{d \Phi}_{op} \leq (1 + \norm{  I_D - d \Phi }_{op})$-Lipschitz, $\Phi^{-1}$ is $\norm{d \Phi^{-1}}_{op} \leq ({1-\norm{ I_D - d \Phi}_{op}})^{-1}$-Lipschitz, and $d \Phi$ is $\norm{d^2 \Phi }_{op}$-Lipschitz. 
Hence, \Cref{lem:reach_stability} applied with $S = M$ yields
\begin{align*}
\tau_{M'} 
&\geq 
\frac{\tau_{M}(1 - \normop{  I_D - d \Phi })^2}{ \norm{d^2 \Phi }_{op} \tau_{M} + (1 + \norm{  I_D - d \Phi }_{op})} 
\geq
\tau_{M}/2
\geq
\tau_{\min}/2
,
\end{align*}
where the second inequality used that $\normop{I_D - d \Phi} \leq 1/10$ and $\normop{d^2 \Phi}\tau_{M} \leq 1/2$.
Similarly, if the boundary $S = \partial M$ is not empty and $\normop{d^2 \Phi}\tau_{\partial M} \leq 1/2$, we get
\begin{align*}
\tau_{\partial M'}
=
\tau_{\Phi(\partial M)}
\geq
\tau_{\partial M}/2
\geq
\tau_{\partial,\min}/2,
\end{align*}
and otherwise $\tau_{\partial M'} = \tau_{\emptyset} = \infty \geq \tau_{\partial,\min}/2$,
which concludes the proof.
\end{proof}

\subsubsection{Strict convexity}

To prove \Cref{prop:diffeomorphism_convexity_stability}, we will use the following non-standard characterization of convexity for full-dimensional domains.

\begin{lemma}
	\label{lem:convexity_characterisation}
	Let $C \subset \R^d$ be a compact domain with $\interior{C} \neq \emptyset$, that has a $\mathcal{C}^2$ boundary $\bar{\partial} C$. Assume that:
	\begin{itemize}
	\item
	for all $x \in \bar{\partial} C$,
	$\bar{\partial} C \setminus \{x\}$ is connected;
	\item
	for all $x, y \in \bar{\partial} C$,
	$
	\dd(y-x,T_x \bar{\partial} C) 
	> 
	0
	$
	as soon as $x \neq y$.
	\end{itemize}
	Then $C$ is convex.
\end{lemma}

\begin{proof}[Proof of \Cref{lem:convexity_characterisation}]
	Let us prove the contrapositive. To this aim, assume that $C$ is not convex, meaning that $\tau_C < \infty$. We will prove the existence of points $x,\tilde{y} \in \bar{\partial} C$ such that $\dd(\tilde{y}-x,T_x \bar{\partial} C) = 0$.
	
	From \cite[Theorem 4.18]{Federer59}, there exist $x \neq y \in C$ such that $\dd(y-x,Tan(x,C)) > 0$. But for all $x \in \interior{C}$, $Tan(x,C) = \R^d$, so that $x \in \bar{\partial} C$ necessarily.
	From here, \Cref{prop:cones_of_manifold} asserts that $Tan(x,C)$ is a half-space with
	$
	\Span(Tan(x,C))
	=
	\R^d
	=
	T_x \bar{\partial} C \Operp \Span(\eta_x)
	$
	and
	$Tan(x,C) = \set{ \inner{\eta_x}{.} \leq 0}$,
	for some unit vector $\eta_x \in \R^d$.
	Using this representation, for all $z \in C$, we have 
	$
	\dd(z-x,Tan(x,C)) 
	=
	\inner{z-x}{\eta_x}_+ 
	$
	and
	$
	\dd(z-x,T_x \bar{\partial} C) 
	= 
	|\inner{z-x}{\eta_x}|
	$
	.
	
	On one hand, we have seen that the continuous map $C \ni y \mapsto \inner{y-x}{\eta_x}_+$ takes a positive value. Hence, by compactness of $C$, it attains its maximum at some $y_0 \in C$ with $\inner{y_0-x}{\eta_x}_+ = \inner{y_0-x}{\eta_x} > 0$. 
	But for $\delta \in \R^d$ small enough, $\inner{y_0+\delta-x}{\eta_x}_+ = \inner{y_0+\delta-x}{\eta_x} = \inner{y_0-x}{\eta_x} + \inner{\delta}{\eta_x}$, so $y_0$ must belong to $\bar{\partial} C$ as otherwise, $y_0$ would belong to $\interior{C}$ and one could increase the value of $\inner{\cdot-x}{\eta_x}_+$ locally around $y_0$ and still stay in $C$.
	
	On the other hand, if we assumed that for all $y \in \bar{\partial} C$, $\inner{y-x}{\eta_x} \geq 0$ this would lead to a contradiction. Indeed, this inequality would extend to all the points $z \in C$: since $C$ is compact, for all $z \in \interior{C}$ and $v \in \R^d \setminus \{0\}$, $\{z + \lambda v, \lambda \in \R \} \cap C$ is a non-empty compact set, so there exist $\lambda_- < \lambda_+$ such that for all $\lambda \in [\lambda_-,\lambda_+]^c$, $z + \lambda v \notin C$ and $y_\pm = z + \lambda_\pm v \in C$. In particular, 
	$y_\pm \in \bar{\partial} C$ and $z \in [y_-,y_+] \subset \R^d$. This shows that $z \in C$ can be written as linear combination of elements $y_\pm \in \bar{\partial} C$ and as a result the assumption $\inner{y_\pm - x}{\eta_x} \geq 0$ would yield $\inner{z-x}{\eta_x} \geq 0$. This is a contradiction, since by definition of $Tan(x,C) \ni -\eta_x$ (\Cref{def:tangent_normal_cones}), there exists $\tilde{z} \in C \setminus \{x \}$ such that $\norm{-\eta_x - \frac{\tilde{z}-x}{\norm{\tilde{z}-x}}} < \frac{1}{2}$ and in particular, $\inner{\tilde{z}-x}{\eta_x} < 0$.
	This ends proving that there exists $y_1 \in \bar{\partial} C$ such that $\inner{y_1-x}{\eta_x} < 0$.
	
	Summing everything up, we have shown that the continuous map $\bar{\partial} C \setminus \{x\} \ni y \mapsto \inner{y-x}{\eta_x}$ takes both a positive and a negative value on its connected domain $\bar{\partial} C \setminus \{x\}$. Hence, it must vanish at some point $\tilde{y} \in \bar{\partial} C \setminus \{x\}$, meaning that $x \neq \tilde{y} \in \bar{\partial} C$ and $\dd(y-x,T_x \bar{\partial} C) = 0$, which concludes the proof.
\end{proof}

\propdiffeomorphismconvexitystability*

\begin{proof}[Proof of \Cref{prop:diffeomorphism_convexity_stability}]
First, from \Cref{lem:convexity_characterisation}, we get that $C$ is convex. Furthermore, as in the proof of \Cref{prop:diffeomorphism_stability}, note that the assumptions $\normop{d \Phi - I_d} < 1$ and $\lim_{\norm{x} \rightarrow \infty}\norm{\Phi(x)} = \infty$ yield that $\Phi$ is a global diffeomorphism of $\R^d$, using the Hadamard-Cacciopoli theorem \cite{DeMarco94}.
Hence, writing $C' = \Phi(C)$, we get that $C'$ is a compact domain with $\interior{C'} \neq \emptyset$, that has a connected $\mathcal{C}^2$ boundary $\bar{\partial} C'$. In addition, $\bar{\partial} C' = \Phi(\bar{\partial} C)$ and for all $x' = \Phi(x) \in \bar{\partial} C'$, $T_{x'} \bar{\partial} C' = d_x \Phi \bigl( T_x \bar{\partial} C \bigr)$.

Now, for all $x,y \in \bar{\partial} C$ and $u \in T_x \bar{\partial} C$, Taylor's theorem and the assumption $\dd(y-x,T_x \bar{\partial} C) \geq A \norm{y-x}^2$ yield
\begin{align*}
\norm{d_x \Phi.(y-x) - d_x \Phi.u}
&\geq
\normop{d \Phi^{-1}}^{-1} \norm{(y-x)-u}
\\
&\geq
\normop{d \Phi^{-1}}^{-1}
\dd(y-x,T_x \bar{\partial} C)
\\
&\geq
\normop{d \Phi^{-1}}^{-1}
A\norm{y-x}^2
\\
&\geq
(1- \normop{I_d-d \Phi})
A\norm{y-x}^2
\\
&\geq
(9A/10)\norm{y-x}^2.
\end{align*}
At second order, Taylor's theorem writes
\begin{align*}
\norm{\Phi(y)-\Phi(x) - d_x \Phi.(y-x)}
\leq
{\normop{d^2 \Phi}}\norm{y-x}^2/2
.
\end{align*}
As a result, for all $x' \neq y' \in \bar{\partial} C'$, writing $x' = \Phi(x)$ and $y' = \Phi(y)$ we have $x \neq y$ as $\Phi^{-1}$ is one-to-one, and
\begin{align*}
\dd(y'-x',T_{x'} \bar{\partial} C')
&=
\inf_{u \in T_{x} \bar{\partial} C}
\norm{ \Phi(y) - \Phi(x) - d_x \Phi.u}
\\
&\geq
\inf_{u \in T_{x} \bar{\partial} C}
\left\{
\norm{d_x \Phi.(y-x) - d_x \Phi.u}
-
\norm{\Phi(y)-\Phi(x) - d_x \Phi.(y-x)}
\right\}
\\
&\geq
\left(
9A/10
-
\normop{d^2 \Phi}/2
 \right)
\norm{y-x}^2
\\
&>
0,
\end{align*}
since $\normop{d^2 \Phi} \leq A < 9A/5$.
From \Cref{lem:convexity_characterisation}, $C'$ is hence convex.
\end{proof}

\subsection{Construction of hypotheses}
\label{sec:construction-of-hypotheses}

Throughout this section, we will use a smooth localizing bump-type function $\phi : \R^D \to \R$ to build local variations of manifolds.
The following result gathers differential estimates, and can be shown using elementary differential calculus.

\begin{proposition}
	\label{prop:kernel_properties}
	The localizing function defined as
	\begin{align*}
		 \phi: \R^D &\longrightarrow \R
		 \\
		 x &\longmapsto \exp\left(-{\norm{x}^2}/{(1-\norm{x}^2)} \right) \indicator{\B(0,1)}(x)
	\end{align*}
is $\mathcal{C}^\infty$ smooth, equal to $0$ outside  $\B(0,1)$, satisfies $0 \leq \phi \leq 1$, $\phi(0) = 1$,
	\begin{align*}
	\normop{d \phi}
	:=
	\sup_{x \in \R^D}
	\normop{d_x \phi}
	\leq 5/2
	\text{ and }
	\normop{d^2 \phi}
	:=
	\sup_{x \in \R^D}
	\normop{d_x^2 \phi}
	\leq 23
	.
	\end{align*}
\end{proposition}

\subsubsection{Hypotheses with empty boundary}

	The proof of \Cref{prop:hypotheses_empty_boundary} follows that of~\cite[Lemma 5]{Aamari19b}, and provides a result similar to \cite[Theorem~6]{Genovese12a} in essence. 
	We include it below for sake of completeness and to keep track of explicit constants.
	
\prophypothesesemptyboundary*
	
	\begin{proof}[Proof of \Cref{prop:hypotheses_empty_boundary}]
	
		We let $R = 2 \tau_{\min}$, and $M_0 = \Sphere^{d}(0,R) \times \{0\}^{D-d-1}$ be a $d$-dimensional sphere of radius $R$ embedded in $\R^{d+1}\times \{{0}\}^{D-(d+1)}$.
	Clearly, $\partial M_0 = \emptyset$ (meaning that $\tau_{\partial M_0} = \infty$) and $\tau_{M_0} = R = 2 \tau_{\min}$.
	
		Let $e_1 = (1, 0, \ldots, 0)$ denote the first vector of the canonical basis of $\R^D$, and $x_0 = R e_1 \in M_0$.
	For $\delta>0$ to be specified later, consider the probability distribution $P_0$ having the following density with respect to the $d$-dimensional Hausdorff measure $\Haus^d$:
	\begin{align*}
		f_0(x)
		&=
		2 f_{\min} \indicator{M_0 \cap \B(x_0,\delta)}(x)
		+
		\frac{1-2 f_{\min} \Haus^d(M_0 \cap \B(x_0,\delta))}{\Haus^d(M_0 \cap \B(x_0,\delta)^c)
}
		\indicator{M_0 \cap \B(x_0,\delta)^c}(x)
		,
	\end{align*}
	for all $x \in \R^D$.
	Clearly, $P_0$ has support $M_0$ as soon as $2 f_{\min} \Haus^d(M_0 \cap \B(x_0,\delta))< 1$.
	In addition, writing $\sigma_d$ for the volume of the $d$-dimensional unit Euclidean sphere,
	\begin{align*}
		\frac{1-2 f_{\min} \Haus^d(M_0 \cap \B(x_0,\delta))}{\Haus^d(M_0 \cap \B(x_0,\delta)^c)}
		&\geq
		\frac{
			1-2 f_{\min} \Haus^d(M_0 \cap \B(x_0,\delta))
			}{
			\Haus^d(M_0)
		}
		\\
		&=
		\frac{
			1-2 f_{\min}
			R^d \Haus^d\left( \B_{\Sphere^{d}}(0,2\arcsin(\delta/(2R))\right)
			}{
			\sigma_d R^d
		}
		\\
		&\geq
		\frac{1}{\sigma_d R^d}
		-
		2 f_{\min} \left(\frac{\delta}{R}\right)^d
		.
	\end{align*}
	As a result, $f_0 \geq 2 f_{\min}$ over $M_0$ as soon as $(\sigma_d (2 \tau_{\min})^d)^{-1} \geq 4 f_{\min}$ and $\delta \leq 2 \tau_{\min}$.
	To upper bound $f_0$ on $M_0$, we note that $2f_{\min} \leq f_{\max}/2$ as soon as $2c_d \leq c'_d/2$, and that similarly to above, we derive
	\begin{align*}
		\frac{1-2 f_{\min} \Haus^d(M_0 \cap \B(x_0,\delta))}{\Haus^d(M_0 \cap \B(x_0,\delta)^c)}
		&\leq
		\frac{1}{\sigma_d(R^d - \delta^d)}
		\leq
		\frac{2}{\sigma_dR^d}		
	\end{align*}
	as soon as $\delta \leq \tau_{\min}$, which is further upper bounded by $f_{\max}/2$ as soon as $2/(2^d\sigma_d) \leq c'_d/2$.	
	This ends proving that $P_0 \in \mathcal{P}^{d,D}_{\tau_{\min},\infty}(2f_{\min},f_{\max}/2)$.

	We now build $P_1$ by small and smooth ambient perturbation of $P_0$. Namely, for $\eta > 0$ to be specified later, write
	\begin{align*}
	\Phi(x)
	=
	x
	+
	\eta \phi\left( \frac{x - x_0}{\delta} \right) e_1
	,
	\end{align*}
	where $\phi : \R^D \to \R$ is the localizing function of \Cref{prop:kernel_properties}. 
		We let $P_1 = \Phi_\ast P_0$ be the pushforward distribution of $P_0$ by $\Phi$, and $M_1 = \supp(P_1)$.
		
		From \Cref{prop:kernel_properties}, we get that 
	$\Phi$ is $\mathcal{C}^\infty$ smooth, 
	$\normop{d \Phi - I_D} = \frac{\eta}{\delta} \normop{d \phi} \leq \frac{5 \eta}{2 \delta}$,
	and
	$\normop{d^2 \Phi} = \frac{\eta}{\delta^2} \normop{d^2 \phi} \leq \frac{23 \eta}{\delta^2}$.	
	Recalling that $\tau_{M_0} \geq 2 \tau_{\min}$, \Cref{prop:diffeomorphism_stability} asserts that
	$M_1 \in \mathcal{M}^{d,D}_{\tau_{\min},\infty}$ as soon as $\frac{5\eta}{2 \delta} \leq \frac{1}{10}$ and $\frac{23 \eta}{\delta^2} \leq \frac{1}{4 \tau_{\min}}$.
	Furthermore, from \cite[Appendix, Lemma A.6]{Aamari19b}, $P_1$ admits a density $f_1$ with respect to $\Haus^d$ that satisfies 
	\begin{align*}
	f_{\min}
	=
	\inf_{M_0} f_0/2
	\leq
	\inf_{M_1} f_1
	\leq
	\sup_{M_1} f_1
	\leq
	2 \sup_{M_0} f_0
	\leq
	f_{\max}
	\end{align*}
	as soon as $\frac{5\eta}{2 \delta} \leq \frac{1}{3d} \wedge \frac{1}{3(2^{d/2} - 1)}$.
	Hence, under all the above requirements, we finally get that $P_1 \in \mathcal{P}^{d,D}_{\tau_{\min},\infty}(f_{\min},f_{\max})$.

	Now, notice that by construction, $x_0 + \eta e_1 = \Phi(x_0)$ belongs to $M_1=\Phi(M_0)$.
	As a result,
	\begin{align*}
	\dHaus(M_0,M_1)
	&\geq
	\dd(x_0+\eta e_1,M_0)
	=
	\eta
	.
	\end{align*}

	In addition, under the same requirements on $\delta$ and $\eta$ as above, $\Phi$ is a global diffeomorphism of $\R^D$ (\Cref{prop:diffeomorphism_stability}). As it coincides with the identity map on $\B(x_0,\delta)^c$, this implies that $P_0$ and $P_1 = \Phi_\ast P_0$ coincide outside $\B(x_0,\delta)$. 
	Hence,
	\begin{align*}
	\TV(P_0,P_1)
	&=
	\sup_{A \in \mathcal{B}(\R^D)}
	|P_1(A\cap \B(x_0,\delta)) - P_0(A \cap \B(x_0,\delta))|
	\\
	&\leq
	\sup_{A \in \mathcal{B}(\R^D)}
	P_0(A\cap \B(x_0,\delta)) \vee P_1(A \cap \B(x_0,\delta))
	\\
	&\leq
	P_0(\B(x_0,\delta))
	\vee
	P_1(\B(x_0,\delta))
	\\
	&=
	P_0(\B(x_0,\delta))
	\\
	&=
	2f_{\min} \Haus^d\bigl(M_0 \cap \B(x_0,\delta)\bigr)
	\\
	&=
	2f_{\min} R^d \Haus^d\left( \B_{\Sphere^{d}}(0,2\arcsin(\delta/(2R))\right)
	\\
	&\leq
	2\sigma_d f_{\min}\delta^d
	.
	\end{align*}
	Setting $2\sigma_d f_{\min}\delta^d = 1/n$ and $\eta = \frac{\delta}{2^{d+10}} \wedge \frac{\delta^2}{92 \tau_{\min}}$ (which satisfy all the above requirements) then yields the result, since with that choice, $\delta \leq \tau_{\min}$ and $\eta = \frac{\delta^2}{92 \tau_{\min}}$ as soon as $n \geq C_d /(f_{\min} \tau_{\min}^d)$ for some large enough $C_d>0$.
	\end{proof}

\subsubsection{Convex hypotheses (with boundary)}
	The proof of \Cref{prop:hypotheses_with_boundary} is similar to that of \Cref{prop:hypotheses_empty_boundary}.

\prophypotheseswithboundary*

	\begin{proof}[Proof of \Cref{prop:hypotheses_with_boundary}]
	
	Let $R = 2 \tau_{\partial, \min}$, and $M_0 = \B_{\R^d}(0,R) \times \{0\}^{D-d}$ be a $d$-dimensional ball of radius $R$ embedded in $\R^{d}\times \{{0}\}^{D-d}$.
	Clearly, $M_0$ is convex, meaning that $\tau_{M_0} = \infty$, and $\partial M_0 = \Sphere^{d-1}(0,R) \times \{0\}^{D-d}$ has reach $\tau_{\partial M_0} = R$.
	
	Let $e_1 = (1, 0, \ldots, 0)$ denote the first vector of the canonical basis of $\R^D$, and $x_0 = R e_1 \in M_0$.
	For $\delta>0$ to be specified later, consider the probability distribution $P_0$ having the following density with respect to the $d$-dimensional Hausdorff measure $\Haus^d$:
	\begin{align*}
		f_0(x)
		&=
		2 f_{\min} \indicator{M_0 \cap \B(x_0,\delta)}(x)
		+
		\frac{1-2 f_{\min} \Haus^d(M_0 \cap \B(x_0,\delta))}{\Haus^d(M_0 \cap \B(x_0,\delta)^c)
}
		\indicator{M_0 \cap \B(x_0,\delta)^c}(x)
		,
	\end{align*}
	for all $x \in \R^D$.
	We see that $P_0$ has support $M_0$ if $2 f_{\min} \Haus^d(M_0 \cap \B(x_0,\delta))< 1$.
	Denoting by $\omega_d$ the volume of the $d$-dimensional unit Euclidean ball, we derive
	\begin{align*}
		\frac{1-2 f_{\min} \Haus^d(M_0 \cap \B(x_0,\delta))}{\Haus^d(M_0 \cap \B(x_0,\delta)^c)}
		&\geq
		\frac{
			1-2 f_{\min} \Haus^d(M_0 \cap \B(x_0,\delta))
			}{
			\Haus^d(M_0)
		}
		\\
		&\geq
		\frac{
			1-2 f_{\min}
			(\omega_d \delta^d/2)
			}{
			\omega_d R^d
		}
		\\
		&\geq
		\frac{1}{\omega_d R^d}
		-
		f_{\min} \left(\frac{\delta}{R}\right)^d
		.
	\end{align*}
	As a result, $f_0 \geq 2 f_{\min}$ over $M_0$ as soon as $(\omega_d (2 \tau_{\min})^d)^{-1} \geq 4 f_{\min}$ and $\delta \leq 2 \tau_{\min}$,
	To upper bound $f_0$ on $M_0$, we note that $2f_{\min} \leq f_{\max}/2$ as soon as $2c_d \leq c'_d/2$, and that similarly to above, we derive
	\begin{align*}
		\frac{1-2 f_{\min} \Haus^d(M_0 \cap \B(x_0,\delta))}{\Haus^d(M_0 \cap \B(x_0,\delta)^c)}
		&\leq
		\frac{1}{\omega_d(R^d - \delta^d/2)}
		\leq
		\frac{2}{\omega_dR^d}	
	\end{align*}
	as soon as $\delta \leq R = 2 \tau_{\min}$, which is further upper bounded by $f_{\max}/2$ as soon as $2/(2^d\omega_d) \leq c'_d/2$.	
	In all, we have $P_0 \in \mathcal{P}^{d,D}_{\infty,\tau_{\partial, \min}}(2f_{\min},f_{\max}/2)$.
	
	Now, to build $P_1$, let $\eta > 0$ be a parameter to be specified later, and write
	\begin{align*}
	\Phi(x)
	=
	x
	+
	\eta \phi\left( \frac{x - x_0}{\delta} \right) e_1
	,
	\end{align*}
	where $\phi : \R^D \to \R$ is the localizing function of \Cref{prop:kernel_properties}. 
		We let $P_1 = \Phi_\ast P_0$ be the pushforward distribution of $P_0$ by $\Phi$, and $M_1 = \supp(P_1)$.
		Note by now that if $\delta \leq R$, we have $M_0 \subset M_1$.
		
		From \Cref{prop:kernel_properties}, 
	we get that $\Phi$ is $\mathcal{C}^\infty$ smooth, 
	$\normop{d \Phi - I_D} = \frac{\eta}{\delta} \normop{d \phi} \leq \frac{5 \eta}{2 \delta}$,
	and
	$\normop{d^2 \Phi} = \frac{\eta}{\delta^2} \normop{d^2 \phi} \leq \frac{23 \eta}{\delta^2}$.	
	It is also clear that $\lim_{\norm{x} \to \infty} \norm{\Phi(x)} = \infty$.
	Hence, recalling that $\tau_{\partial M_0} \geq 2 \tau_{\partial, \min}$, \Cref{prop:diffeomorphism_stability} asserts that
	$\tau_{\partial M_1}  \geq \tau_{\partial, \min}$ as soon as $\frac{5\eta}{2 \delta} \leq \frac{1}{10}$ and $\frac{23 \eta}{\delta^2} \leq \frac{1}{4 \tau_{\partial,\min}}$.
	In addition, as $\Phi$ preserves $\R^{d} \times \{0\}^{D-d}$, both $M_0$ and $M_1$ can be seen as compact domains of $\R^d$ with non-empty interior. In this $d$-plane $\R^d \times \{0\}^{D-d} \cong \R^d$, $M_0$ has a $\mathcal{C}^2$ (topological) boundary $\bar{\partial} M_0 = \Sphere^{d-1}(0,R)$,  the set $\bar{\partial} M_0 \setminus \{x\}$ is connected for all $x \in \bar{\partial} M_0$ (note that for $d=1$, this set is only reduced to a point), and for all $x,y \in \bar{\partial} M_0$,
	$
	\dd(y-x,T_x \bar{\partial} M_0) 
	=
	\frac{1}{4 \tau_{\partial, \min}} \norm{y-x}^2
	$
	.
	As a result, \Cref{prop:diffeomorphism_convexity_stability} applied with $k=d$ asserts that $M_1 = \Phi(M_0)$ remains convex as soon as  $\frac{5\eta}{2 \delta} \leq \frac{1}{10}$ and $\frac{23 \eta}{\delta^2} \leq \frac{1}{4 \tau_{\partial,\min}}$.
	This ends proving that $M_0,M_1 \in \mathcal{M}^{d,D}_{\infty,\tau_{\partial,\min}}$ under the above requirements.

	Furthermore, from \cite[Appendix, Lemma A.6]{Aamari19b}, we get that $P_1$ admits a density $f_1$ with respect to $\Haus^d$ that satisfies 
	\begin{align*}
	f_{\min}
	=
	\inf_{M_0} f_0/2
	\leq
	\inf_{M_1} f_1
	\leq
	\sup_{M_1} f_1
	\leq
	2 \sup_{M_0} f_0
	\leq
	f_{\max}
	\end{align*}
	as soon as $\frac{5\eta}{2 \delta} \leq \frac{1}{3d} \wedge \frac{1}{3(2^{d/2} - 1)}$.
	Hence, under all the above requirements, we have that both $P_0$ and  $P_1$ belong to the model $\mathcal{P}^{d,D}_{\infty,\tau_{\partial, \min}}(f_{\min},f_{\max})$.
	
	Further analyzing the properties of $f_1$, let $y \in M_1 \cap \B(x_0,\delta)$. As the diffeomorphism $\Phi$ maps $\B(x_0,\delta)$ onto itself, $y = \Phi(x)$ for a unique $x \in M_0 \cap \B(x_0,\delta)$. 
	Hence, applying \cite[Appendix, Lemma A.6]{Aamari19b} again we get
	\begin{align*}
		|f_1(y) - 2f_{\min}|
		&=
		|f_1(y) - f_0(x)|
		\\
		&\leq
		f_0(x) \left(\frac{3d}{2} \vee 3(2^{d/2}-1)\right) 
		\normop{d \Phi - I_D}
		\\
		&=
		\frac{2^{d+10} f_{\min} \eta}{\delta}
		,
	\end{align*}	
	provided that $\frac{5\eta}{2 \delta} < \frac{1}{3}$.
	From this bound, we also read that $f_1 \leq 3 f_{\min}$ on $M_1 \cap \B(x_0,\delta)$ as soon as $\frac{2^{d+10} \eta}{\delta} \leq 1$.	
	We can now move forward and prove the result.

	First, notice that by construction $x_0 + \eta e_1 = \Phi(x_0)$ belongs to $\partial M_1=\Phi(\partial M_0)$.
	As a result,
	\begin{align*}
	\dHaus(\partial M_0,\partial M_1)
	=
	\dHaus(M_0,M_1)	
	&\geq
	\dd(x_0+\eta e_1,M_0)
	=
	\eta
	.
	\end{align*}
	Second, under the same requirements on $\delta$ and $\eta$ as above, $\Phi$ is a global diffeomorphism of $\R^D$ (\Cref{prop:diffeomorphism_stability}). As it coincides with the identity map on $\B(x_0,\delta)^c$, it implies that $P_0$ and $P_1 = \Phi_\ast P_0$ coincide outside $\B(x_0,\delta)$.
	Applying the second formula of \Cref{def:total_variation} with the $\sigma$-finite dominating measure $\mu = \indicator{\R^d \times \{0\}^{D-d}}\Haus^d$, we hence get	
	\begin{align*}
	\TV(P_0,P_1)
	&=
	\frac{1}{2}
	\int_{\B(x_0,\delta)\cap (M_0 \cup M_1)}
	|f_1 - f_0|
	d \Haus^d
	\\
	&=
	\frac{1}{2}
	\int_{\B(x_0,\delta)\cap M_0}
	|f_1 - 2 f_{\min}|
	d \Haus^d
	+
	\frac{1}{2}
	\int_{\B(x_0,\delta)\cap (M_1\setminus M_0)}
	f_1
	d \Haus^d
	\\
	&\leq
	\frac{2^{d+10} f_{\min} \eta}{2 \delta}
	\Haus^d\bigl(\B(x_0,\delta)\cap M_0\bigr)
	+
	\frac{3 f_{\min}}{2}
	\Haus^d\bigl(\B(x_0,\delta)\cap (M_1 \setminus M_0)\bigr)
	.
	\end{align*}
	Furthermore, by construction, 
	$
	\Haus^d\bigl(\B(x_0,\delta)\cap M_0\bigr)
	\leq 
	\omega_d \delta^d/2
	$ 
	and
	$
	\Haus^d\bigl(\B(x_0,\delta)\cap (M_1 \setminus M_0)\bigr)
	\leq
	C'_d \delta^{d-1} \eta
	$
	,
	so that
	\begin{align*}
	\TV(P_0,P_1)
	&\leq
	C''_d f_{\min}
	\delta^{d-1}
	\eta
	.
	\end{align*}
	Finally, setting $C''_d f_{\min}	\delta^{d-1} \eta = 1/n$ and $\eta = \frac{\delta}{2^{d+10}} \wedge \frac{\delta^2}{92 \tau_{\partial, \min}}$ (which satisfy all the above requirements) then yields the result, since with that choice, $\delta \leq \tau_{\min}$ and $\eta = \frac{\delta^2}{92 \tau_{\partial,\min}}$ as soon as $n \geq \tilde{C}_d /(f_{\min} \tau_{\partial, \min}^d)$ for some large enough $C_d>0$.
	\end{proof}

%% file: main.bbl
\begin{thebibliography}{10}

\bibitem{Aamari19}
Eddie Aamari, Jisu Kim, Fr\'{e}d\'{e}ric Chazal, Bertrand Michel, Alessandro
  Rinaldo, and Larry Wasserman.
\newblock Estimating the reach of a manifold.
\newblock {\em Electron. J. Stat.}, 13(1):1359--1399, 2019.

\bibitem{Aamari18}
Eddie Aamari and Cl\'{e}ment Levrard.
\newblock Stability and minimax optimality of tangential {D}elaunay complexes
  for manifold reconstruction.
\newblock {\em Discrete Comput. Geom.}, 59(4):923--971, 2018.

\bibitem{Aamari19b}
Eddie Aamari and Cl\'{e}ment Levrard.
\newblock Nonasymptotic rates for manifold, tangent space and curvature
  estimation.
\newblock {\em Ann. Statist.}, 47(1):177--204, 2019.

\bibitem{Aaron16b}
Catherine Aaron and Olivier Bodart.
\newblock Local convex hull support and boundary estimation.
\newblock {\em J. Multivariate Anal.}, 147:82--101, 2016.

\bibitem{Aaron16a}
Catherine Aaron and Alejandro Cholaquidis.
\newblock On boundary detection.
\newblock {\em Ann. Inst. Henri Poincar\'{e} Probab. Stat.}, 56(3):2028--2050,
  2020.

\bibitem{Aaron20}
Catherine Aaron, Alejandro Cholaquidis, and Ricardo Fraiman.
\newblock {Estimation of surface area}.
\newblock {\em Electronic Journal of Statistics}, 16(2):3751 -- 3788, 2022.

\bibitem{Aizenbud21}
Yariv {Aizenbud} and Barak {Sober}.
\newblock {Non-Parametric Estimation of Manifolds from Noisy Data}.
\newblock {\em arXiv e-prints}, page arXiv:2105.04754, May 2021.

\bibitem{Belkin06}
Mikhail Belkin, Partha Niyogi, and Vikas Sindhwani.
\newblock Manifold regularization: a geometric framework for learning from
  labeled and unlabeled examples.
\newblock {\em J. Mach. Learn. Res.}, 7:2399--2434, 2006.

\bibitem{Berenfeld20}
Cl{\'e}ment Berenfeld, John Harvey, Marc Hoffmann, and Krishnan Shankar.
\newblock Estimating the reach of a manifold via its convexity defect function.
\newblock {\em Discrete {\&} Computational Geometry}, Jun 2021.

\bibitem{Berry17}
Tyrus Berry and Timothy Sauer.
\newblock Density estimation on manifolds with boundary.
\newblock {\em Comput. Statist. Data Anal.}, 107:1--17, 2017.

\bibitem{Boissonnat14}
Jean-Daniel Boissonnat and Arijit Ghosh.
\newblock Manifold reconstruction using tangential {D}elaunay complexes.
\newblock {\em Discrete Comput. Geom.}, 51(1):221--267, 2014.

\bibitem{Boissonnat09}
Jean-Daniel Boissonnat, Leonidas~J. Guibas, and Steve~Y. Oudot.
\newblock Manifold reconstruction in arbitrary dimensions using witness
  complexes.
\newblock {\em Discrete Comput. Geom.}, 42(1):37--70, 2009.

\bibitem{Boissonnat19}
Jean-Daniel Boissonnat, Andr\'{e} Lieutier, and Mathijs Wintraecken.
\newblock The reach, metric distortion, geodesic convexity and the variation of
  tangent spaces.
\newblock {\em J. Appl. Comput. Topol.}, 3(1-2):29--58, 2019.

\bibitem{Boissonnat20}
Jean-Daniel Boissonnat and Mathijs Wintraecken.
\newblock {The Topological Correctness of PL-Approximations of Isomanifolds}.
\newblock In Sergio Cabello and Danny~Z. Chen, editors, {\em 36th International
  Symposium on Computational Geometry (SoCG 2020)}, volume 164 of {\em Leibniz
  International Proceedings in Informatics (LIPIcs)}, pages 20:1--20:18,
  Dagstuhl, Germany, 2020. Schloss Dagstuhl--Leibniz-Zentrum f{\"u}r
  Informatik.

\bibitem{Bredon93}
Glen~E. Bredon.
\newblock {\em Topology and geometry}, volume 139 of {\em Graduate Texts in
  Mathematics}.
\newblock Springer-Verlag, New York, 1993.

\bibitem{Burago01}
Dmitri Burago, Yuri Burago, and Sergei Ivanov.
\newblock {\em A course in metric geometry}, volume~33 of {\em Graduate Studies
  in Mathematics}.
\newblock American Mathematical Society, Providence, RI, 2001.

\bibitem{Calder22}
Jeff Calder, Sangmin Park, and Dejan Slep{\v{c}}ev.
\newblock Boundary estimation from point clouds: Algorithms, guarantees and
  applications.
\newblock {\em Journal of Scientific Computing}, 92(2):56, Jul 2022.

\bibitem{Chazal13}
Fr\'{e}d\'{e}ric Chazal, Marc Glisse, Catherine Labru\`ere, and Bertrand
  Michel.
\newblock Convergence rates for persistence diagram estimation in topological
  data analysis.
\newblock {\em J. Mach. Learn. Res.}, 16:3603--3635, 2015.

\bibitem{Chazal17}
Fr{\'e}d{\'e}ric {Chazal} and Bertrand {Michel}.
\newblock {An introduction to Topological Data Analysis: fundamental and
  practical aspects for data scientists}.
\newblock {\em arXiv e-prints}, page arXiv:1710.04019, October 2017.

\bibitem{Cuevas04}
Antonio Cuevas and Alberto Rodr{\'{\i}}guez-Casal.
\newblock On boundary estimation.
\newblock {\em Adv. in Appl. Probab.}, 36(2):340--354, 2004.

\bibitem{DeMarco94}
Giuseppe De~Marco, Gianluca Gorni, and Gaetano Zampieri.
\newblock Global inversion of functions: an introduction.
\newblock {\em NoDEA Nonlinear Differential Equations Appl.}, 1(3):229--248,
  1994.

\bibitem{Dey09b}
Tamal~K Dey, Kuiyu Li, Edgar~A Ramos, and Rephael Wenger.
\newblock Isotopic reconstruction of surfaces with boundaries.
\newblock In {\em Computer Graphics Forum}, volume~28, pages 1371--1382. Wiley
  Online Library, 2009.

\bibitem{Divol20}
Vincent {Divol}.
\newblock {Minimax adaptive estimation in manifold inference}.
\newblock {\em arXiv e-prints}, page arXiv:2001.04896, January 2020.

\bibitem{Divol21}
Vincent {Divol}.
\newblock {Reconstructing measures on manifolds: an optimal transport
  approach}.
\newblock {\em arXiv e-prints}, page arXiv:2102.07595, February 2021.

\bibitem{DoCarmo92}
Manfredo Perdig\~{a}o do~Carmo.
\newblock {\em Riemannian geometry}.
\newblock Mathematics: Theory \& Applications. Birkh\"{a}user Boston, Inc.,
  Boston, MA, 1992.
\newblock Translated from the second Portuguese edition by Francis Flaherty.

\bibitem{Dumbgen96}
Lutz D{\"u}mbgen and G{\"u}nther Walther.
\newblock Rates of convergence for random approximations of convex sets.
\newblock {\em Adv. in Appl. Probab.}, 28(2):384--393, 1996.

\bibitem{Federer59}
Herbert Federer.
\newblock Curvature measures.
\newblock {\em Trans. Amer. Math. Soc.}, 93:418--491, 1959.

\bibitem{Federer69}
Herbert Federer.
\newblock {\em Geometric measure theory}.
\newblock Die Grundlehren der mathematischen Wissenschaften, Band 153.
  Springer-Verlag New York Inc., New York, 1969.

\bibitem{Fefferman19}
Charles {Fefferman}, Sergei {Ivanov}, Matti {Lassas}, and Hariharan
  {Narayanan}.
\newblock {Fitting a manifold of large reach to noisy data}.
\newblock {\em arXiv e-prints}, page arXiv:1910.05084, October 2019.

\bibitem{Genovese12b}
Christopher~R. Genovese, Marco Perone-Pacifico, Isabella Verdinelli, and Larry
  Wasserman.
\newblock Manifold estimation and singular deconvolution under {H}ausdorff
  loss.
\newblock {\em Ann. Statist.}, 40(2):941--963, 2012.

\bibitem{Genovese12a}
Christopher~R. Genovese, Marco Perone-Pacifico, Isabella Verdinelli, and Larry
  Wasserman.
\newblock Minimax manifold estimation.
\newblock {\em J. Mach. Learn. Res.}, 13:1263--1291, 2012.

\bibitem{Har11}
Sariel Har-Peled.
\newblock {\em Geometric approximation algorithms}, volume 173 of {\em
  Mathematical Surveys and Monographs}.
\newblock American Mathematical Society, Providence, RI, 2011.

\bibitem{Hastie09}
Trevor Hastie, Robert Tibshirani, and Jerome Friedman.
\newblock {\em The elements of statistical learning}.
\newblock Springer Series in Statistics. Springer, New York, second edition,
  2009.
\newblock Data mining, inference, and prediction.

\bibitem{Hatcher02}
Allen Hatcher.
\newblock {\em Algebraic topology}.
\newblock Cambridge University Press, Cambridge, 2002.

\bibitem{Hirsch76}
Morris~W. Hirsch.
\newblock {\em Differential topology}.
\newblock Graduate Texts in Mathematics, No. 33. Springer-Verlag, New
  York-Heidelberg, 1976.

\bibitem{Kim2015}
Arlene K.~H. Kim and Harrison~H. Zhou.
\newblock Tight minimax rates for manifold estimation under {H}ausdorff loss.
\newblock {\em Electron. J. Stat.}, 9(1):1562--1582, 2015.

\bibitem{Lee07}
John~A. Lee and Michel Verleysen.
\newblock {\em Nonlinear dimensionality reduction}.
\newblock Information Science and Statistics. Springer, New York, 2007.

\bibitem{Lee11}
John~M. Lee.
\newblock {\em Introduction to topological manifolds}, volume 202 of {\em
  Graduate Texts in Mathematics}.
\newblock Springer, New York, second edition, 2011.

\bibitem{Maggioni16}
Mauro Maggioni, Stanislav Minsker, and Nate Strawn.
\newblock Multiscale dictionary learning: non-asymptotic bounds and robustness.
\newblock {\em J. Mach. Learn. Res.}, 17:Paper No. 2, 51, 2016.

\bibitem{Tsybakov95}
E.~Mammen and A.~B. Tsybakov.
\newblock Asymptotical minimax recovery of sets with smooth boundaries.
\newblock {\em Ann. Statist.}, 23(2):502--524, 1995.

\bibitem{Moller89}
J.~M\o~ller.
\newblock Random tessellations in {${\bf R}^d$}.
\newblock {\em Adv. in Appl. Probab.}, 21(1):37--73, 1989.

\bibitem{Niyogi08}
Partha Niyogi, Stephen Smale, and Shmuel Weinberger.
\newblock Finding the homology of submanifolds with high confidence from random
  samples.
\newblock {\em Discrete Comput. Geom.}, 39(1-3):419--441, 2008.

\bibitem{Puchkin19}
Nikita {Puchkin} and Vladimir {Spokoiny}.
\newblock {Structure-adaptive manifold estimation}.
\newblock {\em arXiv e-prints}, page arXiv:1906.05014, June 2019.

\bibitem{Rineau08}
Laurent Rineau and Mariette Yvinec.
\newblock Meshing 3d domains bounded by piecewise smooth surfaces*.
\newblock In Michael~L. Brewer and David Marcum, editors, {\em Proceedings of
  the 16th International Meshing Roundtable}, pages 443--460, Berlin,
  Heidelberg, 2008. Springer Berlin Heidelberg.

\bibitem{Casal07}
Alberto Rodr\'{\i}guez~Casal.
\newblock Set estimation under convexity type assumptions.
\newblock {\em Ann. Inst. H. Poincar\'{e} Probab. Statist.}, 43(6):763--774,
  2007.

\bibitem{Sharma07}
Alok Sharma and Kuldip~K. Paliwal.
\newblock Fast principal component analysis using fixed-point algorithm.
\newblock {\em Pattern Recognition Letters}, 28(10):1151--1155, 2007.

\bibitem{Sheehy15}
Donald~R Sheehy.
\newblock An output-sensitive algorithm for computing weighted
  $\alpha$-complexes.
\newblock In {\em CCCG}, 2015.

\bibitem{Wasserman18}
Larry Wasserman.
\newblock Topological data analysis.
\newblock {\em Annu. Rev. Stat. Appl.}, 5:501--535, 2018.

\bibitem{Yu97}
Bin Yu.
\newblock Assouad, {F}ano, and {L}e {C}am.
\newblock In {\em Festschrift for {L}ucien {L}e {C}am}, pages 423--435.
  Springer, New York, 1997.

\end{thebibliography}
